\documentclass[12pt,a4paper]{article}
\usepackage{amsfonts}
\usepackage{amsmath}
\usepackage{amssymb}
\usepackage{hyperref}

\usepackage{color}

\usepackage{graphicx}%
\usepackage{ifthen}
\newtheorem{theorem}{Theorem}

\newtheorem{condition}[theorem]{Condition}

\newtheorem{corollary}[theorem]{Corollary}

\newtheorem{definition}[theorem]{Definition}
\newtheorem{example}[theorem]{Example}

\newtheorem{lemma}[theorem]{Lemma}

\newtheorem{proposition}[theorem]{Proposition}
\newtheorem{remark}[theorem]{Remark}

\newenvironment{proof}[1][Proof]{\noindent\textbf{#1.} }{\ \rule{0.5em}{0.5em} \\}

\newcommand{\sesq}[2]{< #1 \, , \, #2>} %
\newcommand{\E}[1]{E_{#1}} %
\newcommand{\dualE}[1]{E^{'}_{#1}} %
\newcommand{\sequiv}[0]{\mathcal{S}}  %
\newcommand{\timeh}[0]{\bar{T}}   %
\newcommand{\numrand}[0]{\bar{m}} %
\newcommand{\wienerp}[2]{W^{#1}_{#2}} %
\newcommand{\wienerq}[2]{\tilde{W}^{#1}_{#2}} %
\newcommand{\zcpx}[1]{\hat{p}_{#1}}   %
\newcommand{\zcpxd}[1]{p_{#1}} %
\newcommand{\prtfpxd}[1]{V_{#1}}   %
\newcommand{\vol}[2]{\sigma^{#1}_{#2}} %
\newcommand{\dualvol}[2]{\sigma^{'{#1}}_{#2}} %
\newcommand{\drift}[1]{m_{#1}} %
\newcommand{\ltrans}[1]{\mathcal{L}_{#1}}   %
\newcommand{\prtfs}[0]{\mathsf{P}} %
\newcommand{\sfprtfs}[0]{\mathsf{P}_{sf}} %
\newcommand{\gammapx}[2]{\gamma^{#1}_{#2}} %
\newcommand{\fwrate}[1]{R_{#1}} %
\newcommand{\spotrate}[1]{r_{#1}} %
\newcommand{\derprod}[2]{\mathsf{D}^{#1}_{#2}} %
\newcommand{\mder}[1]{D_{#1}} %
\newcommand{\norm}[1]{\|#1\|_{\ell^{2}}} %
\newcommand{\multsp}[0]{E^{s}} %
\newcommand{\multspd}[0]{E^{s+1}} %
\newcommand{\ltransm}[1]{\mathcal{L}_{#1}}   %

\date{May 2005 Version 2007.04.03}
\author{Ivar Ekeland\footnote{ekeland@math.ubc.ca; Canada Research Chair in Mathematical Economics,
University of  British Columbia,
Department of Mathematics,
1984 Mathematics Road,
V6T 1Z2 Canada}
\and Erik Taflin\footnote{taflin@eisti.fr; Chair in Mathematical Finance, EISTI,
Ecole International des Sciences du Traitement de l'Information,
Avenue du Parc, 95011 Cergy, France}}
\title{Optimal Bond Portfolios}

\begin{document}
\maketitle

\begin{abstract}
\noindent
We aim to construct a general framework for portfolio
management in continuous time, encompassing both stocks and bonds.
In these lecture notes we give an overview of the state of the art of
optimal bond portfolios and we re-visit main results and mathematical constructions
introduced in our previous publications (Ann. Appl. Probab. \textbf{15}, 1260--1305 (2005)
and Fin. Stoch.   {\bf9}, 429--452 (2005)).

A solution of the optimal bond portfolio problem is given
for general utility functions and volatility operator processes, provided that the 
market price of risk process
has certain Malliavin differentiability properties or is finite dimensional.

The text is essentially self-contained.
\end{abstract}

\noindent {\bf Keywords:} Bond portfolios, optimal portfolios, utility optimization, roll-overs,
          Hilbert space valued processes \\
\noindent {\bf JEL Classification:}  C61, C62, G10, G11 \\
\noindent {\bf Mathematical Subject Classification:}  91B28, 49J55, 60H07, 90C46

\section{Motivation}

The literature on portfolio management starts with the Markowitz portfolio
and  the CAPM 
(\cite{Lintner65}, \cite{Markowitz52}, \cite{Sharp64}). It is a
one-period  %
model, where the information on assets is minimal. Every asset is
characterized by two numbers, its expected return and its covariance with
respect to the market portfolio. With such poor information, one cannot hope
to distinguish between stocks and bonds, and indeed part of the beauty of the
CAPM lies in its generality: it applies to any type of financial assets.

On the other hand, as soon as one tries to make use of all the information
available on assets, important differences appear between stocks and bonds.
Bonds mature, that is they are eventually converted into cash, whereas stocks
do not. The price of bonds depends on interest rates, and the price of stocks,
at least in the academic literature, does not. The bond market is notoriously
incomplete, much more so than the stock market, as is observed in practice. %
As a result, the classical
results on portfolio management, such as Merton's (\cite{Mert69}, \cite{Mert71}),
 concern stock portfolios.
This paper and the papers \cite{I.E.-E.T bond th} and \cite{E.T Bond  Completeness} were
born from a desire to extend them to bond portfolios.

More generally, we aim to construct a general framework for portfolio
management in continuous time, encompassing both stocks and bonds.

The first difficulty to overcome (and, in our opinion, the main  financial one) is the
fact that such a theory should encompass two very different kinds of financial
assets: bonds, which have a finite life, and stocks, which are permanent. We
do it by introducing a new type of financial asset, the \emph{rollovers}. A
rollover of time to maturity $x$ is a bank deposit and which can be cashed at
any time, with accrued interest, provided notice be given time $x$ in advance.
Roll-overs have constant time to maturity (as opposed to zero-coupon bonds,
for instance), and are similar to stocks, in the sense that their main
characteristics do not change with time. By decomposing bonds into rollovers,
instead of decomposing them into zero-coupons, we can hope to incorporate
bonds and stocks into a unified theory of portfolio management.
Rollovers were considered in \cite{Rutkowski99} under the name
``rolling-horizon bond''.

This implies that the time to maturity $x$, rather than the maturity date $T,$
becomes the relevant characteristic of bonds.
Thus, we shall describe bonds using a moving maturity-time frame, where at time $t,$
the origin is the time to maturity $x=0,$ corresponding to the maturity date $T=t.$
 As we shall see very soon, there
will be a mathematical price to pay for that.

At any time $t,$ denote by $p_{t}(  x)  $ the price of a unit
zero-coupon with time to maturity $x.$ The function $x\mapsto p_{t}(
x)  $ will be called the zero-coupon  (price) curve at time $t$; note that the
actual time when that zero-coupon matures is $T=t+x$, and that $T$, is fixed
while $x$ changes with $t$. The zero-coupon curve $p_{t}$ will be understood
to move randomly, and the second difficulty we face is to describe its motion
in some reasonable way. One solution is to decide that $p_{t}$ belongs to a
fixed family of curves, depending on finitely many parameters, so that
\[
p_{t}(  x)  =f(  t,x;r_{1},...,r_{d})
\]
and the random motion of $p_{t}$ is the image of a random motion of the
$r_{i}$, which could, as in spot-rates models, be modelled, for instance by diffusions. This is the
\emph{parametric} approach, which exhibits the classical difficulty of all
parametric approaches, namely that there is no theoretical reason why the
$p_{t}$ should be written in that way, so that the choice of the function $f$
has to be dictated by observational fit. One then has to strike the right
balance between two evils:\ if the number of parameters is too small, the
model will be unrealistic, and if it his higher, it becomes very difficult to calibrate.

We will operate in a \emph{non-parametric} framework: we will make no
assumption on $p_{t}$, beyond some very rough ones, regarding smoothness and
behavior at infinity, nothing that would much constrain their shape.
Mathematically speaking, we will let the curve $p_{t}$ move freely in a linear
space $E$, which will typically be an infinite-dimensional Banach space, of
functions from $[0, \infty [ \,$ to $\mathbb{R}.$

In order to reflect adequately known financial facts, the correct definition
of $E$ must incorporate some basic constraints:

\begin{enumerate}
\item At any time $t$, the zero-coupon prices $p_{t}(  x)  \ $must
depend continuously on the time to maturity $x$. In order for forward  %
interest rates to be well-define\.{d}, they must also have some degree of
differentiability with respect to $x$. So $E$ must consist of continuous
curves with some degree of differentiability.

\item The degree of differentiability of functions in $E$ will determine which  basic
interest rates derivatives can be modelled. If $p_{t}$ is continuous, for
instance, then we can %
introduce bonds. The price of a unit zero-coupon bond with
time to maturity $x$ is $p_{t}(  x)  $; the bond itself, i.e. the value of a portfolio including exactly one bond is
represented by the linear form $p_{t}\mapsto p_{t}(  x).$
Mathematically speaking, this is just the Dirac mass $\delta_{x}$ at $x.$
Now other derivatives such as Call's and Put's on zero-coupon bonds can be introduced,
since the pay-off for each of them is a continuous function of the zero-coupon bond
price $\zcpxd{T}(x),$ with a given time to maturity $x.$
If $p_{t}$ is continuously differentiable, then the forward interest rate
with time to maturity $x,$
$-\frac{\partial}{\partial x}p_{t}(x)/p_{t}(x)$ is well-defined, and further
contingent claims can be defined, such as caps, floors and swaps.
\item The curve $p_{t}$ %
will be understood to move randomly in
$E$, the randomness being driven by a Brownian motion. We will therefore need
to define Brownian motions in the infinite-dimensional space $E$, which for
all practical purposes will require $E$ to be a Hilbert space.
\item The accepted standard in mathematical modelling of zero-coupon
prices (the Heath-Jarrow-Morton model, henceforth HJM) is to decide that the
real-valued process %
$t \mapsto p_{t}(  T-t),  $ the price at time
$t$ of a unit zero-coupon maturing at a given time $T$, is an It\^o
process satisfying an stochastic ordinary  differential equation (SODE).
As is well-known, for fixed $x,$ the real-valued process %
$t\mapsto p_{t}(  x) , $  which is also an It\^o process,
then no longer satisfies an SODE. Indeed, if  $f(t,T) \equiv p_{t}(  T-t),$
then we have $p_{t}(x)= f(t,t+x)$  %
so that for fixed $x:$
\begin{equation}
d_{t}  p_{t}(  x)  =[d_{t}  f(t, T)+\frac{\partial f(t,T)}{\partial T}  dT]_{T=t+x}
=[d_{t}  f(t, T)]_{T=t+x}+\frac{\partial p_{t}(  x)}{\partial x} dt.
\label{1}%
\end{equation}
Here the right-hand side (r.h.s) depends, not only on $p_{t}(  x)  ,$ but also on
its partial derivative with respect to $x$. So,  equation (\ref{1}) for $p$ is a SPDE,
stochastic partial  differential equation, where the first term on the r.h.s depends
only on the un-known $p_{t}(  x) ,$ since $f(\cdot,T)$ satisfies an SODE. This is the well-known difficulty
of the Musiela parametrization (see \cite{Musi93}), and the space $\E{}$ shall permit a simple mathematical
formulation of the SPDE (\ref{1}).
\item
At any time $t$, the zero-coupon prices $p_{t}(  x) $ should
go to zero as  the time to maturity $x$ goes to $ \infty.$ To include also the
trivial case, where all interest rates vanish, and also cases where the forward
rates converges  rapidly to zero as $x \rightarrow  \infty,$ we only require
that $\lim p_{t}(x) $ exists as $x \rightarrow  \infty.$ N.B. We will chose $\E{}$
such that the elements $f \in E$ satisfying
$\lim_{x \rightarrow \infty} f(x)=0$ form a closed
sub-space of $\E{},$ in order to cover easily the case where $p_{t}(  x) \rightarrow 0. $
\end{enumerate}
Formula (\ref{1}) is really an infinite family of coupled
 equations, one for each
$x \geq  0,$ describing the motion of the random variable $p_{t}(x),$
which we write
\begin{equation}
dp_{t}(x)=p_{t}(x)m_{t}(x)dt+p_{t}(x)\sigma_{t}(x)dW_{t}+\frac{\partial p_{t}(  x)}{\partial x}  dt,
\label{2'}%
\end{equation}
where for the moment $W$ is thought of as being a high dimensional  Brownian motion.
Let us rewrite it as a single stochastic evolution equation for the motion of the random curve
$p_{t}$ in $E,$ i.e. as a SODE in $E:$
\begin{equation}
dp_{t}=p_{t}m_{t}dt+p_{t}\sigma_{t} dW_{t}+(  \partial p_{t})  dt
\label{2}%
\end{equation}
where
$\partial$ is the differentiation operator with respect to \textit{time to maturity},
i.e. it is defined by $(\partial u)(x)=\frac{d u(  x)}{d x},$
for differentiable $u \in \E{}.$
Since the left-hand side ``belongs'' to $E,$ so must the right-hand side, and then
$\partial p_{t}$ must belong to $E$.
There are ways to achieve that. One
is to choose a framework where the operator $\partial $ is
continuous over all of $E$. Then so is its $n$-th iterate $\partial ^{n}$,
so that the space $E$ must consist of functions which have infinitely many
derivatives. Unfortunately, the natural topology of such spaces cannot be
defined by a single norm, except for very particular cases,
 and the mathematics become more demanding. %
A second more standard way  to proceed is to consider $\partial $ as an unbounded
operator in a Hilbert space $E$, so that $\partial $ is defined only on a
subspace $\mathcal{D}(  \partial )  \subset E$, called the domain of the
operator. One would then hope to define the solution of equation (\ref{2}) in
such a way that, if the initial condition $p_{0}$ lies in $\mathcal{D}(
\partial )  ,$ then $p_{t}$ remains in $\mathcal{D}(
\partial )  $ for every $t,$ so that $t \mapsto p_{t}$ is a trajectory in
$\mathcal{D}(\partial).$ In other words, if
$p_{0}(  x)  \ $is differentiable with respect to $x$, so should
the functions $x \mapsto p_{t}(x)$ be for all $t>0.$

To summarize, the introduction of rollovers and a moving frame forces us to complicate the
equations for price dynamics, by incorporating an additional term,
$\partial p_{t}$. To be able to solve the relevant equations, we have to
treat $\partial $ as an unbounded operator in Hilbert space. The definition
of the relevant Hilbert space has to incorporate basic properties which we
expect of zero-coupon curves.

This suits our purpose well, for it enables us to work in a non-parametric
framework, where no particular shape is assigned to the the zero-coupon
curves. On the other hand, we then have to use the theory of Brownian motion
in infinite-dimensional Hilbert spaces and the corresponding stochastic
integrals, which creates some additional difficulties. We do not limit the
number of sources of noise, indeed in our paper there can be infinitely many.
This is natural, since the already mentioned experimental fact,
that even using a large number of bonds, not all interest rate derivatives can be hedged.
The third difficulty to overcome, is the mathematically significant fact that
such a market can not be complete in the usual sense, i.e. every (sufficiently integrable)
contingent claim being hedgeable. This has important implications for the solution
of the portfolio optimization problem. The now classical two-step solution, so successfully
applied to the case of a
finite number of stocks (cf.  \cite{Kr-Scha}, \cite{Pliska86}), consisting of first determining the optimal final wealth
by  duality methods and then determining a hedging portfolio, does not (yet at least) apply
to the general infinite dimensional bond markets.
In this paper (see \cite{I.E.-E.T bond th} and \cite{E.T Bond  Completeness})
 we give, within the considered general  It\^o process model, the optimal final wealth
for every case it exists (Proposition \ref{exist unique X}).
The existence of an optimal portfolio, is then established by the construction
of a hedging portfolio for two cases :
The first
is for deterministic $\E{}$-valued drift $m$ and volatility operator $\sigma,$
where we give a necessary and sufficient condition for the existence of an optimal portfolio. Here
there can exist several equivalent martingale measures (e.m.m.), so the market
can clearly be incomplete in every sense of the word.
The second is for certain stochastic
$m$ and $\sigma,$ for which there is a unique market price of risk process $\gamma.$
There is then a unique e.m.m. $Q.$ Now, certain integrability conditions on the $\ell^{2}$-valued
Malliavin derivative of the Radon-Nikodym density $dQ/dP$ leads to the construction of a
hedging portfolio.

We have tempted to make these notes self-contained, with exception of the general hedging result
in Theorem \ref{th completeness l^2}. The notes first recall some basic facts concerning
linear operators and semi-groups in Hilbert spaces, Sobolev spaces and stochastic
integration in Hilbert spaces. The theory of bond portfolios and hedging of
interest rate derivatives are then introduced.
Once this theory is explained, the paper proceeds to a short solution of the
optimization problem, leading to the results %
of \cite{I.E.-E.T bond th} and \cite{E.T Bond  Completeness}.
In particular, under the assumption that the market prices
of risk are deterministic, some explicit formulas are given,
very similar in spirit to those who are known in the case of stock portfolios,
and a mutual fund theorem is formulated.
We conclude by stating an alternative formulation, of the optimization problem,
within a Hamilton-Jacobi-Bellman approach.

\section{Mathematical preliminaries}
\subsection{Hilbert spaces and bounded maps}
We shall be working with separable infinite-dimensional real Hilbert spaces.
Let  $E$ be a Hilbert space with scalar product $(\;,\;)_{E}$
and norm $\| \;\;\|_{E},$ simply denoted $(\;,\;)$ and $\| \;\;\|$  if no risk for confusion.
The topology and convergence in $E$ is w.r.t. this norm, if not otherwise stated, i.e.
the strong topology and convergence.  By definition $E$ is, \emph{separable} if it has
a countable dense subset. One shows easily that $E$ is separable iff it has a
countable orthonormal  basis $e_{n},n\in \mathbb{N},$ i.e.  $(  e_{i},e_{j})
=0$ for $i\neq j$ and $\Vert e_{i} \Vert =1,$ so that every $x \in E$
can be written:%
\[
x=\sum_{n=0}^{\infty}\left(  x,e_{n}\right)  e_{n},%
\]
where the right-hand side converges in $E.$ Since the $e_{n}$ are orthonormal,
we have Parseval's equality:%
\[
\left\Vert x\right\Vert^{2} =\sum_{n=0}^{\infty}\left\vert \left(  x,e_{n}\right)
\right\vert ^{2}.
\]
A typical separable Hilbert space is $\ell^{2}$, which is the space of all
real sequences $a_{n},n\in \mathbb{N},$ such that $\sum\left\vert a_{n}\right\vert
^{2}<\infty$. The scalar product in $\ell^{2}$ is given by $\left(
a,b\right)  =\sum a_{n}b_{n}$. In fact, every infinite dimensional separable Hilbert space $E$ is
isomorphic to $\ell^{2}$. The map
\begin{equation}
x\mapsto a_{n}=(  x,e_{n})_{E}  ,n\in \mathbb{N},
\label{3}%
\end{equation}
of $E$ into $\ell^{2}$ is a linear bijection and it preserves norms on both sides.

A linear map $ L:E_{1}\rightarrow E_{2}$ is continuous if and only if it is
\emph{bounded}, that is if there exists a constant $c$ such that $\left\Vert
Lx\right\Vert _{E_{2}}\leq c\left\Vert x\right\Vert _{E_{1}}$ for every $x\in E_{1}.$
The (operator) norm of $L$ is then defined to be the infinimum of all such $c$:%
\[
\left\Vert L\right\Vert =\inf\left\{  c\ |\left\Vert Lx\right\Vert _{E_{2}}\leq
c\left\Vert x\right\Vert _{E_{1}}\ \forall x\right\}.
\]
For example, the linear map in (\ref{3}) of $E$ onto $\ell^{2}$ as well as
its inverse has norm $1.$
The linear space  of all continuous linear maps from $E_{1}$
to $E_{2},$ $L(E_{1},E_{2}),$ is a Banach space when given this norm. One writes $L(E)$
as a shorthand for $L(E,E).$
Linear maps
are also called linear operators or just operators.
A \emph{bounded operator} on $E$ is a bounded linear map from $E$ into itself.
The \emph{dual} space  $E'$ of $E,$ i.e. the space of all linear continuous functionals on $E,$
is given by $E'=L(E,\mathbb{R}).$ By the F. Riesz representation theorem,
\begin{equation} \label{riesz}
F \in E' \; \text{iff} \; \exists f \in E \; \text{such that} \; F(x)=(f,x) \; \forall x \in E.
\end{equation}
Also $\|F\|_{E'}=\|f\|_{E},$ so $E'$ and $E$ are isomorphic. In this paper we will
often use, in the context of Sobolev spaces, other representations of the dual $E'.$

By duality, every operator in  $L(E_{1},E_{2})$ corresponds to an operator in
$L(E_{2}',E_{1}').$ Using the representation of the dual space given by (\ref{riesz}), the adjoint
operator $A^{*}$ of $A \in L(E_{1},E_{2})$ is defined by $A^{*}y=y^{*},$ where for $y \in E_{2}$
the element $y^{*} \in E_{1}$ is defined by
\begin{equation} \label{adj}
(y^{*},x)_{E_{1}} =(y,Ax)_{E_{2}}\; \forall x \in E_{1}.
\end{equation}
This defines an operator $A^{*} \in L(E_{2},E_{1})$. One easily checks that $(A^{*})^{*}=A$ and
$\left\Vert A^{*}  \right\Vert=\left\Vert A \right\Vert.$
Let us consider a simple example, which will be relevant in the sequel
of this paper:
\begin{example}[Left-translation in $L^{2}$] \label{translation} \text{} \\ \normalfont
i) Let $E=L^{2}(\mathbb{R})$ and let $a$ be a given real number. Define the operator
$A$ on $E$ by $(Af)(x)=f(x+a).$ Then $\Vert A \Vert=1$ and $(A^{*}f)(x)=f(x-a).$
We note that $A$ has a bounded inverse $A^{-1}$ given by $(A^{-1}f)(x)=f(x-a),$
so $A A^{*}=A^{*} A =I,$ where $I$ is the identity operator. \\
ii) Let $E=L^{2}([0,\infty[)$ and let $a>0$ be a given real number. Define the operator
$A$ on $E$ by $(Af)(x)=f(x+a).$ Here we find that $\Vert A \Vert=1,$ that a.e. $(A^{*}f)(x)=0$ if $0\leq x < a$
and that $(A^{*}f)(x)=f(x-a)$ if $a \leq x.$ In this case $A^{*}$ is one-to-one and
$A A^{*}=I.$ But $ A^{*} A$ is the orthogonal projection on the (non-trivial) closed subspace of
$E$ of functions with support in $[a,\infty[ \, .$ So $A^{*}A \neq I.$
\end{example}
An operator $S \in L(E_{1},E_{2})$ is called unitary if $S S^{*}=S^{*} S =I.$
This is the case of $A$ in (i) of Example \ref{translation}.
An operator $S \in L(E_{1},E_{2})$ is called isometric if $S^{*} S =I.$
This is the case of $A^{*}$ in (ii) of Example \ref{translation}.

We will be interested in a particular class of bounded operators on $E.$
 We begin with an easy result

\begin{lemma}
Suppose $L \in L(E_{1},E_{2})$ and that  we have:
\[
\sum_{n=0}^{\infty}\left\Vert Le_{n}\right\Vert ^{2}<\infty
\]
for an orthonormal basis $e_{n},n\in \mathbb{N}$ in $E_{1}.$ Let $f_{n},n\in \mathbb{N}$ be another
orthonormal basis. Then:%
\[
\sum_{n=0}^{\infty}\left\Vert Le_{n}\right\Vert ^{2}=\sum_{n=0}^{\infty
}\left\Vert Lf_{n}\right\Vert ^{2}%
\]
\end{lemma}

\begin{definition} \label{H-S}
An operator $L$ on $E_{1}$ into $E_{2}$ is Hilbert-Schmidt if $\ \sum_{n=0}^{\infty
}\left\Vert Le_{n}\right\Vert ^{2}<\infty$ for some orthonormal basis
$e_{n},n\in \mathbb{N}$, in $E_{1}.$ Its Hilbert-Schmidt norm is defined to be:
\[
\left\Vert L\right\Vert _{\mathcal{HS}}=\left(  \sum_{n=0}^{\infty}\left\Vert
Le_{n}\right\Vert ^{2}\right)  ^{1/2}.%
\]
It does not depend on the choice of the orthonormal basis $e_{n},n\in \mathbb{N}$, in
$E$. The linear space  of Hilbert-Schmidt operators from $E_{1}$ into $E_{2}$ is denoted
$\mathcal{HS}(E_{1},E_{2}).$
\end{definition}

Hilbert-Schmidt operators are bounded (in fact, $\left\Vert L\right\Vert
\leq\left\Vert L\right\Vert _{\mathcal{HS}}$) and even compact:\ they map
bounded subsets of $E_{1}$ into relatively compact subsets of $E_{2}.$ In other words,
if $L$ is Hilbert-Schmidt and $\left(  x_{n}\right)  _{n\in \mathbb{N}}$ is a bounded
sequence, then one can extract from $\left(  Lx_{n}\right)  _{n\in \mathbb{N}}$ a
norm-convergent subsequence.
This property of a Hilbert-Schmidt operator $L$
follows from the fact that $L$ is the limit in the operator norm of finite rank operators.
The space $\mathcal{HS}(E_{1},E_{2})$ endowed with the Hilbert-Schmidt norm defines a
Hilbert space.

Some general references for this subsection are: \cite{Kato66}, \cite{Lax 02}, \cite{Rudin 1}, \cite{Rudin 2}.
\subsection{Linear semi-groups and unbounded operators.}
Let $L$ be a bounded linear operator on $E$. For every $t\in \mathbb{R}$, define:%
\[
\Phi\left(  t\right)  =e^{tL}=\sum_{i=0}^{\infty}\frac{1}{n!}t^{n}L^{n},
\]
which converges in the operator norm.
Then $\Phi\left(  t\right)  $ is a bounded linear operator for every $t$, and
we have the relation:%
\begin{equation}
\Phi\left(  t+s\right)  =\Phi\left(  t\right)  \Phi\left(  s\right)
\; \forall s,t \in \mathbb{R} \; \; \text{and} \;\; \Phi(  0)  =I,
\label{43}%
\end{equation}
where $I$ is the identity operator on $E,$
from which it follows that $\Phi\left(  t\right)  $ and $\Phi\left(  s\right)
$ commute %
and that $\Phi\left(  t\right)  $
is invertible for every $t$. Relation (\ref{43}) states that the map
$t\mapsto \Phi\left(  t\right)  $ is a group homomorphism. Note that it is
continuous in the norm topology for operators:%
\begin{equation}
\left\Vert \Phi\left(  t\right)  -I\right\Vert \rightarrow0 \;\text{ when}\; t\rightarrow0.
\label{43'}
\end{equation}
The solution of the Cauchy problem:%
\begin{align}
\frac{dx(t)}{dt}  &  =Lx(t),\label{31}\\
x\left(  0\right)   &  =x_{0} \label{32}%
\end{align}
is given by $x\left(  t\right)  =\Phi\left(  t\right)  x\left(  0\right)  $.
In other words, $\Phi\left(  t\right)  $ is the flow associated with the
ordinary differential equation (\ref{31}). We can recover $L$ from
$\Phi\left(  t\right)  $ by writing:%
\begin{equation}
Lx=\lim_{h\rightarrow0}\frac{1}{h}\left[  \Phi\left(  h\right)  x-x\right], \:\: x \in E.
\label{34}%
\end{equation}
The norm continuity of the mapping $t \mapsto \Phi(  t)$
is exceptional and has to be replaced by a more useful weaker property
(cf. Definition 1, Sect. 1, Chap. IX of \cite{Yosida}):
\begin{definition}
A family $\Phi\left(  t\right)  ,t\geq0$, of bounded operators  on $E$ is
called a one parameter semi-group if $\Phi\left(  0\right)  =I$, and for all $t\geq0$
and $s\geq0$ we have:%
\begin{equation}
\Phi\left(  t+s\right)  =\Phi\left(  t\right)  \Phi\left(  s\right)
=\Phi\left(  s\right)  \Phi\left(  t\right).  \label{33}%
\end{equation}
It is said to be strongly continuous or to be of class $(C_{0})$ if, for every $x\in E$, we have:%
\begin{equation}
\lim_{t\rightarrow0}\Phi\left(  t\right)  x=x. \label{35}
\end{equation}
It is said to be a contraction semi-group if $\|\Phi(t) \| \leq 1$
for all $t \geq 0.$
\end{definition}
Note that, since equality (\ref{33}) is supposed to hold only for positive $s$ and $t$,
the operators $\Phi\left(  t\right)  $ are no longer necessarily invertible, as in the case
of a group. It can be proved easily that, if the semi-group $\Phi\left(  t\right)
$ is strongly continuous, then $\lim_{s\rightarrow t}\Phi\left(  s\right)
x=\Phi\left(  t\right)  x$ and there are constants $c$ and $C$ such that
$\left\Vert \Phi\left(  t\right)  \right\Vert \leq C\exp\left(  ct\right) .$
We also note that if  $[0, \infty [ \; \ni t \mapsto \Phi(t)$ is a
one parameter semi-group, so is the family of adjoint operators $[0, \infty [ \; \ni t \mapsto \Phi^{*}(t),$
where we define $\Phi^{*}(t)=(\Phi(t))^{*}.$
\begin{example} \label{translation 1} \text{} \\ \normalfont
In the situation of (i) (resp. of (ii)) of Example \ref{translation}, for given $a,$
let $\Phi_{1} (a)=A$ (resp.  $\Phi_{2} (a)=A$).
Then  $\mathbb{R} \ni t \mapsto \Phi_{1} (t)$ is a strongly continuous contraction group.
However $[0, \infty [ \; \ni t \mapsto \Phi_{2}(t)$ is only a strongly continuous contraction
semi group, which can not be extended to a group. In fact, $\Phi_{2}(t)$ is not
invertible for $t>0.$
\end{example}

We now try to extend formula (\ref{34}). It turns out that when $\Phi$ is no
longer norm-continuous, but only strongly continuous, the
right-hand side does not converge for every $x$, and if the limit exists, it
does not depend continuously on $x.$
The set of $x$ for which the limit exists
is obviously a linear subspace of $E$ and on this subspace the limit is a linear
function, let's say $G$ of $x.$ More formally, let $\mathcal{D}(G)$ be the subset of $E$
\textit{of all elements} $x \in E$ for which the strong limit
\begin{equation}
Gx=\lim_{h\rightarrow0}\frac{1}{h}\left[  \Phi\left(  h\right)  x-x\right]
\label{40}%
\end{equation}
exists.
\begin{theorem} \label{sem-group gen}
Assume
$\Phi  $ is a strongly continuous semi-group. The set $\mathcal{D}(G)$
is then a  dense linear subspace of $E$ and $G$ given by (\ref{40}) defines a linear map
$G: \mathcal{D}(G)  \rightarrow E.$ This map is closed, i.e. if
 $x_{n}$ is a sequence in $\mathcal{D}(G)$ such that $x_{n} \rightarrow \bar{x} \in E$
and $Gx_{n}\rightarrow\bar{y} \in E$ then $\bar{x} \in \mathcal{D}(G)$ and $\bar{y}=G\bar{x}.$

For every $x\in \mathcal{D}(G)$ and $t\geq 0$ we have
 $\Phi\left(  t\right)  x \in \mathcal{D}(G),$
\begin{equation}
G\Phi\left(  t\right)  x    =\Phi\left(  t\right)  Gx\label{41}
\end{equation}
and
\begin{equation}
\frac{d}{dt}\Phi\left(  t\right)  x    =G\Phi\left(  t\right)  x. \label{42}%
\end{equation}
\end{theorem}
\begin{proof}
By definition $\mathcal{D}(G)  $ is the set of $x$ where the limit in
formula (\ref{40}) exists (note that this is a strong limit, meaning that we
should have norm-convergence), and $Gx$ then is the value of that limit.
Clearly $G:\mathcal{D}(G)  \rightarrow E$ is a linear map.

Given any $x\in E$ and $t>0$, consider the integral:%
\[
X\left(  t\right) =\int_{0}^{t}\Phi\left(  s\right)  xds.
\]
It is well-defined since the integrand is a continuous function from $\left[
0,t\right]  $ into $E$. Using the semi-group property, we have:%
\begin{align*}
\frac{1}{h}\left[  \Phi\left(  h\right)  X\left(  t\right)  -X\left(
t\right)  \right]   &  =\frac{1}{h}\left[  \Phi\left(  h\right)  \int_{0}%
^{t}\Phi\left(  s\right)  xds-\int_{0}^{t}\Phi\left(  s\right)  xds\right] \\
&  =\frac{1}{h}\left[  \int_{0}^{t}\Phi\left(  s+h\right)  xds-\int_{0}%
^{t}\Phi\left(  s\right)  xds\right] \\
&  =\frac{1}{h}\int_{0}^{h}\Phi\left(  s+h\right)  xds-\frac{1}{h}\int
_{0}^{h}\Phi\left(  s\right)  xds\\
&  \rightarrow\Phi\left(  t\right)  x-x.
\end{align*}
This proves that $X\left(  t\right)  $ belongs to $\mathcal{D}(G) $. Then
so does $\frac{1}{t}X\left(  t\right)  $, and when $t\rightarrow0$, we have
$\frac{1}{t}X\left(  t\right)  \rightarrow x$, so $\mathcal{D}(G)  $ is dense in $H$,
as announced.

Now write:%
\[
\frac{1}{h}\left[  \Phi\left(  t+h\right)  -\Phi\left(  t\right)  \right]
x=\Phi\left(  t\right)  \frac{\Phi\left(  h\right)  -I}{h}x=\frac{\Phi\left(
h\right)  -I}{h}\Phi\left(  t\right)  x.
\]
If $x\in \mathcal{D}(G)  $, the second term converges to $\Phi\left(
t\right)  Gx$ and the third one to $G\Phi\left(  t\right)  x$.
Formulas (\ref{41}) and (\ref{42}) now follow, since these two
terms must be equal.

To prove the last condition, note that:%
\begin{equation}
\forall x\in \mathcal{D}(G)  ,\ \ \ \Phi\left(  t\right)  x-x=\int_{0}%
^{t}\Phi\left(  s\right)  Gxds. \label{44}%
\end{equation}
Indeed, we have two functions of $t$, with values in $E\,$, which are zero for
$t=0$ and which have the same derivative, namely $\Phi\left(  t\right)  Gx$,
for every $t>0$. So they must be equal. Now take a sequence $x_{n}%
\rightarrow\bar{x}$, and assume that $Gx_{n}=y_{n}\rightarrow\bar{y}$ in $E$.
Writing $x=x_{n}$ in formula (\ref{44}), we get:%
\[
\Phi\left(  t\right)  \bar{x}-\bar{x}=\int_{0}^{t}\Phi\left(  s\right)
\bar{y}ds.
\]
Dividing by $t$ and letting $t\rightarrow0$, we find that $\bar{x}\in \mathcal{D}(G)$
 and that $\bar{y}=G\bar{x}$.
%\AAAqed
\end{proof}

\begin{definition}
In the situation of Theorem \ref{sem-group gen},
$G$ is called the \emph{infinitesimal generator} of the semi-group $\Phi.$
\end{definition}
A linear map $L:\mathcal{D}(  L)  \rightarrow E_{2}$, where $\mathcal{D}\left(  L\right)  $
is a  subspace of $E_{1},$ is called an operator from $E_{1}$
to $E_{2}$ with domain $\mathcal{D}\left(  L\right).$ That two operators are equal, $L_{1}=L_{2},$
means that they have the same domain $\mathcal{D}( L_{1}) =\mathcal{D}( L_{2})$
and that $L_{1}x=L_{2}x$ for all $x$ in the domain.
 The operator $L$ is densely defined
if $\mathcal{D}(  L)$ is dense in $E_{1}.$
It is called a bounded operator if there exists a finite constant
 $C \geq 0$ such that for all $x \in \mathcal{D}(  L)$ one has
$ \|Lx\| \leq C \|x\| $ and it is called an \emph{unbounded operator} if
such $C$ does not exist. It
is \emph{closed} if its graph $\{(x,Lx) \, | \, x \in \mathcal{D}\left(  L\right) \}$
is a closed subset of $E_{1} \times E_{2},$ which extends the definition in the preceding  theorem.
With these definitions, we can
rephrase part of the preceding theorem by saying that every strongly
continuous semi-group in $E$ has a unique infinitesimal generator, which is a densely defined closed operator in $E.$
The problem to determine if a given densely defined closed operator $L$ in $E$ is the
infinitesimal generator of a  strongly continuous semi-group is more difficult  and
we refer the interested reader to the references mentioned in the end of this subsection.

The definition of the \emph{adjoint} of an operator can be extended to unbounded operators.
Let $L$ be a densely defined operator from $E_{1}$ to $E_{2}.$ We introduce the adjoint operator
$L^{*}$ to $L.$ The domain of $\mathcal{D}(L^{*})$ consists of all $y \in E_{2}$ for which the linear
functional
\begin{equation} \label{adj gen 1}
x \mapsto (y,Lx)
\end{equation}
is continuous on $\mathcal{D}(L),$ endowed with the strong topology of $E_{1}.$
For $y \in \mathcal{D}(L^{*})$ we define
$L^{*}y$ by
\begin{equation} \label{adj gen 2}
 (L^{*}y,x)=(y,Lx)  \; \forall x \in \mathcal{D}(  L).
\end{equation}
This defines $L^{*}y$ uniquely, since $\mathcal{D}(  L)$ is dense in $E_{1}.$ One proves
that $\mathcal{D}(L^{*})$ is dense in $E_{2}$ if $L$ is also closed.

 An operator $L$ in $E$ is called \emph{selfadjoint} if $L^{*}=L$ and skew-adjoint if $L^{*}=-L.$
We have the following clear-cut result (Stone's theorem): $L$ is the infinitesimal
generator of a group of unitary operators iff $L$ is skew-adjoint.
\begin{example} \label{translation 2} \text{} \\ \normalfont
In the situation of Example \ref{translation 1}, let $L_{1}$ and $L_{2}$ be the infinitesimal
generators  of $\Phi_{1}$ and $\Phi_{2}$ respectively.  $L_{1}$ is given by
\[
\mathcal{D}( L_{1})=\{f \in L^{2}(\mathbb{R}) \; | \; f' \in L^{2}(\mathbb{R}) \},
\]
and $ (L_{1}f)(x)=f'(x),$ where $f'$ is the derivative of $f.$
 $L_{2}$ is given by
\[
\mathcal{D}( L_{2})=\{f \in L^{2}([0,\infty[ \,) \; | \; f' \in L^{2}([0,\infty[ \,) \},
\]
and $ (L_{2}f)(x)=f'(x).$
Since $\Phi_{1}$ is a group of unitary operators, we have that $L_{1}^{*}=-L_{1}.$
$\Phi_{2}$ is not a semi-group of unitary operators, so $L_{2}^{*} \neq -L_{2}.$
A simple calculation shows that
\[
\mathcal{D}( L_{2}^{*})=\{f \in L^{2}([0,\infty[ \,) \; | \; f(0)=0 \; \text{and} \;  f' \in L^{2}([0,\infty[ \,) \}
\]
and $(L_{2}^{*}f)(x)=-f'(x).$
So here $\mathcal{D}( L_{2}^{*}) \subset \mathcal{D}( L_{2}),$ with strict inclusion.
One checks that $ \Phi_{2}^{*}$ is a strongly continuous semi-group in $L^{2}([0,\infty[ \,).$
It represents  right translations of functions. Its infinitesimal generator is $L_{2}^{*}.$
\end{example}

Some general references for this subsection are:
\cite{Kato66}, \cite{Lax 02}, \cite{Rudin 2}, \cite{Yosida}.
\subsection{Sobolev spaces}
For any integer $n\geq0$, the Sobolev space $H^{n}(  \mathbb{R})  $ is
defined to be the set of functions $f$ which are square-integrable together
with all their derivatives of order up to $n$:
\[
f\in H^{n}(\mathbb{R})  \Longleftrightarrow\int_{-\infty}^{\infty}\left[
f^{2}+\sum_{k=1}^{n}\left(  \frac{d^{k}f}{dx^{k}}\right)^{2}\right]
dx\leq\infty.
\]
This is a linear space, and in fact a Hilbert space with  norm given by:
\[
\|f\|_{H^{n}}=\left(\int_{-\infty}^{\infty}\left[  f^{2}+\sum_{k=1}^{n} (
\frac{d^{k}f}{dx^{k}})  ^{2}\right]  dx \right)^{1/2}.
\]
It is a standard fact that this norm of $f$  can be expressed in terms of
the Fourier transform $\hat{f}$ (appropriately normalized) of $f$ by:
\[
\left\Vert f\right\Vert _{H^{n}}^{2}=\int_{-\infty}^{\infty}\left(
1+y^{2}\right)  ^{n}\left\vert \hat{f}\left(  y\right)  \right\vert ^{2}dy.
\]
The advantage of that new
definition is that it can be extended to non-integral and non-positive values.
For any real number $s$, not necessarily an integer nor positive, we define
the Sobolev space $H^{s}(\mathbb{R})  $ to be the Hilbert space of
functions associated with the following norm:%
\begin{equation} \label{Hn}
\left\Vert f\right\Vert _{H^{s}}^{2}=\int_{-\infty}^{\infty}\left(
1+y^{2}\right)  ^{s}\left\vert \hat{f}\left(  y\right)  \right\vert ^{2}dy.
\end{equation}
Clearly, $H^{0}(\mathbb{R})=L^{2}(\mathbb{R})$ and
$H^{s}(\mathbb{R})\subset H^{s^{\prime }}(\mathbb{R})$ for $s\geq
s^{\prime }$ and in particular
$H^{s}(\mathbb{R})\subset L^{2}(\mathbb{R})\subset H^{-s}(\mathbb{R}),$ for $s \geq 0.$
$H^{s}(\mathbb{R})$ is, for general $s \in \mathbb{R},$ a space of (tempered) distributions.
For example $\delta^{(k)},$ the $k$-th derivative of a delta Dirac distribution, is
in $H^{ -k-1/2-\epsilon}(\mathbb{R})$ for $\epsilon >0.$

In the case when $s>1/2$, there are two classical results.
\begin{theorem}
[Continuity of multiplication] If $s>1/2$, if $f$ and $g$ belong to
$H^{s}(  \mathbb{R})  $, then $fg$ belongs to $H^{s}(  \mathbb{R})  $, and
the map $(  f,g)  \rightarrow fg$ from $H^{s}\times H^{s}$ to
$H^{s}$ is continuous.
\end{theorem}
Denote by $C_{b}^{n}(\mathbb{R})$ the space of $n$ times continuously differentiable
real-valued functions which are bounded together with all their $n$ first derivatives.
Let $C_{b0}^{n}(\mathbb{R})$ the  closed subspace of $C_{b}^{n}(\mathbb{R})$
of functions which  converges to $0$ at $\pm \infty$ together with
all their $n$ first derivatives. These are Banach spaces for the norm:
\[
\left\Vert f\right\Vert _{C_{b}^{n}}=\max_{0\leq k\leq n}\ \sup_{x}%
\ \left\vert f^{\left(  k\right)  }\left(  x\right)  \right\vert =\max_{0\leq
k\leq n}\left\Vert f^{\left(  k\right)  }\right\Vert _{C_{b}^{0}}.
\]
\begin{theorem} [Sobolev embedding]
If $s>n+1/2$ and if $f \in H^{s}(\mathbb{R}),$
then \ there is a function $g\ $in $C_{b0}^{n}(  \mathbb{R})  $ which is
equal to $f$ almost everywhere. In addition, there is a constant $c_{s}$,
depending only on $s$, such that:
\[
\Vert g\Vert _{C_{b}^{n}}\leq c_{s}  \Vert f\Vert _{H^{s}}.
\]
\end{theorem}
From now on we shall no longer distinguish between $f$ and $g$, that is, we
shall always take the continuous representative of any function in
$H^{s}(  \mathbb{R})  $. As a consequence of the Sobolev embedding theorem,
if $s>1/2$, then any function $f$ in $H^{s}(  \mathbb{R})  $ is continuous
and bounded on the real line
and converges to zero at $\pm \infty,$ so that its value is defined everywhere.

We define, for $s\in \mathbb{R},$ a continuous bilinear form on $H^{-s}(\mathbb{R}) \times H^{s}(\mathbb{R})$ by:
\begin{equation}
\sesq{f}{g}=\int_{-\infty}^{\infty} \overline{\left( \hat{f}(y)\right) }\text{ } \hat{g}(y)dy,
\label{sesqprod}
\end{equation}
where $ \overline{z}$ is the complex conjugate of $z.$ Schwarz inequality and (\ref{Hn})
give that
\begin{equation}
|\sesq{f}{g}| \leq \Vert f\Vert _{H^{-s}} \Vert g\Vert _{H^{s}},
\label{sesqprod 1}
\end{equation}
which indeed shows that the bilinear form in (\ref{sesqprod}) is continuous.
We note that formally the bilinear form (\ref{sesqprod}) can be written
\[
\sesq{f}{g}=\int_{-\infty}^{\infty} f(x) g(x)dx,
\]
where, if $s \geq 0,$ $f$ is in a space of distributions $H^{-s}(\mathbb{R})$ and
$g$ is in a space of ``test functions'' $H^{s}(\mathbb{R}).$

Any continuous linear form $g\rightarrow u\left( g\right) $ on $H^{s}(\mathbb{R})$
is, due to (\ref{Hn}), of
the form $u(g) =\sesq{f}{g}$ for some $f\in H^{-s}(\mathbb{R}),$ with
$\Vert f\Vert _{H^{-s}}=\Vert u\Vert _{(H^{s})^{'}}$, so
that henceforth we can identify the dual $( H^{s}(\mathbb{R}))^{'}$
of $H^{s}(\mathbb{R})$ with $H^{-s}(\mathbb{R}).$
In particular, if $s>1/2$ then $H^{s}(\mathbb{R}) \subset C_{b0}^{0}( \mathbb{R}),$
so $H^{-s}(\mathbb{R})$ contains all bounded Radon measures.

In the sequel, we will also be interested in functions defined only on the
half-line $[0,\infty[ \,.$ Let $s \geq 0.$ We define the space $H^{s}([0,\infty [ \,)$
to be the set of restrictions to $[0,\infty [ \,$ of functions in $H^{s}(\mathbb{R}).$
This is clearly a linear space. To turn it into a Hilbert space,
we have to use the following norm:
\begin{equation} \label{H norm}
\Vert f\Vert _{H^{s}(  [0,\infty))  }=\inf \left \{
\Vert g\Vert _{H^{s}(  \mathbb{R})  }\ |\ g(  x)
=f(  x)  \ \text{a.e. on }[0,\infty)\right\}.
\end{equation}
This is a Hilbert space norm on $H^{s}([0,\infty [ \,)  ,$ which is the
natural restriction of the norm on $H^{s}(  \mathbb{R})  $. For instance, if
$f$ is a function in $H^{s}(\mathbb{R})  $ such that $f\left(
x\right)  =0$ for $x\leq0$, then its restriction $f_{0}$ to $[0,\infty [ \,$
belongs to $H^{s}([0,\infty [ \,)  $, and we have:%
\[
\left\Vert f_{0}\right\Vert _{H^{s}([0,\infty [ \,) }=\left\Vert
f\right\Vert _{H^{s}(  \mathbb{R})  }%
\]
If $s=n$ is an integer, the norm on $H^{s}([0,\infty [ \,)$ turns
out to be equivalent to the following one:%
\[
|||f|||_{H^{s}}^{2}=\int_{0}^{\infty}\left[  f^{2}+\sum_{k=1}^{n}\left(
\frac{d^{k}f}{dx^{k}}\right)  ^{2}\right]  dx.
\]
To establish  properties of translations in $H^{s}([0,\infty [ \,),$ we need
to know if there is a continuous linear embedding of $H^{s}([0,\infty [ \,)$ into
$H^{s}(\mathbb{R}),$ i.e. to know if the restriction operator has a continuous
right-inverse.  Fortunately, as we are in a Hilbert space setting,
this problem is easy to solve.
Let $s \geq 0$ and let $H_{-}^{s}$ be the subset
 of functions in $H^{s}(\mathbb{R})$ with support in $]-\infty,0 ]$,
   so that $f\in H_{-}^{s}$ \ if and only if $f \in H^{s}(\mathbb{R})$ and $f( x) =0$
   for all $x >0.$ $H_{-}^{s}$ is a closed subspace of $H^{s}(\mathbb{R}).$ Two functions
$f_{1}, f_{2} \in H^{s}(\mathbb{R})$ have the same restriction to $[0,\infty [ \,$
iff $f_{1}- f_{2} \in H_{-}^{s}.$ This means exactly that $H^{s}([0,\infty [ \,)$
is a quotient space:
$H^{s}([0,\infty [ \,)=H^{s}(\mathbb{R})/H_{-}^{s}.$
Introducing the notation $\oplus$ for the Hilbert space direct sum, we have the
following result, which proof we omit since its trivial:
\begin{proposition} \label{prop Hs decomp}
For $s \geq 0$ we have: \\
i) $H^{s}(\mathbb{R})=H^{s}([0,\infty [ \,) \oplus H_{-}^{s}.$ \\
ii) Let $M$ be the orthogonal complement of $H_{-}^{s}$ in $H^{s}(\mathbb{R})$ w.r.t.
the scalar product in $H^{s}(\mathbb{R}),$ let $\kappa$ be the canonical projection
of $H^{s}(\mathbb{R})$ on $H^{s}([0,\infty [ \,)$ and let $\iota$ be the canonical
bijection of $H^{s}([0,\infty [ \,)$ onto $M.$ Then $\kappa$ is continuous, $\iota$
is a Hilbert space isomorphism, $\kappa \iota $ is the identity map on
$H^{s}([0,\infty [ \,)$ and $\iota \kappa $ is the orthogonal projection map in $H^{s}(\mathbb{R})$ on
 $M.$
\end{proposition}
We note that $\iota$ is a continuous operator extending functions on $[0,\infty [ $
to functions on $\mathbb{R}$ and that
$\left\Vert f\right\Vert _{H^{s}([0,\infty [ \,)}
   =\left\Vert \iota f \right\Vert _{H^{s}(\mathbb{R})}.$

The dual space of $H^{s}([0,\infty [ \,)$ can easily be characterized in terms
of distributions. For $s \geq 0,$ $H^{s}([0,\infty [ \,)=H^{s}(\mathbb{R})/H_{-}^{s},$
so
\begin{equation}
(H^{s}([0,\infty [ \,))'=\left\{ f\in H^{-s}(\mathbb{R}) \, | \, \sesq{f}{g}=0 \;\;\forall \,
g \in H_{-}^{s}\right\}.
\end{equation}
For $s \geq 0,$ we define $H^{-s}([0,\infty [ \,)$ to be the closed subspace
of all distributions in $H^{-s}(\mathbb{R})$ with support in $[0,\infty [ \, .$
It then follows that $(H^{s}([0,\infty [ \,))'$ can be identified with  $H^{-s}([0,\infty [ \,).$
Since $(H^{s}([0,\infty [ \,))''=H^{s}([0,\infty [ \,),$ it then follows that
\begin{equation} \label{H'}
(H^{s}([0,\infty [ \,))'=H^{-s}([0,\infty [ \,) \;s \in \mathbb{R}.
\end{equation}
If $s  \in \mathbb{R},$ then the constant function taking the value $1$ is not in
$H^{s}([0,\infty [ \,).$
If $s > 1/2,$ then even every function in $H^{s}([0,\infty [ \,) $ converges to zero at $\infty.$
For this reason, we will need a larger class of distributions containing the constant
functions. Let  $s  \in \mathbb{R}$ and let $f$ be a distribution with support
in $[0,\infty [$ such that it admits the decomposition $f=g+a,$ where
$g \in H^{s}([0,\infty [ \,)$ and $a \in \mathbb{R}.$  This decomposition of $f$
is then unique and the set of all such distributions is naturally given the Hilbert space
structure $H^{s}([0,\infty [ \,) \oplus \mathbb{R}.$ The norm of $f=g+a$ is then
given by
\[
\Vert f\Vert^{2}=\Vert g\Vert _{H^{s}([0,\infty [ \,)}^{2}+a^{2}.
\]
This unique decomposition property leads us to the following
\begin{definition} \label{Es}
For $s \in \mathbb{R},$ set $E^{s}(  [0,\infty [ \,)  =H^{s}([0,\infty [ \,)  \oplus
\mathbb{R}$ with the corresponding Hilbert space norm.
If $f\in E^{s}(  [0,\infty [ \,)$ and if
$g\in H^{s}(  [0,\infty [ \,)$ and  $a\in \mathbb{R}$ are related by the unique decomposition
$f =g+a,$ then the norm of $f$ is given by
\[
\left\Vert f\right\Vert _{E^{s}}^{2}=\left\Vert g\right\Vert _{H^{s}}^{2}+a^{2}.
\]
\end{definition}
The dual $(E^{s}(  [0,\infty [ \,))',$ of $E^{s}(  [0,\infty [ \,)$ is identified
with $(H^{s}(  [0,\infty [ \,))'  \oplus  \mathbb{R} \approx E^{-s}(  [0,\infty [ \,)$
by extending the bi-linear form, defined in (\ref{sesqprod}), to
$E^{-s}(  [0,\infty [ \,) \times E^{s}(  [0,\infty [ \,):$
\begin{equation}
\sesq{F}{G}=ab + \sesq{f}{g},
\label{ext-sesqprod}
\end{equation}
where $F=a+f \in E^{-s}(  [0,\infty [ \,),$ $G=b+g \in E^{s}(  [0,\infty [ \,),$
$a,b \in \mathbb{R},$ $f \in H^{-s}(  [0,\infty [ \,)$ and $g \in H^{s}(  [0,\infty [ \,).$

For all the Sobolev spaces $H^{s}$ we have introduced, and also for the spaces $E^{s},$
there are  two natural realizations of the dual space. Let us consider only the case
of $E^{s}(  [0,\infty [ \,),$ the other being similar.  One possibility,
the canonical one, is to identify  $(E^{s}(  [0,\infty [ \,))'$ with
$E^{s}(  [0,\infty [ \,)$ by the scalar product in $E^{s}(  [0,\infty [ \,).$
This gives the Riesz representation in (\ref{riesz}).
Another possibility, is, as we have seen, to identify $(E^{s}(  [0,\infty [ \,))'$
with $E^{-s}(  [0,\infty [ \,),$ by the bi-linear form defined in (\ref{ext-sesqprod}).
There is a linear continuous map
$\sequiv : E^{s}(  [0,\infty [ \,) \rightarrow (E^{s}(  [0,\infty [ \,))'$
with continuous inverse,  relating the two realizations. It is defined by:
\begin{equation}
  (f,g)_{E^{s}(  [0,\infty [ \,)}=\sesq{\sequiv f}{g}, \; \forall f,g \in E^{s}(  [0,\infty [ \,).
  \label{S equiv}
\end{equation}
Now,  different realizations of the dual space leads to different realizations of adjoint
operators. Let $A$ be a closed and densely defined operator from a Hilbert-space
$H$ to $E^{s}(  [0,\infty [ \,).$
We have already defined in (\ref{adj gen 1}) its adjoint operator $A^{*}$ from $E^{s}(  [0,\infty [ \,)$
to $H$ w.r.t. the duality defined by the scalar product. Let the dual $H'$ of $H$
be realized by $H_{1}$ and the continuous bi-linear form
$\sesq{\;}{\;}_{1}: H_{1} \times H \rightarrow \mathbb{R}.$
The adjoint $A',$ w.r.t. the duality realized by  $\sesq{\;}{\;}_{1}$ and $\sesq{\;}{\;}$
is the operator from $E^{-s}(  [0,\infty [ \,)$ to $H_{1},$  defined by: 
The domain of $\mathcal{D}(A^{'})$ consists of all $y \in E^{-s}(  [0,\infty [ \,)$
for which the linear functional
\begin{equation} \label{adj gen 3}
x \mapsto \ \sesq{y}{Ax}
\end{equation}
is continuous on $\mathcal{D}(A).$ For $y \in \mathcal{D}(A^{'})$ we define
$A^{'}y$ by
\begin{equation} \label{adj gen 4}
 \sesq{A^{'}y}{x}_{1}=\sesq{y}{Ax}  \; \forall x \in \mathcal{D}(A).
\end{equation}
This defines $A^{'}y$ uniquely, since $\mathcal{D}(A)$ is dense in $H.$

We now study translation semi-groups in the different spaces we have introduced.
It follows directly from the definition (\ref{Hn}) of the norm in $H^{s}(\mathbb{R})$
and by dominated convergence that that left-translations defines a strongly
continuous group of unitary operators $\tilde{\mathcal{L}}$ in $H^{s}(\mathbb{R})$
for $s \in \mathbb{R}$ 
(similarly as to  the case of $\Phi_{1}$ in Example \ref{translation 1}):
\begin{equation} \label{lefttr 0}
(\tilde{\mathcal{L}}_{t}f)(x)=f(x+t),  \; \forall f \in H^{s}(\mathbb{R}) \; \text{and} \; t \in \mathbb{R}.
\end{equation}
Since, for $s \geq 0,$ the closed subspace $H_{-}^{s}$ of $H^{s}(\mathbb{R})$
is invariant under the semi-group $\tilde{\mathcal{L}}_{t},$ $t \geq 0,$ it defines
a semi-group $\ltrans{}$ in $H^{s}([0,\infty [ \,).$ Defining $\ltrans{}$ also on
constants $a \in \mathbb{R}$ by $\ltrans{t}a=a$ we extend the semi-group $\ltrans{}$
to $E^{s}([0,\infty [ \,),$ $s \geq 0:$
\begin{equation}
( \ltrans{t}f) (x) =f(x+t), \; \forall f \in E^{s}([0,\infty [ \,) \; \text{and} \; t \geq 0.
\label{lefttr}
\end{equation}
\begin{proposition} \label{prop L}
If $s \geq 0,$ then $\ltrans{}$ is a strongly continuous contraction semi-group
on $E^{s}([0,\infty [ \,).$ Its
infinitesimal generator, denoted $\partial,$ has domain
$\mathcal{D}(\partial)=$ $E^{s+1}([0,\infty [ \,).$
If $f\in E^{s+1}([0,\infty [ \,)$ then $\partial f =f',$
where $f'$ is the derivative of $f.$
\end{proposition}
\begin{proof}
We first observe that, in the canonical decomposition
$E^{s}(  [0,\infty [ \,)  =H^{s}([0,\infty [ \,)  \oplus \mathbb{R}$ in Definition \ref{Es},
$\ltrans{a}$ leaves the subspace $H^{s}([0,\infty [ \,)$ invariant
and acts trivially on $\mathbb{R}.$ It is therefore sufficient
to prove the statement with $E^{s}([0,\infty [ \,)$ replaced by $H^{s}([0,\infty [ \,).$

We use the notations of Proposition \ref{prop Hs decomp} and let $P=\iota \kappa$
be the orthogonal projection on $M.$
Since $\tilde{\mathcal{L}}_{t}H_{-}^{s} \subset H_{-}^{s},$
for $t \geq 0,$ it follows that $P \tilde{\mathcal{L}}_{t}(I-P)= 0.$ The group
composition law $\tilde{\mathcal{L}}_{t}\tilde{\mathcal{L}}_{u}=\tilde{\mathcal{L}}_{t+u},$
then gives for $t,u \geq 0:$
\[
(P\tilde{\mathcal{L}}_{t}P)(P\tilde{\mathcal{L}}_{u}P)=P\tilde{\mathcal{L}}_{t+u}P
\]
So, $[0,\infty [ \, \ni t \mapsto P\tilde{\mathcal{L}}_{t}P$ is a semi-group of bounded
operators on $M.$ It is a strongly continuous contraction semi-group since this is the
case for $\tilde{\mathcal{L}}_{}$ and $\|P\|=1.$

We have that $\ltrans{t}=\kappa \tilde{\mathcal{L}}_{t} \iota,$ for $t \geq 0.$
Using that $\ltrans{t}=\kappa P \tilde{\mathcal{L}}_{t} P \iota$ it easily follows
from the semi-group property of $ P \tilde{\mathcal{L}}_{t} P$ that $\ltrans{}$
is a semi-group on $H^{s}([0,\infty [ \,).$ It is a strongly continuous contraction semi-group,
since this is the case for $P\tilde{\mathcal{L}}_{}P$ and since $\|\kappa\|=\|\iota\|=1.$
 Let
$\partial$ be the infinitesimal generator of $\ltrans{}.$ By the definition of  $\ltrans{}$
it follows that
$\mathcal{D}(\partial)=\{f \in H^{s}([0,\infty [ \,) \; | \; f' \in H^{s}([0,\infty [ \,)\}$
and $\partial f=f'$ for $f \in \mathcal{D}(\partial).$ But
$H^{s+1}([0,\infty [ \,)=\{f \in H^{s}([0,\infty [ \,) \; | \; f' \in H^{s}([0,\infty [ \,)\},$
which proves the proposition.
%\AAAqed
\end{proof}
\begin{example} \label{translation 3} \text{} \\ \normalfont
Let $\ltrans{t}' : E^{-s}([0,\infty [ \,) \rightarrow E^{-s}([0,\infty [ \,)$ be the adjoint
of $\ltrans{t},$ $t \geq 0$ in Proposition \ref{prop L},
w.r.t. duality defined by the bilinear form $\sesq{\;}{\;}.$
$\ltrans{}'$ is then a semi-group of right-translations on the space of distributions
$E^{-s}([0,\infty [ \,).$ Loosely speaking $(\ltrans{t}'f)(x)=f(x-t).$
Let $s \geq 1.$
hen the generator $\partial'$ has domain $E^{-s+1}([0,\infty [ \,)$ and
$-\partial'$ is the derivative of distributions, so
$(\partial'f)(x)=-df(x)/dx$ if $f$ is a differentiable function.
One is easily convinced that the expressions for $\ltrans{t}^{*}$ and $\partial^{*}$
are more complicated. 
\end{example}
Some general references for this subsection are: \cite{Adams 03}, \cite{Calderon}, \cite{Horm}.
\subsection{Infinite-dimensional Brownian motion}
In this sub-section we consider a separable Hilbert space $E$ and an index-set
$\mathbb{I}$ with the cardinality equal to the dimension of $E.$
The space $E$ can be infinite-dimensional or finite-dimensional.
There is given  a family $W^{i},$
$i\in\mathbb{I}$ of standard independent Brownian motions on a complete filtered
probability space $(\Omega,P,\mathcal{F},\mathcal{A}).$ The filtration
$\mathcal{A}=\{\mathcal{F}_{t}\}_{0\leq t \leq T},$ is generated by the
$W^{i},$ and $\mathcal{F}=\mathcal{F}_{T}.$
\begin{definition}
A standard cylindrical Brownian motion $W_{t},\ 0\leq t\leq T$, on $E$ is a
 sequence $e_{i}W_{t}^{i}  ,\ i\in\mathbb{I}$ of $E$-valued processes,
 where the $e_{i}$ are the elements of an orthonormal basis of $E$ and the
$ W_{t}^{i},$ $i\in\mathbb{I},$  are
independent real-valued standard Brownian motions on a filtered probability space
$(\Omega,P,\mathcal{F},\mathcal{A}).$
\end{definition}
From now on, given a standard cylindrical Brownian motion %
$W,$ we shall
write informally $W_{t}=\sum_{i\in\mathbb{I}}W_{t}^{i}e_{i}$.
If $\mathbb{I}$ is finite, we have:%
\[
\Vert W_{t} \Vert ^{2}=\sum_{i\in
\mathbb{I}}\left\Vert W_{t}^{i}\right\Vert^{2} < \infty \; \; \text{a.s.}
\]
and $W_{t}$ is a stochastic process with values in $E$.
If $\mathbb{I}$ is infinite, then for every $t$ the right-hand side is the sum
of infinitely many i.i.d. positive random variables, which does not converge
in any reasonable way. In that case, the formula $W_{t}=\sum_{i\in\mathbb{I}%
}W_{t}^{i}e_{i}$ cannot be understood as an equality in $E$, and must be given
another meaning.
\begin{proposition}
If $W_{t}=\sum_{i\in\mathbb{I}}W_{t}^{i}e_{i}$ is a standard cylindrical
Brownian motion, then, for every $f\in E$ with $\left\Vert f\right\Vert =1$,
the real-valued stochastic process $W_{t}^{f}$ defined by
\begin{equation}
W_{t}^{f}=\sum_{i\in\mathbb{I}}\left(  e_{i},f\right)  W_{t}^{i} \label{13}%
\end{equation}
is a standard Brownian motion on the real line.
\end{proposition}
\begin{proof}
If $\mathbb{I}$ is finite, the result is obvious. Let us then consider the
case when $\mathbb{I}$ $=\mathbb{N}$. We first have to check if the right-hand side is
well-defined. By Doob's inequality for martingales:%
\begin{align*}
E\left[  \sup_{0\leq t\leq T}\left\vert \sum_{i=n}^{n+p}\left(  e_{i}%
,f\right)  W_{t}^{i}\right\vert ^{2}\right]   &  \leq4E\left[  \left\vert
\sum_{i=n}^{n+p}\left(  e_{i},f\right)  W_{T}^{i}\right\vert ^{2}\right] \\
&  \leq4T\sum_{i=n}^{n+p}\left(  e_{i},f\right)  ^{2}\rightarrow0.
\end{align*}
This implies that the right-hand side of (\ref{13}) converges in probability
to a continuous process.
Since each finite sum is Gaussian, so is the limit, and the result follows.
%\AAAqed
\end{proof}

So, in the case when $\mathbb{I}$ is infinite, the r.h.s. of
$W_{t}=\sum_{i\in\mathbb{I}}W_{t}^{i}e_{i}$ makes no sense in $E$, but every projection
does. Equation (\ref{13}) can be rewritten as:%
\[
\forall f\in E,\ \ \ \left(  W_{t},f\right)  =\sum_{i\in\mathbb{I}}\left(
e_{i},f\right)  W_{t}^{i}.
\]

We will now show that the stochastic integrals with respect to cylindrical
Brownian motion make sense, provided the integrand satisfies a strong
integrability condition.
Consider the space $\mathcal{HS}(E,F) $ of all Hilbert-Schmidt operators from $E$
into a Hilbert space $F.$ Let the space $\mathcal{L}^{2}\left(  \mathcal{HS}(E,F) \right) $
consist of all progressively measurable processes $A$ with values in
the Hilbert space $\mathcal{HS}(E,F) $, such that:
\[
E\left[\int_{0}^{T}\left\Vert A_{t}\right\Vert _{\mathcal{HS}}^{2}dt\right]<\infty.
\]
Recall that we have, according to Definition \ref{H-S}:
\[
\left\Vert A_{t}\right\Vert _{\mathcal{HS}}^{2} %
=\sum_{n=0}^{\infty}\left\Vert A_{t}e_{n}\right\Vert ^{2}, %
\]
where $\left(e_{n}\right)_{n\in \mathbb{N}}$ is any orthonormal basis of $E.$
\begin{theorem} \label{stoch int}
The stochastic integral:
\[
\int_{0}^{T}A_{t}dW_{t}
\]
is well-defined for every %
 process $A\in\mathcal{L}^{2}\left(
\mathcal{HS}(E,F) \right)  .$ It is a continuous martingale with values in $E$, and
we have the usual isometry:
\[
\left\Vert \int_{0}^{T}A_{t}dW_{t}\right\Vert _{L^{2}}^{2}
=\int_{0}^{T} E\left[ \left\Vert A_{t}\right\Vert _{\mathcal{HS}}^{2}\right]dt.
\]
\end{theorem}
In other words, the random variable $\int_{0}^{T}A_{t}dW_{t}$ has mean
$0$ and its variance is
$\sum_{n=0}^{\infty}\int_{0}^{T} E\left[ \left\Vert A_{t}e_{n}\right\Vert ^{2}\right]dt,$
 the sum of the variances of the independent sources of Gaussian noise.

As usual, by localization the stochastic integral can be extended to a wider class of
processes. Denote by $\mathcal{L}_{loc}^{2}\left(  \mathcal{HS}\right)  $ the
set of all progressively measurable processes with values in $\mathcal{HS}$, such that:%
\[
P\left[  \int_{0}^{T}\left\Vert \Phi\right\Vert _{\mathcal{HS}}^{2}
dt<\infty\right]  =1.
\]
Then the stochastic integral defines a continuous local martingale.

Some general references for this subsection are: 
\cite{DaPrato-Zabczyk},  \cite{Kall-Xiong},  \cite{Mikul-Rozo98},
\cite{Mikul-Rozo99}, \cite{Nual71}.

\section{The dynamics of bond prices}
\subsection{The non-parametric framework} \label{non-param frame}
From now on, and for the rest of the paper, we are given a finite time interval
of possible trading times $\mathbb{T}=[0, \timeh]$ and we are given  a family $W^{i},$
$i\in\mathbb{I\ }$ of standard independent Brownian motions on a complete filtered
probability space $(\Omega,P,\mathcal{F},\mathcal{A}),$ the filtration
$\mathcal{A}=\{\mathcal{F}_{t}\}_{0\leq t \leq \timeh},$ is generated by the
$W^{i},$ and $\mathcal{F}=\mathcal{F}_{\timeh}.$  The family $\mathbb{I}$
itself can be finite or infinite, in
which case we take $\mathbb{I}=\mathbb{N}.$ Let $\ell^{2}(\mathbb{I}),$ be
the Hilbert space of all real sequences $x=(x_{i})_{i \in \mathbb{I}},$ such that
$\|x\|_{\ell^{2}(\mathbb{I})}=(\sum_{i\in \mathbb{I}} (a_{n})^{2})^{1/2}<\infty$.
So, when $\mathbb{I}$ has a finite number $\numrand$ of elements, then
$\ell^{2}(\mathbb{I}) =\mathbb{R}^{\numrand}.$ Often we write just $\ell^{2}$
for $\ell^{2}(\mathbb{I}).$

Heath, Jarrow and Morton (henceforth HJM) were the first to study the term
structure of interest rates in a non-parametric framework. Their basic idea
(see \cite{HJM92}) consists of writing one equation for the price of every
zero-coupon at time $t$. Denoting by $\hat{B}_{t}(T)$ the price at time $t$ of
a zero-coupon bond maturing at time $T\geq t,$ the HJM equation has the following
form:
\begin{equation}
\hat{B}_{t}(T)=\hat{B}_{0}(T)+\int_{0}^{t}\hat{B}_{s}(T)a_{s}(T)ds+\int
_{0}^{t}\sum_{i\in\mathbb{I}}\hat{B}_{s}(T)v_{s}^{i}(T)dW_{s}^{i},\ \ 0\leq
t\leq T \label{a2}%
\end{equation}
There are infinitely many such equations, one for each maturity $T\geq t.$

The trend $a_{t}\left(  T\right)  $ and the volatilities $v_{t}^{i}(T)$ are
supposed to be progressively measurable processes, which means, for instance, that they
could be functions of all the $\hat{B}_{s}(S),$ $S\geq s$ and $s \leq t.$ In due course, we will
make further assumptions so as to ensure that equations such as (\ref{a2})
make mathematical sense.

Let us discount all prices to $t=0,$ by the spot interest rate  $r_{t},$
which in terms of the zero-coupon bond price is given by
\begin{equation} \label{spot rate}
\spotrate{t}=-\frac{\partial \hat{B}_{t}(T)}{\partial T}\Big |_{T=t}.
\end{equation}
The discounted prices of zero-coupons are now:
\begin{equation} \label{discounted B}
B_{t}(T)=\hat{B}_{t}(T)\exp(-\int_{0}^{t}\spotrate{s}ds)
\end{equation}
and the equations (\ref{a2}) become:
 \begin{equation}        \label{a2'}
B_{t}(T)=B_{0}(T)+\int_{0}^{t}B_{s}(T)(a_{s}(T)-r_{s})ds+\int_{0}^{t}%
\sum_{i\in\mathbb{I}}B_{s}(T)v_{s}^{i}(T)dW_{s}^{i},\ \ 0\leq t\leq T
\end{equation}
and, again, there is one such equation for every maturity $T\geq t$. Note the
boundary condition $\hat{B}_{T}(T)=1,$ and hence, from (\ref{discounted B}):
\begin{equation} \label{52}
B_{t}(t)=\exp(-\int_{0}^{t}r_{s}ds) .
\end{equation}
\subsection{The bond dynamics in the moving frame}
For every $x\geq0$, we denote by $\zcpx{t}(x)$ the price and by  $\zcpxd{t}(x)$ the discounted
price at time $t$ of a zero-coupon maturing at time $t+x$. The stochastic
processes $B_{t}(T)  $ and $p_{t}(x)  $ are related
by:
\[
p_{t}(x)=B_{t}(t+x).
\]
In other words, as explained in the introduction, instead of dating events by their distance from a fixed
origin, defined to be $t=0$, we are dating them by their distance from today:
we are using a time frame which moves with the observer. The equation for
$p_{t}$ in the moving frame, is easily obtained from (\ref{a2'}). For every
$x\geq0$, we have:\
\begin{equation}
\begin{split}
p_{t}(x)=&p_{0}(t+x)+\int_{0}^{t}p_{s}(t-s+x)m_{s}(t-s+x)ds \\
   &+\int_{0}^{t}
\sum_{i\in\mathbb{I}}p_{s}(t-s+x)\sigma_{s}^{i}(t-s+x)dW_{s}^{i}, \label{51}
\end{split}
\end{equation}
where
\begin{equation}
\drift{t}(x)=a(t,t+x)-\spotrate{t} \; \; \text{and} \; \; \vol{i}{t}(x)=v^{i}(t,t+x),
\label{drift-vol}
\end{equation}
for all $0 \leq t \leq \timeh$ and $x \geq 0.$ 
Here, again, the trends $t \mapsto m_{t}(x)$ and the volatilities $t \mapsto \sigma_{t}^{i}(x)$
are progressively measurable processes.

Instead of looking at (\ref{51}) as an infinite family of coupled equations, one for
each $x\geq0$, we shall interpret it as a single equation describing the dynamics
of an infinite-dimensional object, the curve $x \mapsto p_{t}(x),  $
 which will be seen as a vector $p_{t}$ in the Hilbert space
$E^{s}([0, \infty[\,),$ for some fixed $s>1/2$, chosen so that the functions $m_{t}$ and
$\sigma_{t}^{i}$ belong to $E^{s}([0, \infty[\,).$

Let $\ltrans{}$ be the  semi-group  left translations on  $E^{s}([0, \infty[\,)$
(see formula (\ref{lefttr}) and Proposition \ref{prop L}). From now on we shall just wright
$E^{s}$ instead of $E^{s}([0, \infty[\,),$ when there is no risk of confusion.
The equations in
(\ref{51}) can be rewritten as one equation in $E^{s}$:%
\begin{equation} \label{dynam discount p}
\ p_{t}=\mathcal{L}_{t}p_{0}+\int_{0}^{t}(\mathcal{L}_{t-s}(p_{s}%
m_{s}))ds+\int_{0}^{t}\sum_{i\in\mathbb{I}}(\mathcal{L}_{t-s}(p_{s}\sigma
_{s}^{i}))dW_{s}^{i}.
\end{equation}
\begin{theorem} \label{mild}
Let $s> 1/2.$ Assume that $p_{0} \in E^{s}$ and assume that $m_{t}$
and the $\sigma_{t}^{i},$ $i\in\mathbb{I},$ are progressively measurable processes in
$E^{s}$ satisfying:%
\begin{equation}
\int_{0}^{\timeh}(\Vert m_{t}\Vert_{E^{s}}+\sum_{i\in\mathbb{I}}\Vert
\sigma_{t}^{i}\Vert_{E^{s}}^{2})dt<\infty\;\ \ \text{a.s.} \label{55}%
\end{equation}
Then equation (\ref{dynam discount p})  defines a unique process  $p$ in
$E^{s}$ satisfying:
\begin{equation}  \label{56}
\int_{0}^{\timeh}(\Vert p_{t}\Vert_{E^{s}}+\Vert p_{t}m_{t}\Vert_{E^{s}}%
+\sum_{i\in\mathbb{I}}\Vert p_{t}\sigma_{t}^{i}\Vert_{E^{s}}^{2}%
)dt<\infty\;\text{a.s.}
\end{equation}
The process $p$ has continuous trajectories in $E^{s},$
\begin{equation} \label{bond dyn sol p}
p_{t}=\exp\left\{ \int_{0}^{t}\mathcal{L}_{t-s}\left((m_{s}-\frac{1}{2}\sum
_{i\in\mathbb{I}}(\sigma_{s}^{i})^{2})ds+\sum_{i\in\mathbb{I}}\sigma_{s}
^{i}dW_{s}^{i}\right) \right\}  \mathcal{L}_{t}p_{0}.
\end{equation}
and if $p_{0} \in H^{s}$ then the process $p$ takes its values in $H^{s}.$
If $p_{0} \in E^{s}$ satisfies $p_{0} \geq 0$ (resp. $p_{0} > 0$),
i.e. $p_{0}(x) \geq 0$ (resp. $p_{0}(x) > 0$) for all $x \geq 0,$
then so does $p_{t}.$
\end{theorem}
For a proof of this theorem see Lemma A.1 of \cite{I.E.-E.T bond th}, which is reproduced
in the appendix of this article (Lemma \ref{existence lemma}).
Note that equation (\ref{dynam discount p}) implies that $p_{0}$ is the value
of $p_{t}$ for $t=0$.

A word here about the choice of function spaces. Assuming that $p_{t}$ belongs
to $H^{s}$ for some $s>1/2$ is minimal: it is basically saying that the
zero-coupon prices depend continuously on time to maturity and go to zero at
infinity. This, however, is too strong a requirement for $m_{t}$ and the
$\sigma_{t}^{i}$: we cannot expect the trend and the volatilities to go to
zero when the time to maturity increases to infinity. This is why we are
assuming that $m_{t}$ and the $\sigma_{t}^{i}$ belong to $E^{s}$.
To simplify the mathematical formalism and also to include interest rate models,
with vanishing long term rates, we have permitted that $p_{t} \in E^{s}.$
Now according to Theorem \ref{mild}, $p_{t}$ is  in-fact in   $H^{s}$ if  $p_{0} \in H^{s}.$

Condition (\ref{55}) implies that $\sum_{i\in\mathbb{I}}\Vert\sigma_{t}%
^{i}\Vert_{E^{s}}^{2}$ is finite for almost every
$(t,\omega) \in \mathbb{T} \times \Omega.$ This means, when $\mathbb{I}=\mathbb{N},$ that the
operator $\sigma_{t}$ from $\ell^{2}(\mathbb{I})$ to $E^{s}$ defined by:%
\begin{equation} \label{vol hs}
\sigma_{t}e_{i}=\sigma^{i}_{t},\ \ \ i\in\mathbb{I},
\end{equation}
where $e_{i}$ are the elements of the standard basis of $\ell^{2}(\mathbb{I}),$
is Hilbert-Schmidt a.e. $(t,\omega).$ We have
\[
\left\Vert \sigma_{t}\right\Vert _{\mathcal{HS}(\ell^{2},E^{s})}^{2}=\sum_{i\in\mathbb{I}%
}\Vert\sigma_{t}^{i}\Vert_{E^{s}}^{2}.
\]
We shall refer to $\sigma$ as the \emph{volatility operator} process. It takes its values in
$\mathcal{HS}(\ell^{2},E^{s})$ and when we say that it is progressively measurable,
it is meant that all the $\sigma^{i}$ are progressively measurable.

We can now, using the stochastic integral introduced in Theorem \ref{stoch int},
rewrite equation (\ref{dynam discount p}) on a more compact form in $E^{s},$ where  $s> 1/2:$
\begin{equation} \label{dynam discount p 1}
\ p_{t}=\mathcal{L}_{t}p_{0}+\int_{0}^{t}\mathcal{L}_{t-s}(p_{s} m_{s})ds
+\int_{0}^{t}\mathcal{L}_{t-s}(p_{s}\sigma_{s})dW_{s}.
\end{equation}
This makes sens in $E^{s}.$ Indeed, the only difference with equation (\ref{dynam discount p})
is the last term on the r.h.s. When condition (\ref{55}) is satisfied then the
volatility operator $\sigma_{u},$ defined by (\ref{vol hs}), from $\ell^{2}$ to $E^{s},$
is Hilbert-Schmidt a.e. $(u,\omega).$ Since pointwise  multiplication of functions in $E^{s}$
is a continuous operation for $s> 1/2$ it follows that the linear operator $x \mapsto p_{u}\sigma_{u}x,$
 from $\ell^{2}$ to $E^{s},$ is Hilbert-Schmidt a.e. $(u,\omega).$ $\ltrans{v}$ is bonded for every $v \geq 0,$
so the integrand is a progressively measurable $\mathcal{HS}(\ell^{2},E^{s})$-valued process
satisfying the conditions of Theorem \ref{stoch int}.

A process $p$ with values in $E^{s}$ satisfying (\ref{dynam discount p 1})
(or equivalently (\ref{dynam discount p})) and (\ref{56})
will be called a \emph{mild solution} of the bonds dynamics. %

Note that we are not worrying about the boundary condition (\ref{52}) at this
time, because it does not make mathematical sense: how do we define $r_{t}%
$\ ?\ This will be taken care of in the next section.
\subsection{Smoothness of the zero-coupon curve.}
Another way to proceed is to write (\ref{51}) in differentiated form. For fixed $x \geq 0,$
a formal calculation using It\^o's lemma and which can be rigorously justified  gives:
\begin{equation} \notag%
\begin{split}
dp_{t}(x)-&p_{t}\left(  x\right)  m_{t}\left(  x\right)  dt-\sum_{i\in\mathbb{I}}p_{t}(x)\sigma_{t}^{i}(x)dW_{t}^{i} \\
&=\bigg(\frac{\partial}{\partial t}p_{0}(t+x) 
+\int_{0}^{t}\frac{\partial}{\partial t}\left(  p_{s}(t-s+x)m_{s}(t-s+x)\right)  \ ds \\
&+\int_{0}^{t}\frac{\partial
}{\partial t}  \sum_{i\in\mathbb{I}}p_{s}(t-s+x)\sigma_{s}%
^{i}(t-s+x) dW_{s}^{i}\bigg) \, dt.
\end{split}
\end{equation}
In the expression on r.h.s. we can replace $\partial/\partial t$ by $\partial/\partial x,$
since $p_{0}$ and the integrands on the r.h.s. are functions of $t+x.$ Derivation w.r.t.
$x$ under the integral then gives:
\begin{equation} \notag%
\begin{split}
&dp_{t}(x)-p_{t}\left(  x\right)  m_{t}\left(  x\right)  dt-\sum_{i\in\mathbb{I}}p_{t}(x)\sigma_{t}^{i}(x)dW_{t}^{i} \\
&=\bigg(\frac{\partial}{\partial x} \Big( p_{0}(t+x)
+\int_{0}^{t} p_{s}(t-s+x)m_{s}(t-s+x)  \ ds \\
&+\int_{0}^{t} \sum_{i\in\mathbb{I}}p_{s}(t-s+x)\sigma_{s}^{i}(t-s+x) dW_{s}^{i} \Big)\bigg) \, dt.
\end{split}
\end{equation}
The l.h.s. is equal to $((\partial/\partial x) p_{t}(x)) \, dt,$ according to (\ref{51}), so
\begin{equation} \label{diff form3}
dp_{t}(x)-p_{t}\left(  x\right)  m_{t}\left(  x\right)  dt-\sum_{i\in\mathbb{I}}p_{t}(x)\sigma_{t}^{i}(x)dW_{t}^{i}
=\big(\frac{\partial}{\partial x}  p_{t}(x)\big) \, dt,
\end{equation}
for all $x \geq 0$ and $t \in \mathbb{T}.$

Introducing the infinitesimal generator $\partial$ of the semi-group
$\ltrans{}$ (see Proposition \ref{prop L}), this can be understood as an equation in $E^{s}:$
\begin{equation}
dp_{t}=(\partial p_{t}+p_{t}  m_{t})dt+\sum_{i\in\mathbb{I}}p_{t}\sigma_{t}^{i}dW_{t}^{i}  \label{SPDE pd}
\end{equation}
or equivalently: %
\begin{equation}
p_{t}=p_{0}+\int_{0}^{t}(  \partial p_{s}+p_{s}m_{s})
ds+\int_{0}^{t} \sum_{i\in\mathbb{I}}p_{s}\sigma_{s}^{i}dW_{s}^{i}.
\label{SPDE p}
\end{equation}

Equation (\ref{dynam discount p}) is the integrated version of (\ref{SPDE p}),
w.r.t. the semi-group $\ltrans{}.$
The connection between formulas (\ref{SPDE p}) and (\ref{dynam discount p})
is similar to the \emph{variations of constants} formula for ODE's in finite dimension.

We now have to give some mathematical meaning to equation (\ref{SPDE p}). This
will require beefing up the existence conditions given in Theorem \ref{mild}.
The following corollary follows from applying Theorem \ref{mild} with $s+1$
instead of $s:$
\begin{corollary} \label{strong}
Let $s> 1/2.$ Assume that $p_{0} \in \mathcal{D}(\partial) =  E^{s+1}$ and assume
that $m_{t}$ and the $\sigma_{t}^{i},$ $i\in\mathbb{I}$ are progressively measurable
processes with values in $E^{s+1}$ satisfying
\begin{equation}
\int_{0}^{\timeh}(\Vert m_{t}\Vert_{E^{s+1}}+\sum_{i\in\mathbb{I}}\Vert
\sigma_{t}^{i}\Vert_{E^{s+1}}^{2})dt<\infty\;\ \ \text{a.s} \label{54}.
\end{equation}
Then the mild solution $p,$ in Theorem \ref{mild}, of the bonds dynamics satisfies
the following condition:
\begin{equation} \label{strong cond}
 p_{t}\in E^{s+1} \; \text{and} \; 
 \int_{0}^{\timeh}(\Vert p_{t}\Vert_{E^{s+1}}+\Vert p_{t}m_{t}%
\Vert_{E^{s}}+\sum_{i\in\mathbb{I}}\Vert p_{t} \sigma_{t}^{i}\Vert_{E^{s}}^{2}) \ dt<\infty\;\ \text{a.s.}%
\end{equation}
Equation (\ref{SPDE p}) holds for every $t.$ In addition $p$ has
continuous paths in $E^{s+1}$ and $p_{t} \in H^{s+1}$ if $p_{0} \in H^{s+1}.$
\end{corollary}
By definition a solution of equation (\ref{dynam discount p}) is called a strong solution
of the equation (\ref{SPDE p}), when  condition (\ref{strong cond}) is satisfied.
Here we shall say that $p$ is a \emph{strong solution} of the bonds dynamics.

As a consequence, in the situation of Corollary \ref{strong}, the term structure
$x \mapsto p_{t}(x)$ is
$C^{1}$ for every $t,$ and interest rates are well defined. The instantaneous forward
rate $\fwrate{t}(x)$ contracted at $t \in \mathbb{T}$ for time to maturity $x$ and
the spot rate $\spotrate{t}$ at time $t,$ for instance, are defined by:
\begin{equation} \label{intrest}
\fwrate{t}(x) =-\frac{\partial \log p_{t}(x)}{\partial x}=-\frac{(\partial p_{t})(x)}{p_{t}(x)}
\quad \text{and} \quad \spotrate{t}=\fwrate{t}(0) =-\frac{\left(  \partial p_{t}\right)  \left(  0\right)  }{p_{t}\left(
0\right)  }.
\end{equation}

By  Corollary  \ref{strong}, $p$ is a strong solution and the maps
$t\mapsto p_{t}$ and $t\mapsto\partial p_{t}$ are continuous from
$\mathbb{T}$ into $E^{s}$, and hence into $C^{0}(
[0,\infty [\, )$ endowed with the topology of uniform convergence.
So $p_{s}(  0)  $ and $(\partial p_{s})\left(
0\right)  $ converge to $p_{t}(  0)  $ and $(\partial p_{t})\left(  0\right),  $
when $s\rightarrow t.$ In other words, $r_{t}$ is a continuous function of $t,$
when $p_{t}(  0)  > 0$ for all $t\in \mathbb{R}.$

We are now able to make sense of the boundary condition (\ref{52}), which we rewrite
in terms of $p:$
\begin{equation}
\zcpxd{t}(0)=\exp (\int_{0}^{t} \frac{(\partial \zcpxd{s})(0)}{\zcpxd{s}(0)}ds),
 \label{p bound cond}
\end{equation}
for every  $ t \in \mathbb{T}.$ 
\begin{proposition} \label{prop. boundary cond}
Let $s> 1/2.$
Assume that $m_{t}$ and the $\sigma_{t}^{i}$ are progressively measurable processes with
values in $E^{s+1}$ satisfying (\ref{54}) and
\begin{equation}
m_{t}(0)=0,\;\sigma_{t}^{i}(0)=0\;\ \ \forall i\in\mathbb{I} \label{63}%
\end{equation}
and assume that $p_{0}$ satisfies
\begin{equation} \label{p0}
p_{0}  \in E^{s+1}, \; \;  p_{0}(0)=1,  \; \; p_{0}(x) >0 \ \ \ \forall x \geq 0.
\end{equation}
Then the solution of
the bond dynamics, given by Corollary \ref{strong}, satisfies the boundary
condition (\ref{p bound cond}).
\end{proposition}
\begin{proof}
Since $m_{t}$ and the $\sigma_{t}^{i}$ take values in $E^{s+1}$, they are
continuous function on $[0, \infty[ \, ,$ and condition (\ref{63}) makes sense.
As $p_{0} >0$ it follows from Proposition \ref{mild} that $p_{t} >0.$ We have
shown that, if $p_{t}$ is a strong and strictly positive solution of the bond dynamics,
then $r_{t}$ given by (\ref{intrest}) is a continuous function of $t$.
Writing conditions (\ref{63}) into equation
(\ref{SPDE p}), we get:%
\begin{align*}
p_{t}(  0)   &  =p_{0}(  0)  +\int_{0}^{t}( (\partial p_{s})(  0)  +p_{s}(  0)  m_{s}(
0)  )  ds+\int_{0}^{t} \sum_{i\in\mathbb{I}}p_{s}(  0) \sigma_{s}^{i}(  0)  dW_{s}^{i}\\
&  =1+\int_{0}^{t}(\partial p_{s})(  0)  ds
  =1-\int_{0}^{t}r_{s}p_{s}(  0)  ds.
\end{align*}
In other words, $\varphi(t)  =p_{t}(0)  $ must
satisfy the differential equation $\varphi^{\prime}(t)
=-r_{t}\varphi (t)  $, with the initial condition $\varphi (0)  =1$. The result follows.
%\AAAqed
\end{proof}

When we get to optimizing portfolios, we will need $L^{p}$ estimates on the
solutions of the bond dynamics. They are provided by the following result:
\begin{theorem}\label{Th price compatible strong cond.}
Let  $q(t)=\zcpxd{t}/\ltrans{t}\zcpxd{0}$ and $\hat{q}(t)=\zcpx{t}/\ltrans{t}\zcpx{0}.$
If $\zcpxd{0},$ $\vol{}{}$ and $\drift{}$ in Proposition \ref{prop. boundary cond}  %
also satisfy the following additional conditions:
\begin{equation}
 E( (\int_{0}^{\timeh}\| \vol{}{t} \|^{2}_{\mathcal{HS}(\ell^{2},E^{s+1})}dt)^{a}
          + \exp(a \int_{0}^{\timeh}\| \vol{}{t} \|^{2}_{\mathcal{HS}(\ell^{2},E^{s})}dt)) < \infty,
  \;  \forall a \in [1,\infty[
 \label{domain sigma i}
\end{equation}
and
\begin{equation}
 E(( \int_{0}^{\timeh}\| \drift{t} \|_{E^{s+1}}dt)^{a}
       + \exp(a \int_{0}^{\timeh}\| \drift{t} \|_{E^{s}}dt)) < \infty,
     \forall a \in [1,\infty[ \,,
 \label{domain drift}
\end{equation}
then the solution $\zcpxd{}$ in %
Proposition \ref{prop. boundary cond} has the following property:
\begin{equation}
 \zcpxd{},  \zcpx{}, q, \hat{q}, 1/q, 1/\hat{q}  \in
  L^{u}(\Omega, P, L^{\infty}(\mathbb{T},E^{s+1})),   \forall u \in [1,\infty[ \,.
\end{equation}
\end{theorem}
\begin{proof}
We use the notation
\begin{equation}
\tilde{\mathcal{E}}_{t}(L)=\exp (\int_{0}^{t} \ltrans{t-s}
  ((\drift{s}-\frac{1}{2}\sum_{i \in \mathbb{I}}(\vol{i}{s})^{2})ds
                   + \vol{}{s} d\wienerp{}{s}) ),
 \label{exp marting}
\end{equation}
for
\begin{equation}
 L_{t}=\int_{0}^{t}(\drift{s}ds + \vol{}{s} d\wienerp{}{s}), \quad \text{if} \; 0 \leq t \leq \timeh.
 \label{A. notation 1'}
\end{equation}
Conditions $(i)-(iv)$ of Lemma \ref{Lp norms} are satisfied for $p.$ Estimate (\ref{Lp norms 3})
of  Lemma \ref{Lp norms} then shows that
$ \zcpxd{}  \in
  L^{u}(\Omega, P, L^{\infty}(\mathbb{T},E^{s+1}))$ $ \forall u \in [1,\infty[ \,.$
By the  explicit expression (\ref{bond dyn sol p}), $q=\tilde{\mathcal{E}}(L),$
so it follows from Lemma \ref{Lp norms} that the conclusion  holds true also for $q.$

Let %
$N_{t}=\int_{0}^{t}((- \drift{s}+\sum_{i \in \mathbb{I}}(\vol{i}{s})^{2})ds
                - \sum_{i \in \mathbb{I}}\vol{i}{s} d\wienerp{i}{s}).$
Then  $1/q=\tilde{\mathcal{E}}(N).$
According to conditions (\ref{domain sigma i}), (\ref{domain drift}), the
conditions $(i)-(iv)$ of Lemma \ref{Lp norms}
(with $N$ instead of $L$) are satisfied.
We now apply estimate
(\ref{Lp norms 3}) %
to $1/q,$
which proves that  $1/q \in L^{u}(\Omega, P, L^{\infty}(\mathbb{T}, E^{s+1})),$
for all $u \geq 1.$ %

To prove the cases of $\hat{q}^{\alpha},$ $\alpha=1$ or $\alpha=-1,$  we note that $q(t)=\hat{q}(t)\zcpxd{t}(0).$ %
Using that the case of
$q^{\alpha}$ is already proved and H\"olders inequality, it is enough to prove that
$g  \in L^{u}(\Omega, P, L^{\infty}(\mathbb{T}, \mathbb{R})),$ where $g(t)=(\zcpxd{t}(0))^{-\alpha}.$
Since $\zcpxd{t}(0)=(\ltrans{t}\zcpxd{0})(0)(q(t))(0)$ $=\zcpxd{0}(t)(q(t))(0),$
it follows that
\[0\leq g(t) = (\zcpxd{0}(t))^{-\alpha} ((q(t))(0))^{-\alpha}.\]
By Sobolev embedding, $\zcpxd{0}$ is a continuous real valued function on $[0, \infty [$ and it is also strictly positive,
so the function $t \mapsto (\zcpxd{0}(t))^{-\alpha}$ is bounded on $\mathbb{T}.$ Once more by Sobolev embedding,
$((q(t))(0))^{-\alpha} \leq C \|(q(t))^{-\alpha}\|_{E^{s}}.$ The result
now follows, since we have already proved the case of $q^{\alpha}.$
The case of $\zcpx{}$ is so similar to the previous cases that we omit it.
%\AAAqed
\end{proof}
Under the hypotheses of Proposition \ref{prop. boundary cond}, %
$\zcpxd{t}(0)$ satisfies (\ref{p bound cond}), so it is the discount factor (\ref{52}).
It has nice properties, as follows from the second part of the proof of Theorem \ref{Th price compatible strong cond.}
\begin{corollary}\label{discount factor}
Under the hypotheses of Theorem \ref{Th price compatible strong cond.},
 if $\alpha \in \mathbb{R},$ then the discount factor $\zcpxd{t}(0)$ satisfies
     \[E(\sup_{t \in \mathbb{T}}(\zcpxd{t}(0))^{\alpha}) < \infty.\]
\end{corollary}
\begin{remark} \label{p remark}
 It follows from Theorem \ref{Th price compatible strong cond.} that for all $t \in \mathbb{T},$
$\zcpxd{t}$ and $\zcpxd{0}$ have similar asymptotic behavior. In fact for some r.v. $A>0,$
$A^{-1} \zcpxd{0}(t+x) \leq \zcpxd{t}(x) \leq A \zcpxd{0}(t+x),$ for all $t \in \mathbb{T}$ and $x\geq 0,$ where
$A$ is independent of $x$ and $t$ and $A \in L^{u}(\Omega, P)$ for all $u \geq 1.$ \\
\end{remark}

In a different context, Hilbert spaces of
forward rate curves were considered in \cite{Bj-Sv01} and \cite{Filipovic}. %
The space $E^{s},$ with $s> 1/2$ sufficiently small, contains
the image of these spaces, under the nonlinear map of forward rates to zero-coupons
prices. Or more precisely, it contains the image of subsets of forward rate curves
$f$ with positive long term interest rate, i.e. $f(x) \geq 0$ for all $x$ sufficiently
big.

\section{Portfolio theory}
In this section $s>1/2,$  $E^{s}=E^{s}[0,\infty)[ \,)$ and $\mathbb{T}=[0,\ \timeh],$
where $\timeh$ is the time  horizon of the model. We also write $E$ for $E^{s}=E^{s}[0,\infty)[ \,)$
and $E'$ for $E^{-s}[0,\infty)[ \,).$
\subsection{Basic definitions.}
We recall that, by the bilinear form $\sesq{\;}{\;},$ the space  $E^{-s}$ is identified 
with the dual of $E^{s},$ that is, the space of continuous linear functionals on $E^{s}.$
It is important to note that, since
$s>1/2$, the space $E^{s}$ is contained in $C_{b}^{0}([0,\infty)[ \,),$ the space of bounded continuous functions on $[0,\infty [ \,$, so
that $E^{-s}$ contains the dual of $C_{b}^{0}(  [0,\infty[ \,),$
which is the space of bounded Radon measure on $[0,\infty)$. In particular,
all Dirac masses $\delta_{x}$, for $x\geq0$, belong to $E^{-s}$.
\begin{definition} \label{port def}
A portfolio is  progressively measurable process on the time interval $\mathbb{T},$ with values in
$E^{-s}.$ If $\theta$ is a  portfolio, then its discounted value  at time $t \in \mathbb{T}$ is
\begin{equation}
\prtfpxd{t}(\theta)=\sesq{\theta_{t}}{ \zcpxd{t}}.
  \label{wealth *}
\end{equation}
\end{definition}
The basic example is a  portfolio of one zero-coupon:
\begin{example}  \label{Ex ZC 1} \text{} \normalfont  \\
Consider a portfolio containing exactly one zero-coupon bond
with maturity date $T,$ i.e. \textit{time of maturity} $T:$ \\
1) Let $T \geq \timeh$ and let $T$ be fixed.
The portfolio $\theta$ is then defined by
\begin{equation} \label{eq Ex ZC 1}
\theta_{t}=\delta_{T-t}, \; \forall t\leq \timeh.
\end{equation}
Since $T \geq \timeh,$ we have indeed that the support of the distribution $\theta_{t}$
is contained in $[0,\infty[ \,,$ so $\theta_{t} \in E^{-s}.$
With this definition, the value of the zero-coupon is:%
\[
<\delta_{T-t},p_{t}>\ =\ p_{t}\left(  T-t\right)
\]
which is precisely what we had in mind. \\
2) Let $T < \timeh$ and let $T$ fixed. In this case we note that the process in (\ref{eq Ex ZC 1}) does not
continue after time $T$: the zero-coupon is converted into cash. So the
buy-and-hold strategy is not possible for zero-coupon bonds, unless 
the horizon $\timeh$ is less than  the maturity $T.$ \\
3) Let $T=t+x$ and $x \geq 0$ a fixed \textit{time to maturity}. Then the portfolio
is defined by
\begin{equation} \label{eq Ex ZC 2}
\theta_{t}=\delta_{x},\text{ \ \ \ for }t\leq \timeh.
\end{equation}
\end{example}

We note that the higher we choose $s$, the more portfolios
can be incorporated into the model. For instance, if $s>3/2$, all curves in
$E^{s}$ are $C^{1}$, so that the derivative $\delta_{x}^{\prime}$ of the Dirac
mass belongs to $E^{-s}$. The value of $\delta_{T-t}^{\prime}$ is:%
\begin{equation}
<\delta_{T-t}^{\prime},p_{t}>\ =p_{t}^{\prime}(T-t)
=-\fwrate{t}(T-t)  p_{t}(T-t),  \label{65}
\end{equation}
where $p_{t}^{\prime}(x) =\partial p_{t}(x)/\partial x$ and where $\fwrate{t}(x),$
defined in (\ref{intrest}), is the instantaneous
forward rate with time to maturity $x,$ contracted at time $t.$
This also implies that the higher we choose $s$, the more interest
rates derivatives can be incorporated into the model.
If $s >1/2,$ then we can contract directly on the values of zero-coupon bond prices,
and if $s >3/2,$  then we can contract directly on the values of interest rates.

We next introduce the notion of \emph{self-financing} portfolio. We state a definition
such that it will makes sense for mild solutions of the bonds dynamics:
\begin{definition} \label{self-fin def}
A portfolio is called \emph{self-financing} if, for every $t \in \mathbb{T}$
\begin{equation} \label{self-fin prtf}
V_{t}(\theta) =V_{0}(\theta) +\int_{0}^{t} <\theta_{s}\,,\,p_{s}m_{s} \ ds+\sum_{i\in
\mathbb{I}}p_{s}\sigma_{s}^{i}dW_{s}^{i}>.
\end{equation}
\end{definition}
Given a strong solution $p$ of the bonds dynamics, we have for a self-financing portfolio:
\begin{equation}
dV_{t}(\theta)=\ <\theta_{t}\,,\,dp_{t}-\partial p_{t} \ dt>. \label{b5}%
\end{equation}
Note that this is not the standard definition:\ this is because we are in the
moving frame. Changes in portfolio value are due to two causes: changes in
prices, as in the fixed frame, and also to changes in time to maturity.

For the right-hand side of (\ref{self-fin prtf}) to make mathematical sense and to
introduce later arbitrage free markets, we need a further definition.

\begin{definition} \label{def self-fin}
A portfolio $\theta$ is an admissible portfolio if $\| \theta \|_{\prtfs} < \infty,$
where
\[
\| \theta \|^{2}_{\prtfs}=
E\left[  (\int_{0}^{\timeh}|<\theta_{t}\,,\,p_{t}m_{t}>|dt)^{2}+\int
_{0}^{\timeh}\sum_{i\in\mathbb{I}}(<\theta_{t}\,,\,p_{t}\sigma_{t}^{i}%
>)^{2}dt\right].
\]
$\prtfs$ is the linear space of all admissible portfolios and $\sfprtfs$ the subspace
of self-financing portfolios.
\end{definition}
The discounted gains process $G,$ defined by
\begin{equation}
G(t,\theta)=\int_{0}^{t}(\sesq{\theta_{s}}{\zcpxd{s}\drift{s}}ds
       +\sesq{\theta_{s}}{\zcpxd{s}\vol{}{s}d\wienerp{}{s}}) ,
   \label{Gain * explicit}
\end{equation}
is well-defined for admissible portfolios:
\begin{proposition}\label{proposition; square integ gain*}
Assume that  $\zcpxd{0},$ $\drift{}$ and $\vol{}{}$ are as in Proposition \ref{prop. boundary cond}.
If $\theta \in \prtfs,$
then $G( \cdot,\theta)$ is continuous a.s. and 
$E(\sup_{t \in \mathbb{T}}(G( t,\theta))^{2}) < \infty.$
\end{proposition}
\begin{proof}
Let $\theta \in \prtfs$ and introduce $X=\sup_{t \in \mathbb{T}}|G(t,\theta)|,$ $Y(t)=\int_{0}^{t}\sesq{\theta_{s}}{\zcpxd{s}\drift{s}}ds$
and $Z(t)=\int_{0}^{t}\sesq{\theta_{s}}{\zcpxd{s}\vol{}{s} d\wienerp{}{s}}.$
Then $G(t,\theta)=Y(t)+Z(t),$ according to formula (\ref{Gain * explicit}).
Let $ \zcpxd{}$ be  given by Proposition \ref{prop. boundary cond}, of which the
hypotheses are satisfied.

We shall give  estimates for $Y$ and $Z.$ By the definition  of $\prtfs:$
\begin{equation} \begin{split}
 E((\sup_{t \in \mathbb{T}} (Y(t))^{2}) 
  \leq E((\int_{0}^{\timeh} |\sesq{\theta_{s}}{\zcpxd{s}\drift{s}} | ds)^{2}) \leq \|\theta\|^{2}_{ \prtfs}.
   \end{split} \label{wealth; self fin strat proof 1}
\end{equation}
By isometry we obtain
\begin{equation} \begin{split}
E(Z(t)^{2})
=&E(\int_{0}^{t}\sesq{\theta_{s}}{ \zcpxd{s} \sum_{i \in \mathbb{I}}\vol{i}{s}  d\wienerp{i}{s}})^{2} \\
&=E(\int_{0}^{t}\sum_{i \in \mathbb{I}}(\sesq{\theta_{s}}{ \zcpxd{s} \vol{i}{s}})^{2}ds) %
  \leq \|\theta\|^{2}_{ \prtfs}.
   \end{split}
\label{wealth; self fin strat proof 2}
\end{equation}
Doob's $L^{2}$ inequality and inequality (\ref{wealth; self fin strat proof 2})
give $E(\sup_{t \in \mathbb{T}}Z(t)^{2}) \leq 4 \|\theta\|^{2}_{ \prtfs}.$
Inequality (\ref{wealth; self fin strat proof 1}) then gives $E(X^{2}) \leq 10 \|\theta\|^{2}_{ \prtfs},$
which proves the proposition.
%\AAAqed
\end{proof}
\begin{example} \label{Ex ZC 2} \text{} \normalfont  \\
1)
The portfolio in 1) of Example  \ref{Ex ZC 2}  is self-financing and  the portfolios
in 2) and 3) of Example  \ref{Ex ZC 2}  are not self-financing. \\
2) The interest rate portfolio in formula (\ref{65}) is self-financing.
 \\
\end{example} 
\subsection{Rollovers}
\begin{definition}
 Let $S\geq0.$ A $S$-rollover is a self-financing portfolio $\theta$ of a number of zero-coupon bonds
with constant time to maturity $S$ and with initial price $V_{0}(\theta)=\zcpxd{0}(S).$
\end{definition}
It follows directly  from the definition that a $S$-rollover have the same initial price
as a zero-coupon with maturity date $S.$ It also follows that, if $x_{t}$ is the number
of zero-coupon bonds in the portfolio at $t,$ then we must have:%
\[
\theta_{t}=x_{t}\delta_{S},
\]
where the real-valued process  $x$ makes the portfolio self-financing.
\begin{proposition}
If $\theta_{t}$ is a $S$- rollover, then:
\begin{equation} \label{r-o eq}
x_{t}= \exp(\int_{0}^{t} \fwrate{s}(S) \ ds).
\end{equation}
\end{proposition}
\begin{proof}
The portfolio $\theta_{t}$ only contains zero-coupons with time to maturity
$S,$ so that $V_{t}(\theta)  =x_{t}p_{t}(S).$ Assuming the process $x$ to be of
bounded variation it follows that:
\[
dV_{t}(  \theta)  =p_{t}(  S)  dx_{t}+x_{t}dp_{t}(S).
\]
Substituting the expression for $dp_{t}(S)$ this becomes:%
\begin{align*}
dV_{t}(  \theta)   &  =p_{t}(  S)  \frac{dx_{t}}%
{dt}dt+x_{t}\partial_{x}p_{t}(  S)  dt+x_{t}p_{t}(  S)
(  m_{t}(  S)  dt+\sum_{i\in\mathbb{I}}\sigma_{t}^{i}(
S)  dW_{t}^{i}) \\
&  =(  p_{t}(  S)  \frac{dx_{t}}{dt}+x_{t}\partial_{x}%
p_{t}(  S)  +x_{t}p_{t}(  S)  m_{t}(  S)
)  dt+x_{t}p_{t}(  S)  \sum_{i\in\mathbb{I}}\sigma_{t}%
^{i}(  S)  dW_{t}^{i}.
\end{align*}
According to (\ref{self-fin prtf}) the portfolio is then self-financing if and only if:%
\[
p_{t}(  S)  \frac{dx_{t}}{dt}+x_{t}(\partial p_{t})(S)  =0.
\]
This means that:%
\[
\frac{1}{x_{t}}\frac{dx_{t}}{dt}=-\frac{1}{p_{t}(  S)  }%
\frac{\partial p_{t}(S)}{\partial S}=\fwrate{t}(S).
\]
and the formula (\ref{r-o eq}) follows by integration. This proves the proposition
since $x$ then is of bounded variation.
%\AAAqed
\end{proof}

\noindent In particular, if $S=0,$ then we get the usual bank account with spot rate
$r_{t}.$

Henceforth, we will denote by $q_{t}(S)  $ the value (discounted
to $t=0$) at time $t$ of a $S$-rollover. In the preceding notation, $q_{t}(S)=V_{t}(\theta).$

Introducing the price  curve of the roll-over at time $t,$ 
$q_{t}:[0,\infty [ \, \rightarrow \mathbb{R}$, we find that the price
dynamics of roll-overs is given by:
\begin{equation}
q_{t}=p_{0}+\int_{0}^{t}q_{s}m_{s}ds+\int_{0}^{t}q_{s}\sum_{i\in\mathbb{I}%
}\sigma_{s}^{i}dW_{s}^{i}, \label{SPDE q}%
\end{equation}
Note that, compared to the same formula for bond prices, the term in
$\partial$ has disappeared from the right-hand side.

A $S$-rollover is a bank account which needs advance notice to be cashed: if
notice is given at time $t$, the rollover will then pay $x_{t}$ units of
account at time $t+S.$ In other words, at time $t,$ when notice is given, the
rollover is exchanged for $q_{t}(S)/p_{t}(S)=x_{t}$ units of a unit
zero-coupon with time of maturity $t+S.$

As we noted earlier, zero-coupons do not in general allow buy-and-hold
strategies. However rollovers do: a constant portfolio of rollovers is always
self-financing. A general bond portfolio $\theta_{t}$ can be expressed in
terms of a portfolio of rollovers $\eta_{t}$ and vice versa.
\subsection{Absence of arbitrage opportunities.}
Let $p$ be a mild solution of the price dynamics. Suppose that $\theta_{t}$ is a self-financing portfolio such that, for almost
every $(t,\omega) \in \mathbb{T} \times \Omega,$ we have:%
\begin{equation}
\forall i\in\mathbb{I},\text{ \ }<\theta_{t}\,(  \omega)
,\,p_{t}(  \omega)  \sigma_{t}^{i}(  \omega)  >=0.
\label{70}%
\end{equation}
(We note that $p_{t}(\omega) \in E^{s}$ is a function of time to maturity,
$x \mapsto p_{t}(\omega,x),$ and similarly
for $\theta_{t}$ etc.)
Then (\ref{self-fin prtf}) gives $dV_{t}(\theta)=<\theta_{t}\,,\,m_{t}p_{t}>dt$, so that
$\theta_{t}$ is risk-free. Since the spot rate is zero (after discounting
values to $t=0$), in an arbitrage free market it must follow that for almost
every $(  t,\omega)  $:
\begin{equation}
<\theta_{t}(\omega)\,,\,m_{t}(\omega)p_{t}(\omega)>=0. \label{71}
\end{equation}
Comparing (\ref{70}) and (\ref{71}), we find that $p_{t}(\omega)m_{t}(\omega)$
must belong to the closure of the linear span of
$\{p_{t}(\omega)\sigma_{t}^{i}(\omega)\,|\,i\in\mathbb{I}\}.$
In fact this follows rigorously using Lemma \ref{lm hedge eq}, proved independently of this subsection.
There are now two cases:
\begin{itemize}
\item $\mathbb{I}$ is finite. %
Then the linear span is finite-dimensional, and
it coincides with its closure. So there are numbers $\gamma_{t}^{i}(
\omega)  ,i\in\mathbb{I}$ such that
\[ p_{t}(\omega)m_{t}(\omega)
=p_{t}(\omega)\sum_{i\in\mathbb{I}}\gamma_{t}^{i}(\omega)\sigma_{t}^{i}(\omega)
  \;\; \text{ (finite sum)}.\]
 Since $p_{t}(  \omega)  >0$
 for almost every $(  t,\omega)  $, this leads to:%
\[
m_{t}(\omega)=\sum_{i\in\mathbb{I}}\gamma_{t}^{i}(\omega)\sigma_{t}^{i}%
(\omega)
\]
and since the processes $m$ and $\sigma^{i}$ are progressively measurable, so can one choose
the processes $\gamma^{i}$. Note that the preceding equation holds in
$E^{s}$, and that it translates into a family of equations in $[0, \infty [ \,:$%
\[
m_{t}(\omega,x)=\sum_{i\in\mathbb{I}}\gamma_{t}^{i}(\omega)\sigma_{t}%
^{i}(\omega,x)\ \ \ \forall x\geq0
\]
or, as usual, omitting to mention the $\omega$ variable:
\[
m_{t}(x)=\sum_{i\in\mathbb{I}}\gamma_{t}^{i}\sigma_{t}^{i}(x)\ \ \ \forall
x\geq0.
\]
The $\gamma_{t}^{i}$ are the components of a market price of risk, and they do not depend on
the time to maturity $x.$ Using the volatility operator process $\sigma$ the last equality reads
\begin{equation} \label{mkt price 1}
m_{t}=\sigma_{t} \gamma_{t} \; \; \forall t  \in \mathbb{T}
\end{equation}
and any $\gamma,$ progressively measurable with values in
$\ell^{2}(\mathbb{I}),$ satisfying this equation is called a  market price of risk process.
\item $\mathbb{I}=\mathbb{N}.$ Then the linear span is not closed in general;
in fact, it is closed if and only if it is finite-dimensional. In that case,
we shall impose a stronger condition. To prove that the market is arbitrage-free,
we shall use that $m_{t}(\omega)$ is in the range of the volatility operator
$\sigma_{t}(\omega)$ which is a subset of the above closed linear span. So, once more
we impose that the condition (\ref{mkt price 1}) should be satisfied,
but for $\gamma$ with values in  $\ell^{2}(\mathbb{I}).$ If the range  of $\sigma_{t}(\omega)$
is infinite dimensional, then
this condition is indeed stronger, since $\sigma_{t}(\omega)$ is a.e. a compact operator.
\end{itemize}
In both cases, we also need that  $\gamma$ satisfy some integrability condition
in $(\omega,t).$
This leads us to the following
\begin{definition} \label{market cond}
We shall say that the market is strongly arbitrage-free
if there exists a progressively measurable process $\gamma$ with values
in $\ell^{2}(\mathbb{I}),$
such that
\begin{equation}
m_{t}=\sigma_{t} \gamma_{t}, \; \;\forall t \in \mathbb{T}
\label{rel drift-sigma}%
\end{equation}
and
\begin{equation}
E\left[  \exp(a\int_{0}^{\timeh}\norm{\gamma_{t}}^{2}dt)\right]  <\infty,\quad\forall a\geq0.
\label{gamma strong}%
\end{equation}
\end{definition}
If the market is strongly arbitrage-free then, by the Girsanov theorem, a
martingale measure is given by $dQ=\xi_{\timeh}dP$, with:%
\begin{equation}
\xi_{t}=\exp\left(  -\frac{1}{2}\int_{0}^{t} \norm{\gamma_{s}}^{2}ds-\sum_{i\in\mathbb{I}}\gamma_{s}^{i}dW_{s}
^{i}\right).  \label{75}%
\end{equation}
The $\wienerq{i}{},$ $i\in \mathbb{I},$ where
\begin{equation} \label{W Q}
\tilde{W}_{t}^{i}=W_{t}^{i}+\int_{0}^{t}\gamma_{s}^{i}ds,
\end{equation}
are independent Wiener process with respect to $Q.$ The expected
value of a random variable $X$ with respect to $Q$ is given by:
\[
E_{Q}[X]=E[\xi_{\timeh}X].
\]
Under a martingale measure, the discounted zero-coupon price process $p$
satisfies the equation
\begin{equation} \label{SPDE p Q 1}
\ p_{t}=\mathcal{L}_{t}p_{0}
+\int_{0}^{t}\mathcal{L}_{t-s}(p_{s}\sigma_{s})d\wienerq{}{s}
\end{equation}
and also the equation
\begin{equation}
p_{t}=p_{0}+\int_{0}^{t}\partial p_{s}ds+\int_{0}^{t}p_{s}\sigma_{s}\,d\wienerq{}{s}.
\label{SPDE p Q 2}%
\end{equation}
The discounted roll-over price process $q_{t}$ is given by:
\begin{equation}
q_{t}=p_{0}+\int_{0}^{t}q_{s} \sigma_{s}d\tilde{W}_{s}. \label{SPDE q Q}%
\end{equation}
\begin{lemma} \label{lm self-fin}
A portfolio $\theta$ is self-financing if and only if:
\begin{equation}
V_{t}(\theta)=V_{0}(\theta)+\int_{0}^{t}\sum_{i\in\mathbb{I}}<\theta
_{t}\,,\,p_{t}\sigma_{t}^{i}>d\tilde{W}_{s}^{i}. \label{self-fin}%
\end{equation}
\end{lemma}
We note that the integrand is in fact the adjoint operator of the operator
$b_{t}(\omega)= p_{t}(\omega)\sigma_{t}(\omega)$
from $\ell^{2}(\mathbb{I})$ to $E^{s}([0,\infty[ \,):$
\begin{equation} \label{self-fin1} %
(b_{t}(\omega)'\theta_{t})^{i}=<\theta_{t}\,,\,p_{t}\sigma_{t}^{i}>, \; \; \forall i \in \mathbb{I}.
\end{equation}
To see this, with $x_{t}^{i}(\omega)=<\theta_{t}\,,\,p_{t}\sigma_{t}^{i}>,$
rewrite it as follows: \\
for all $(t,\omega)$ and all $z \in \ell^{2}(\mathbb{I})$
\begin{align*}
\left(  z,x_{t}(\omega)\right)_{\ell^{2}}   &  =\sum_{i\in \mathbb{I}}z^{i}<\theta_{t}\left(
\omega\right)  ,\ p_{t}\left(  \omega\right)  \,\sigma_{t}^{i}\left(
\omega\right)  >\\
&  =<\theta_{t}\left(  \omega\right)  ,\ p_{t}\left(  \omega\right)
\sum_{i\in \mathbb{I}}\sigma_{t}^{i}\left(  \omega\right)z^{i}  \,>\\
&  =<\theta_{t}\left(  \omega\right)  ,\ b_{t}\left(  \omega\right)  z\,>
 =<b_{t}\left(  \omega\right)  ^{^{\prime}}\theta_{t}\left(  \omega\right)
,\ z\,>.
\end{align*}

If the market is strongly arbitrage-free and if condition (\ref{domain sigma i})
of Theorem \ref{Th price compatible strong cond.} is satisfied, then also
condition  (\ref{domain drift}) is satisfied and the
Theorem \ref{Th price compatible strong cond.} applies.

\section{Hedging of interest derivatives} \label{sect. hedging}
From now on, it will be a standing assumption that $p_{0}$ satisfies condition (\ref{p0}),
that $\vol{}{}$ satisfy conditions (\ref{63}) and (\ref{domain sigma i}) and that
the market is strongly arbitrage-free according to Definition \ref{market cond}.

Before we solve the optimal portfolio problem, we shall study the problem of
hedging a European interest rates derivative with payoff $X$ at maturity $\timeh.$
$X$ is said to be an attainable contingent claim or derivative if $V_{\timeh}(\theta)=X$
for some admissible self-financing portfolio $\theta.$
Here we are only interested in payoffs, relevant for the  optimal portfolio problem
considered in these notes, i.e. $X \in L^{p}(\Omega,\mathcal{F},P)$ for every $p \geq 1$
(see Lemma \ref{X in Lp}).
We first introduce the hedging equation, the Malliavin derivative and the
Clark-Ocone representation formula, which then permits the reader, if he wish, to
proceed directly to the study of the optimization problem in the case of deterministic
$\sigma$ and $\gamma$ in \S \ref{determ case}

Assume that $X\in L^{2}(\Omega,\mathcal{F},Q),$ where $Q$ is one equivalent martingale
measure given by (\ref{75}). Then, by the martingale
representation theorem, $X$ can be written as a stochastic integral:
\begin{equation}  \label{mart decomp 1}
X=E_{Q}[X]+\int_{0}^{\timeh}\sum_{i\in\mathbb{I}}
x_{t}^{i}d\tilde{W}_{t}^{i},
\end{equation}
with:
\begin{equation} \label{mart decomp 2}
E_{Q}[\int_{0}^{\timeh} \norm{x_{t}}^{2}dt] <\infty.
\end{equation}
Comparing with equations (\ref{self-fin}) and (\ref{self-fin1}) for a self-financing portfolio, we
obtain the hedging equation
\begin{equation} \label{hedge eq 1}
b_{t}(\omega)'\theta_{t}(\omega)=x_{t}(\omega), \; \text{a.e.} \; (t,\omega),
\end{equation}
where the operator $b_{t}(\omega)=\zcpxd{t}(\omega)\vol{}{t}(\omega)$
from $\ell^{2}(\mathbb{I})$ to $E^{s}([0,\infty[ \,)$ was introduced in (\ref{self-fin1}).
Equivalently: for almost every$\;(t,\omega),$
\begin{equation*}
x_{t}^{i}(\omega)=\ <\theta_{t}\left(  \omega\right)  ,\ p_{t}\left(
\omega\right)  \,\sigma_{t}^{i}\left(  \omega\right)  >,\;\forall
i\in\mathbb{I}\text{ }. %
\end{equation*}

We next introduce the Malliavin derivative (c.f. \cite{Nual71}), $D_{t}X,$ with respect
to $\wienerq{}{},$ at time $ t \in \mathbb{T}$ of certain
$\mathcal{F}=\mathcal{F}_{\timeh}$ measurable real random variables $X$ by: \\

\noindent
D1) $D_{t}X=0,$ if $X$ is a constant, \\
D2) $D_{t}X=h_{t},$ if  $h \in L^{2}(\mathbb{T}, \ell^{2}(\mathbb{I}))$ and
  $X=\int_{0}^{\timeh}\sum_{i\in\mathbb{I}}h_{t}^{i}d\tilde{W}_{t}^{i},$ \\
D3) $D_{t}(X Y)=X D_{t}Y+ Y D_{t} X.$ \\

\noindent
The algebra of such random variables is dense in $ L^{2}(\Omega,\mathcal{F},Q),$
which can be used to extend the definition to larger sets.
$D_{t}X$ takes its values in $\ell^{2}(\mathbb{I}).$
The partial derivative, with respect to $\wienerq{i}{},$ $D_{i,t}X,$ is the $i$-th component of
$D_{t}X.$

We will use the following
expression for the Malliavin derivative of an It\^o stochastic integral:
\begin{equation} \label{malliavin 1}
D_{t}\int_{0}^{\timeh}\sum_{i\in\mathbb{I}}x_{s}^{i}d\tilde{W}_{s}^{i}
=x_{t}+\int_{t}^{\timeh}\sum_{i\in\mathbb{I}}(D_{t}x_{s}^{i})d\tilde{W}_{s}^{i},
\end{equation}
when almost all the $x_{s}^{i}$ are Malliavin differentiable and sufficiently integrable.

In the case when $X$ is Malliavin differentiable, the
Clark-Ocone representation formula states that the integrand $x_{t}$ in
(\ref{mart decomp 1}) is given by
\begin{equation}
x_{t}=E_{Q}\left[  D_{t}X\ |\ \mathcal{F}_{t}\right].
\label{Clark-Ocone}%
\end{equation}

We now come back to the hedging equation  (\ref{hedge eq 1}).
The fact that $\theta_{t}=\delta_{0}$ is a solution to the homogeneous %
equation (\ref{hedge eq 1}) permits us to construct  self-financed solutions of
the in-homogeneous equation (\ref{hedge eq 1}), from solutions, which are not self-financed:
\begin{lemma} \label{lm hedge eq}
If $\bar{\theta}$ is an admissible portfolio (not necessarily
self-financed) which satisfies (\ref{hedge eq 1}), then there is a unique
self-financing admissible portfolio $\theta_{t}$ such that the difference
$\theta_{t}-\bar{\theta}_{t}$ is risk-free. It is given by:
\begin{equation} \label{sol hedge eq}
\theta_{t}  =a_{t}\delta_{0}+\bar{\theta}_{t},
\end{equation}
\begin{equation} \label{sol hedge eq1}
a_{t}  =\frac{1}{p_{t}(0)}\left[  E_{Q}[X\,|\,\mathcal{F}_{t}]-V_{t}%
(\bar{\theta})\right].
\end{equation}
\end{lemma}
\begin{proof}
We here omit the argument $\omega.$ Since the portfolio $\theta_{t}%
-\bar{\theta}_{t}$ is risk-free, it must have time to maturity $0$, and the %
formula (\ref{sol hedge eq}) is true by definition.
Substituting \ into equation (\ref{hedge eq 1}), and bearing in mind that
$\sigma_{t}^{i}(0)=0$:
\begin{align*}
((p_{t}\sigma_{t})'&\theta_{t})^{i}=
 <\theta_{t},\ p_{t}\,\sigma_{t}^{i}\;>\ =\ <a_{t}\delta_{0}+\bar{\theta
}_{t},\ p_{t}\,\sigma_{t}^{i}>\text{ }\\
&  =\ a_{t}\ p_{t}(0) \ \sigma_{t}^{i}(0)+<\bar{\theta}_{t},\ p_{t}\,\sigma_{t}^{i}>\\
&  =x_{t}^{i}~\ \ \ \forall i\in\mathbb{I}.
\end{align*}
So $\theta_{t}$ satisfies (\ref{hedge eq 1}). It is then  a hedging portfolio of $X$
if $V_{t}(\theta)=E_{Q}[X\,|\,\mathcal{F}_{t}].$ Substituting again
(\ref{sol hedge eq}) and then (\ref{sol hedge eq1}), we get:
\[
V_{t}(\theta)=a_{t}V_{t}(\delta_{0})+V_{t}(\bar{\theta})
 =a_{t}p_{t}(0)+V_{t}(\bar{\theta})
= E_{Q}[X\,|\,\mathcal{F}_{t}].
\]
If $\bar{\theta}$ is an admissible portfolio, then $\theta$ is also admissible,
since  $\| \theta \|_{\prtfs}=\| \bar{\theta} \|_{\prtfs}.$
%\AAAqed
\end{proof}

By the lemma, the construction of a hedging portfolio for $X$ is reduced to
solve equation (\ref{hedge eq 1}) in $\theta_{t}(  \omega )$
for every $\left(  t,\omega\right),$ in  such a way that $\theta \in \prtfs,$
i.e. $\theta$ is admissible. Any such solution $\theta$ of this equation contains
the risky part of the portfolio.

To solve equation (\ref{hedge eq 1}), for given $(t,\omega ),$ we have to know if
$x_{t}(\omega)$ is in the range of the operator $b_{t}(\omega)'.$
The closure of the range of $b_{t}(\omega)'$ is equal to
the orthogonal complement $(\mathcal{K}(b_{t}(\omega)))^{\perp}$
of the kernel $\mathcal{K}(b_{t}(\omega))$ of $b_{t}(\omega).$

Consider the cases of $\mathbb{I}$ finite:
The  range $\mathcal{R}((b_{t}(\omega))')$ is then closed, since it is finite dimensional.
The kernel $\mathcal{K}(b_{t}(\omega))$ is trivial
 iff the $p_{t}\left(  \omega\right)  \,\sigma_{t}^{i}\left(  \omega\right)$
are linearly independent.
So $(b_{t}(\omega))'$ is
surjective and and there is a (non-unique) solution $\theta_{t}(  \omega ),$ for every $x_{t}(\omega),$
 iff the $p_{t}\left(  \omega\right)  \,\sigma_{t}^{i}\left(  \omega\right)$
are linearly independent.

Consider the cases of $\mathbb{I}$ infinite:
The map $(b_{t}(\omega))'$ from  $E^{-s}([0,\infty[ \,)$ to $\ell^{2}(\mathbb{I}),$
is then never surjective. In fact, $b_{t}(\omega)$ is a Hilbert-Schmidt operator,
so it is compact. The adjoint is then also compact and since $\ell^{2}(\mathbb{I})$
is infinite dimensional, its range must be a proper subspace of  $\ell^{2}(\mathbb{I}).$
This is the  basic reason why there are always non-attainable contingent claims,
when $\mathbb{I}$ is infinite.

We have the following result
(see Th.4.1 and Th.4.2 of \cite{E.T Bond  Completeness} for the case $\mathbb{I}=\mathbb{N}$):
\begin{theorem} \label{Th D0 non complete and approx complete}
Let $\derprod{}{0}=\cap_{p \geq 1} L^{p}(\Omega,P,\mathcal{F}).$ \\
$i)$ If $\mathbb{I}=\mathbb{N},$ then there exists $X \in \derprod{}{0}$
such that $\prtfpxd{\timeh}(\theta) \neq X$ for all $\theta \in \sfprtfs.$ \\
$ii)$ $\derprod{}{0}$ has a dense subspace of attainable contingent claims
if and only if the operator $\vol{}{t}(\omega)$ has
a trivial kernel a.e. $(t,\omega) \in \mathbb{T} \times \Omega.$
\end{theorem}
Statement $ii)$  says by definition that the bond market is approximately complete
(notion introduced in \cite{Bj-Ka-Ru97} and  \cite{Bj-Ma-Ka-Ru97})
if and only if $\vol{}{t}(\omega)$ has a trivial kernel a.e. 

In the sequel of this section, we are interested in the hedging problem for
approximately complete markets, so
 we only consider the solution of the hedging equation (\ref{hedge eq 1})
in the case when $\vol{}{t}(\omega)$ has a trivial kernel a.e.
$(t,\omega) \in \mathbb{T} \times \Omega.$

Consider the case when $\mathbb{I}=\mathbb{N}$ is an infinite and let  $\ell^{2}= \ell^{2}(\mathbb{I}).$
To derive a condition under which (\ref{hedge eq 1}) has a solution
and to derive a closed formula for one of the solutions,
we rewrite the l.h.s. of (\ref{hedge eq 1}) using the notations
\begin{equation}   \label{B and l}
l_{t}=\ltrans{t}\zcpxd{0}, \;  B_{t}(\omega)=l_{t} \vol{}{t}(\omega)
 \; \text{and } \;   \eta_{t}(\omega)=\sequiv^{-1}(\zcpxd{t}(\omega)/l_{t})\theta_{t}(\omega).
\end{equation}
Then
\begin{equation*}
\begin{split}
(\vol{}{t}(\omega))'&\zcpxd{t}(\omega)\theta_{t}(\omega)
=(\vol{}{t}(\omega))'l_{t}(\zcpxd{t}(\omega)/l_{t})\theta_{t}(\omega)
=(l_{t} \vol{}{t}(\omega))'(\zcpxd{t}(\omega)/l_{t})\theta_{t}(\omega) \\
&=(l_{t} \vol{}{t}(\omega))^{*}\sequiv^{-1}(\zcpxd{t}(\omega)/l_{t})\theta_{t}(\omega)
=(B_{t}(\omega))^{*}\eta_{t}(\omega).
\end{split}
\end{equation*}
The linear operator $B_{t}(\omega)$ is given, since $\zcpxd{0}$
and $\vol{}{t}(\omega)$ are supposed given.
Applying  Theorem \ref{Th price compatible strong cond.} to the factor $\zcpxd{}/l,$
it follows that equation (\ref{hedge eq 1}) is equivalent to find a progressive $E^{s}$-valued process
$\eta$ satisfying
the equation
\begin{equation}   \label{prtf eq 2}
  (B_{t}(\omega))^{*}\eta_{t}(\omega)=x_{t}(\omega),  \; \text{a.e. }  (t,\omega) \in \mathbb{T} \times \Omega.
\end{equation}
We define the self-adjoint operator $A_{t}(\omega)$ in $\ell^{2}$ by
\begin{equation}
A_{t}(\omega)=(B_{t}(\omega))^{*}B_{t}(\omega).
  \label{At}
\end{equation}
It is a fact of basic Hilbert space operator theory (cf. \cite{Kato66}) that the range
$\mathcal{R}((B_{t}(\omega))^{*})$ $=\mathcal{R}((A_{t}(\omega))^{1/2}).$
The solvability of each one of equations (\ref{hedge eq 1}) and (\ref{prtf eq 2})
is therefore equivalent to the existence of a progressive $\ell^{2}$-valued process
$z$ satisfying
\begin{equation}  \label{hedge eq l2}
  (A_{t}(\omega))^{1/2}z_{t}(\omega)=x_{t}(\omega),  \; \text{a.e. }  (t,\omega) \in \mathbb{T} \times \Omega.
\end{equation}
The kernel $\mathcal{K}((A_{t}(\omega))^{1/2})$ is trivial since
$\mathcal{K}((A_{t}(\omega))^{1/2})=\mathcal{K}(A_{t}(\omega))=\mathcal{K}(B_{t}(\omega))$ $=\{0\}.$
Now, if $x_{t}(\omega) \in \mathcal{R}((B_{t}(\omega))^{*})$ then the unique solution of (\ref{hedge eq l2})
is $z_{t}(\omega)=(((A_{t}(\omega))^{1/2})^{-1}x_{t}(\omega)$ and a solution of (\ref{prtf eq 2})
is given by
\begin{equation}
 \eta_{t}(\omega)=S_{t}(\omega)(A_{t}(\omega))^{-1/2}x_{t}(\omega),
  \label{prtf eq 3}
\end{equation}
where $S_{t}(\omega),$ the closure of the operator $B_{t}(\omega)(A_{t}(\omega))^{-1/2},$
is isometric (cf. \cite{Kato66}) from $\ell^{2}$ to $E^{s}.$
Let $a$ be  as in (\ref{sol hedge eq1}) and %
\begin{equation}
\theta =a \delta_{0}+\bar{\theta} \; \text{and} \;  \bar{\theta}_{t}=(l_{t}/\zcpxd{t}) \sequiv \eta_{t}. %
  \label{prtf eq 4}
\end{equation}
Then $\theta $ is  a hedging portfolio according to Lemma \ref{lm hedge eq}.

In order to ensure that $x_{t}(\omega)$ of
(\ref{hedge eq 1}) is in the range of $(\dualvol{}{t}\zcpxd{t})(\omega),$
we introduce spaces $\ell^{s,2},$ of vectors decreasing faster (for $s>0$) than those
of $\ell^{2}.$ For $s \in \mathbb{R},$ let $\ell^{s,2}$ be the  Hilbert space of real sequences
endowed with the norm
\begin{equation}
\|x\|_{\ell^{s,2}}=(\sum_{i \in \mathbb{N}}(1+i^{2})^{s}|x^{i}|^{2})^{1/2}.
    \label{ls,2}
\end{equation}
Obviously $\ell^{2}=\ell^{0,2}$ and $\ell^{s',2} \subset \ell^{s,2},$ if $s' \geq s.$
Although $(A_{t}(\omega))^{-1/2}$ is an unbounded operator in $\ell^{2}$ its restriction
to $\ell^{s,2}$ can be a bounded operator for some sufficient large $s>0,$ i.e.
$(A_{t}(\omega))^{-1/2}\ell^{s,2} \subset \ell^{2}.$ This is the idea of our assumption,
which will ensure hedgeability. However a precise formulation of this assumption
must, as in the case of a finite of Bm., take care of integrability properties in $(t,\omega).$

To consider  also  the case of a finite $\mathbb{I},$
we define  after obvious modifications the operator
$A_{t}(\omega))$ in $\ell^{2}(\mathbb{I})$ by formula (\ref{At}). In this case
$A_{t}(\omega))$ has obviously a bounded inverse.
\begin{condition} \label{uniform cond sigma}
$i)$ If $Card(\mathbb{I}) < \infty,$ then there exists $k \in \derprod{}{0},$
such that for all  $x \in \ell^{2}(\mathbb{I}):$
\begin{equation}
 \|x\|_{\ell^{2}} 
 \leq k(\omega)  \|(A_{t}(\omega))^{1/2}x \|_{\ell^{2}} \;
                           \text{a.e.} \; (t,\omega) \in \mathbb{T} \times \Omega.
  \label{uniform cond sigma eq Rm}
\end{equation}
$ii)$ If $\mathbb{I}=\mathbb{N},$ then there exists $s >0$ and
$k \in \derprod{}{0},$ such that for all  $x \in \ell^{2}(\mathbb{I}):$
\begin{equation}
 \|x\|_{\ell^{2}} 
 \leq k(\omega)  \|(A_{t}(\omega))^{1/2}x \|_{\ell^{s,2}} \; \text{a.e.} \; (t,\omega) \in \mathbb{T} \times \Omega.
  \label{uniform cond sigma eq}
\end{equation}
\end{condition}
In the case of a finite number of Bm. Condition \ref{uniform cond sigma} $i)$ leads
to a complete market and one can choose a hedging portfolio such that it is continuous in
the asset to hedge. To state the result let use introduce the notation
$\derprod{}{0}(F)=\cap_{p \geq 1} L^{p}(\Omega,P,\mathcal{F},F),$ where $F$ is a Banach space.
$\derprod{}{0}=\derprod{}{0}(\mathbb{R}).$
\begin{theorem}[Finite number of random-sources, $Card(\mathbb{I}) < \infty$]
\label{th completeness R^m} %
 \text{} \\
If $(i)$ of Condition \ref{uniform cond sigma} is satisfied and 
if $X \in \derprod{}{0},$ then the portfolio given by equation (\ref{prtf eq 4})
satisfies $\theta \in \sfprtfs$ and  $\prtfpxd{\timeh}(\theta)=X.$
Moreover the linear mapping
$\derprod{}{0} \ni X \mapsto \theta \ \in \prtfs \cap \derprod{}{0}(L^{2}(\mathbb{T},\dualE{})),$
is continuous.
\end{theorem}
\begin{proof}
We only outline the proof of the theorem.
Here $\ell^{2}=\ell^{2}(\mathbb{I})=\mathbb{R}^{\numrand}$ is finite dimensional.  \\
Let $X \in \derprod{}{0}$ and let $x$ be given by (\ref{mart decomp 1}).
First one proves (see Lemma 3.1 of \cite{E.T Bond  Completeness})
 that
\begin{equation} \label{D 1}
\derprod{}{0}(F)=\cap_{p \geq 1} L^{p}(\Omega,Q,\mathcal{F},F).
\end{equation}
Applying the BDG inequalities to equation (\ref{mart decomp 1}) it follows that
\begin{equation} \label{proof R^m 1}
x \in \derprod{}{0}(L^{2}(\mathbb{T},\ell^{2})),
\end{equation}
where $x$ is progressively measurable. The definition of $\eta$ in  (\ref{prtf eq 3})
and the condition (\ref{uniform cond sigma eq Rm}) give
\[
 \|\eta_{t}(\omega)\|_{\ell^{2}}
 \leq k_{t}(\omega) \|x_{t}(\omega) \|_{\ell^{2}}.
\]
Inequality (\ref{proof R^m 1}) then leads to $\eta \in \derprod{}{0}(L^{2}(\mathbb{T},\E{})).$
Using the definition (\ref{prtf eq 4}) of $\bar{\theta}$ we then obtain
\begin{equation} \label{proof R^m 2}
\bar{\theta} \in \derprod{}{0}(L^{2}(\mathbb{T},\dualE{})).
\end{equation}
Since $\bar{\theta}$ satisfies equation (\ref{hedge eq 1})
by construction and since formulas (\ref{proof R^m 1}) and (\ref{proof R^m 2})
shows that  $\bar{\theta}$ is admissible,
the hypotheses of Lemma \ref{lm hedge eq} are satisfied, so $\theta \in \sfprtfs.$
This shows that $\theta$ is a hedging portfolio of $X.$

All the linear maps $X \mapsto x \mapsto \eta \mapsto \theta$ are continuous
in the above spaces, which also proves the claimed continuity of the map $X \mapsto  \theta.$
%\AAAqed
\end{proof}

The solution of the hedging problem, given by Theorem \ref{th completeness R^m},
is highly non-unique, since when $Card(\mathbb{I})=\numrand < \infty$ then
the kernel $\mathcal{K}((\dualvol{}{t}\zcpxd{t})(\omega))$ has infinite dimension.
For instance there is a hedging portfolio $\hat{\vartheta}$ consisting of $\numrand+1$
rollovers at any time. %

To state the result in the case of a infinite number of Bm., we first introduce
spaces of contingent claims $\derprod{}{s},$ smaller than $\derprod{}{0}$ if $s>0$
and corresponding to that the integrand $x$ in (\ref{mart decomp 1}) takes values
in $\ell^{s,2}.$ More precisely, for $s>0$ let
\begin{equation} \label{def D_s}
\derprod{}{s} =\{X \in \derprod{}{0} \; | \; x \in \derprod{}{0}(L^{2}(\mathbb{T},\ell^{s,2}))
  \; \text{where $x$ is given by (\ref{mart decomp 1})} \}.
\end{equation}
Condition \ref{uniform cond sigma} $ii)$ leads to a $\derprod{}{s}$\textit{-complete} market,
i.e.  $\derprod{}{s}$ is a space of attainable contingent claims,
$\derprod{}{s}$ is a dense subspace of $\derprod{}{0}$ and $\derprod{}{s}$ is itself
a complete topological vectorspace. This concept gives a natural frame-work to
study  existence and continuity of hedging portfolios.
We have (see Theorem 4.3 of \cite{E.T Bond  Completeness}):
\begin{theorem}[Infinite number of random-sources $\mathbb{I}=\mathbb{N}$] \text{} \\
\label{th completeness l^2} %
If $(ii)$ of Condition \ref{uniform cond sigma} is satisfied
and if $X \in \derprod{}{s},$ where $s>0$ is given by Condition \ref{uniform cond sigma},
then the portfolio given by equation (\ref{prtf eq 4})
satisfies $\theta \in \sfprtfs$ and  $\prtfpxd{\timeh}(\theta)=X.$
Moreover the linear map
$\derprod{}{s} \ni X \mapsto \theta \in \prtfs \cap \derprod{}{0}(L^{2}(\mathbb{T},\dualE{})),$
 is continuous.
\end{theorem}
For the proof, which only uses elementary  spectral properties of self-adjoint
operators and compact operators, the  reader is referred to \cite{E.T Bond  Completeness}.

A Malliavin-Clark-Ocone formalism was adapted recently in reference \cite{Carmona-Tehr},
for the construction of hedging portfolios in a Markovian context, with a Lipschitz
continuous (in the bond price) volatility operator. This guaranties that the
Malliavin derivative of the bond price is proportional to the volatility operator
(formula (30) of \cite{Carmona-Tehr}).
Hedging is then achieved for a restricted class of claims, namely  European claims
being a Lipschitz continuous function in the price of the bond at maturity.

References \cite{DeDonno Pratelli} and \cite{Pham 2003} studies the hedging
problem in a weaker sense of approximate hedging, which in our context simply boils down
to the well-known existence of the integrand $x$ in the decomposition (\ref{mart decomp 1}).
\section{Optimal portfolio management}
We now consider an investor, characterized by a von-Neumann-Morgenstern
utility function $U$, an initial wealth $v,$ and a horizon $\timeh$. The money is
invested in a market portfolio, and the investor seeks to maximize the
terminal (discounted) value $V_{\timeh}(\theta)$ of the portfolio. Transaction costs and
taxes are neglected. The optimal portfolio problem is then to find an
admissible self-financing portfolio $\hat{\theta}$ with $V_{0}(\hat{\theta})=v,$
such that:
\[
\text{($P_{0}$)}\left\{
\begin{array}
[c]{c}%
\sup E_{P}\left[  U\left(  V_{\timeh}(\theta)\right)  \right]  =E_{P}[U(V_{\timeh}(\hat
{\theta}))]\\
V_{0}(\theta)=v\\
\theta\in\mathsf{P}_{sf}.
\end{array}
\right.
\]
We will follow the now classical two-step approach (cf. \cite{Kr-Scha}, \cite{Pliska86})
 towards solving that problem.
If the portfolio is self-financing and is worth $v$ at time $0$, then, by the
martingale property:%
\[
E_{P}\left[  \xi_{\timeh}V_{\timeh}(\theta)\right]  =v
\]
where the random variable $\xi_{\timeh},$ arising from Girsanov's theorem, was
introduced earlier in (\ref{75}). In general there can be several possible $\xi_{\timeh},$
one for each $\gamma$ satisfying the conditions of Definition \ref{market cond}.
The first step (optimization) consists of finding for given $\gamma,$ among %
$\mathcal{F}_{\timeh}$-measurable random variables $X$ such that $E_{P}\left[
\xi_{\timeh}X\right]  =v$, the one(s) that maximize expected utility $E_{P}%
[U\left(  X\right)  ]$. This problem has in our setting a general solution
$\hat{X},$ given by Proposition \ref{exist unique X}. The second one (accessibility)
consists in hedging one of the contingent claims $\hat{X},$ obtained for the different
$\gamma,$ by a self-financing portfolio
$\hat{\theta}.$ This portfolio is then a solution of the optimal portfolio problem
($P_{0}$). By concavity, the final optimal wealth $V_{\timeh}(\hat{\theta})$ is unique.

\subsection{Optimization}
We consider,
for a given $\gamma$ satisfying the conditions of Definition \ref{market cond},
the optimization problem:%
\[
\left\{
\begin{array}
[c]{c}%
\sup E_{P}\left[  U\left(  X\right)  \right] \\
E_{P}\left[  \xi_{\timeh}X\right]  =v\\
X\in L^{2}\left(  \Omega,\mathcal{F}_{\timeh},P\right)
\end{array}
\right.
\]
We can rewrite it in a more geometric way, involving the scalar product
in $L^{2}\left(  \Omega,\mathcal{F}_{\timeh},P\right)  $:%
\[
\text{(P)}\left\{
\begin{array}
[c]{c}%
\sup\int_{\Omega}U\left(  X\right)  dP\\
\int_{\Omega}\xi_{\timeh}XdP=(\xi_{\timeh},X)_{L^{2}}=v \\
X\in L^{2}\left(  \Omega,\mathcal{F}_{\timeh},P\right)
\end{array}
\right.
\]
Problem (P) consists of maximizing a concave function on a closed linear
subspace of $L^{2}$. Assume there is a maximizer $\hat{X}$. If the usual
theory of Lagrange multipliers applies, there will be some $\lambda\in \mathbb{R}$ such
that $\hat{X}$ actually optimizes the functional
\[
\int_{\Omega}\left[  U\left(  X\right)  -\lambda\xi_{\timeh}X\right]  dP
\]
over all of $L^{2}$. Maximizing pointwise under the integral, and bearing in
mind that $U$ is concave, we are led to the equation:%
\begin{equation}
U^{\prime}\left(  \hat{X}\left(  \omega\right)  \right)  =\lambda\xi
_{\timeh}\left(  \omega\right)  \text{ \ }P\text{-a.e.}, \label{5}%
\end{equation}
which fully characterizes the solution $\hat{X}$. Unfortunately this program
cannot be carried through, for the function $E_{P}\left[  U\left(  X\right)
\right]$ has no point of continuity in $L^{2}$ unless $U$ is bounded, so the
constraint qualification conditions do not hold for problem (P), cf. \cite{I.E.-R.T}.
We will
therefore proceed by a roundabout way:\ use (\ref{5}) to define $\hat{X}$, and
then prove that $\hat{X}$ is optimal for a suitable choice of $\lambda$. For
this, we need some conditions on $U$.
\begin{definition} \label{85}
The utility function $U$ will be called \emph{admissible} if it
satisfies the following properties:
\begin{enumerate}
\item $U:\mathbb{R}\rightarrow\left\{  -\infty\right\}  \cup \mathbb{R}$ is concave and upper semi-continuous

\item there is some $a\in\left\{  -\infty\right\}  \cup \ ]- \infty,0],$ such that $U\left(
x\right)  =-\infty$ if $x<a$ and $U\left(  x\right)  >-\infty$ if $x>a$

\item $U$ is twice differentiable on the interval $A= \ ]a,\ \infty\lbrack$; set
$B=U^{\prime}\left(  A\right)  $

\item $\sup B=+\infty;$  $\inf B =0$ or  $\inf B =- \infty.$

\item $U^{\prime}:A\rightarrow B$ is one-to-one, and there are some
positive constants $r, c_{1}, c_{2}$ and  $c_{3}$
such that its inverse $I=\left[  U^{\prime}\right]  ^{-1}$ satisfies the
estimate $\vert I(  y)\vert + \vert y I'(  y)\vert \leq c_{1}+c_{2}\left\vert
y\right\vert ^{r}+c_{3}\left\vert y\right\vert ^{-r}$ for $y\in B$.
\end{enumerate}
\end{definition}
It follows from these assumptions that $I$ is continuous and strictly
decreasing, with:
\begin{align*}
I\left(  \lambda\right)   &  \rightarrow+\infty\text{ when }\lambda
\rightarrow\inf B\\
I\left(  \lambda\right)   &  \rightarrow a\text{ when }\lambda\rightarrow
+\infty.
\end{align*}
We note that the estimate, in point $5)$ of Definition \ref{85}, is satisfied iff
there exist $C \geq 0$ such that
\[
\vert I(  y)\vert + \vert y I'(  y)\vert
\leq C \ ( |y|^{r}+ |y|^{-r}),
\]
 for all $y\in B.$
All usual utility functions are admissible:
\begin{example} \label{U example} \text{} \\
i) Quadratic utility;
Set $U\left(  x\right)  =\mu x-\frac{1}{2} x^{2},$ $\mu \in \mathbb{R}.$ Then
$a=-\infty$, and $U^{\prime}\left(  x\right)  =\mu - x,$ so that $B=\mathbb{R}$ and
$I\left(  y\right)  = \mu   -y.$ The estimate is satisfied with $r=1.$ \\
ii) Exponential utility; Set $U\left(  x\right)  =1- \frac{1}{\mu} \exp\left(  -\mu x\right),$
$\mu > 0.$
Then $a=-\infty$, and $U^{\prime}\left(  x\right)  =\exp\left(
-\mu x\right)  ,$ so that $B=]0,\ \infty [$ and $I\left(  y\right)
=-\frac{1}{\mu}\ln\left( y\right)  $. The estimate is satisfied for any $r>0.$ \\
iii) Power utility; Set $U\left(  x\right)  =\frac{1}{\mu}x^{\mu}$ for some
$\mu<1$ and $\mu \neq 0$ (note that $\mu$ may be negative). Then $a=0$, and $U^{\prime
}\left(  x\right)  =x^{\mu-1}$, so that $B=]0,\ \infty\lbrack$ and
$I\left(  y\right)  = y^{1/( \mu-1)  }$. The
estimate is satisfied with $r=\frac{1}{1-\mu}.$ \\
iv) Logarithmic utility; Set $U\left(  x\right)  =\ln x$. Then $a=0$ and $U^{\prime
}\left(  x\right)  =\frac{1}{x}$, so that $B=]0,\ \infty\lbrack$ and $I\left(
y\right)  =\frac{1}{y}$. The estimate is satisfied with $r=1.$
\end{example}
Take some $\lambda\in B$ and a $\gamma$ satisfying the conditions of
Definition \ref{market cond}, and define a random variable $X_{\lambda}$ by:%
\[
X_{\lambda}\left(  \omega\right) =I\left(  \lambda\xi_{\timeh}\left(
\omega\right)  \right).
\]
$X_{\lambda}$ is $\mathcal{F}_{\timeh}$-measurable. In addition, we have:
\begin{lemma} \label{X in Lp}
$X_{\lambda}\in L^{p}\left(  \Omega,\mathcal{F}_{\timeh},P\right)  $ for every
$p \geq 1$.
\end{lemma}
\begin{proof}
Since $U$ is admissible, we know from condition 4 that, for some $r>0$ we
have:%
\begin{align*}
\left\vert I\left(  \lambda\xi_{\timeh}\right)  \right\vert ^{p}  &  \leq\left(
c_{1}+c_{2}\left\vert \lambda\xi_{\timeh}\right\vert ^{r}+c_{3}\left\vert
\lambda\xi_{\timeh}\right\vert ^{-r}\right)  ^{p}\\
&  \leq k_{1}+k_{2}\left\vert \lambda\right\vert ^{pr}\left\vert \xi
_{\timeh}\right\vert ^{pr}+k_{3}\left\vert \lambda\right\vert ^{-pr}\left\vert
\xi_{\timeh}\right\vert ^{-pr}%
\end{align*}
and the right-hand side is integrable, for we know that $\xi_{\timeh}^{s} \in
L^{1}\left(  \Omega,\mathcal{F}_{\timeh},P\right)  $ for every $s \in \mathbb{R}.$
%\AAAqed
\end{proof}
\begin{lemma} Let $v \in A.$
There is a unique $\hat{\lambda}\in B$ %
such that
$E_{P}\left[X_{\hat{\lambda}}\xi_{\timeh}\right]  =v$
\end{lemma}
\begin{proof}
Consider the map $\varphi:B\rightarrow \mathbb{R}$ defined by $\varphi\left(
\lambda\right)  =E_{P}\left[  X_{\lambda}\xi_{\timeh}\right]  =E_{P}\left[
I\left(  \lambda\xi_{\timeh}\right)  \xi_{\timeh}\right]  $. Since $\xi_{\timeh}$ $>0$
$P$-a.e., and $I$ is strictly decreasing, $\varphi$ is strictly  decreasing. Using the
Lebesgue dominated convergence theorem, we find that it is continuous. Using
Fatou's lemma, we find that:
\begin{itemize}
\item $\varphi\left(  \lambda\right)  \rightarrow+\infty$ when $\lambda
\rightarrow\inf B$
\item $\limsup$ $\varphi\left(  \lambda\right)  \leq a$ when $\lambda
\rightarrow+\infty$
\end{itemize}
Since $v \in A,$ it follows that there is a unique $\hat{\lambda}$ such
that $\varphi\left(  \hat{\lambda}\right)  =v$.
%\AAAqed
\end{proof}

Denote $X_{\hat{\lambda}}$ by $\hat{X}$. We now conclude:
\begin{proposition} \label{exist unique X}
$\hat{X}$ is the unique solution of problem (P).
\end{proposition}
\begin{proof}
Let us show that $\hat{X}$ is indeed a solution of problem (P). Uniqueness
follows from the strict concavity of $U.$

We have shown that $\hat{X}$ is in $L^{2}$, and $E_{P}[  \hat{X}\xi_{\timeh}]=v,$
 so $\hat{X}$ satisfies the constraints. Take another $X\in L^{2}$ such that
$E_{P}[  \hat{X}\xi_{\timeh}]  =v.$ Since $U$ is concave, we have:%
\[
U\left(  X\left(  \omega\right)  \right)  \leq U\left(  \hat{X}\left(
\omega\right)  \right)  +(X\left(  \omega\right)  -\hat{X}\left(
\omega\right)  )U^{\prime}\left(  \hat{X}\left(  \omega\right)  \right)
\text{ \ \ }P\text{-a.e.}%
\]
By definition, $U^{\prime}\left(  \hat{X}\left(  \omega\right)  \right)
=\lambda\xi_{\timeh}\left(  \omega\right)  $. Substituting into the inequality and
integrating, we get:%
\[
\int_{\Omega}U\left(  X\right)  dP\leq\int_{\Omega}U\left(  \hat{X}\right)
dP+\lambda\int_{\Omega}(X-\hat{X})\xi_{\timeh}dP
\]
and the last term vanishes because it is just $\lambda\left(  v-v\right)  $.
So $\hat{X}$ is indeed an optimizer, and the result follows.
%\AAAqed
\end{proof}
\subsection{Hedging}
Once the solution $\hat{X}$ of the optimization problem ($P$) is found, for a given $\gamma,$
 the question
is whether it can be hedged by a self-financing portfolio $\hat{\theta},$ so
that $V_{\timeh}(\hat{\theta})=\hat{X}.$ We note that, if there exists such
$\hat{\theta} \in \sfprtfs,$ then it is a solution of ($P_{0}$).
In fact, let $\theta \in \sfprtfs$ and $V_{0}(\theta)=v$ and set $X=V_{\timeh}(\theta).$
It follows from ($P$) that
\[
E_{P}\left[U( V_{\timeh}(\theta))\right]=E_{P}\left[U(X)\right]
 \leq E_{P}\left[U(\hat{X})\right]=E_{P}\left[U( V_{\timeh}(\hat{\theta}))\right],
\]
so $\hat{\theta}$ is a solution of ($P_{0}$).
\subsubsection{Deterministic case} \label{determ case}
In this paragraph,  we shall use the general hedging results of
\S \ref{sect. hedging} to solve this problem, in the case when the
 $m$ and $\sigma,$
are deterministic (i.e. they do not depend on $\omega$).

Under these conditions, there can be several $\gamma$ that satisfy the conditions
of Definition \ref{market cond} and some  $\gamma$ can even be non-deterministic.
However, as we have supposed that the market is strongly arbitrage free, so equation
(\ref{rel drift-sigma}) has a solution, we can choose $\gamma$ to be the unique
solution with the property of being orthogonal in $\ell^{2}$ to the kernel of
the volatility operator. More precisely, we choose the unique $\gamma$ such that
\begin{equation} \label{gamma orth}%
(\gamma_{t},x)_{\ell^{2}}=0, \; \;  \forall \ x \in \ell^{2}(\mathbb{I}) \; \; \text{s.t.} \; \; \sigma_{t}x=0.
\end{equation}
The $\gamma$ defined by this condition is deterministic. In the sequel of this
paragraph $\gamma$ is given by (\ref{gamma orth}).
In that case, it follows from formula (\ref{75}) that $\xi_{\timeh}$ is Malliavin
differentiable. It follows from formula (\ref{malliavin 1}) that the partial
derivative with respect to $\wienerq{i}{}$ is given by:%
\[
D_{i,t}\xi_{\timeh}=-\gamma_{t}^{i}\xi_{\timeh}%
\]
and $\hat{X}=I\left(  \hat{\lambda}\xi_{\timeh}\right)  $ is Malliavin
differentiable as well, with:
\[
D_{i,t}\hat{X}=-\hat{\lambda}\gamma_{t}^{i}\xi_{\timeh}I^{\prime}(\hat{\lambda}\xi_{\timeh}).
\]
The Clarke-Ocone formula now reads:%
\begin{align}
X  & =E_{Q}[X\,|\,\mathcal{F}_{0}]+\sum_{i\in\mathbb{I}}\int_{0}^{\timeh}E_{Q}
    \left[  D_{i,t}X\ |\ \mathcal{F}_{t}\right]  d\tilde{W}_{t}^{i} \label{co1}\\
& =v-\hat{\lambda}\sum_{i\in\mathbb{I}}\int_{t}^{\timeh}\gamma_{t}^{i}
E_{Q}\left[\xi_{\timeh}I^{\prime}(\hat{\lambda}\xi_{\timeh})\ |\ \mathcal{F}_{t}\right]
       d\tilde{W}_{t}^{i}\label{co2}%
\end{align}
We then write the equation (\ref{hedge eq 1}) for the hedging portfolio
$\hat{\theta},$ and we
substitute the Clark-Ocone formula for $x_{t}^{i}\left(  \omega\right)  $:%
\begin{equation} \label{76}
b_{t}(\omega)'\theta_{t}(\omega)=\ -\hat{\lambda
}E_{Q}\left[  \xi_{\timeh}I^{\prime}(\hat{\lambda}\xi_{\timeh}%
)\ |\ \mathcal{F}_{t}\right] \gamma_{t}.
\end{equation}
This equation has a solution iff $\gamma_{t}$ is in the range of $b_{t}(\omega)'.$
Since $\sigma$ is deterministic, this condition simplifies. In fact, let
$l_{t}$ and $B_{t}$ be given by (\ref{B and l}), which here both are deterministic,
and let $q(t,\omega)=\zcpxd{t}(\omega)/l_{t}.$ Then the expression (\ref{self-fin1}) of $b_{t}(\omega)'$
give:
\[
(b_{t}(\omega)'\theta_{t}(\omega))^{i}=<\theta_{t}(\omega)\,,\,p_{t}(\omega)\sigma_{t}^{i}>
=<\theta_{t}(\omega)q(t,\omega)\,,\, l_{t} \sigma_{t}^{i}>
=(B_{t}'f_{t}(\omega))^{i},
\]
where $f_{t}(\omega) \in E^{-s}$ is given by $f_{t}(\omega)=q(t,\omega)\theta_{t}(\omega).$
So, equation (\ref{76}) has a solution iff $\gamma_{t}$ is in the range of $B_{t}'.$
This is always true when $\mathbb{I}$ is finite, since then the range of $B_{t}'$
is equal to the orthogonal complement of the kernel of $\sigma_{t}$
(we remember that $p_{t}(\omega,x) >0$ for $x\geq 0$). When $\mathbb{I}=\mathbb{N},$
then the range is only  a strictly smaller dense subset.

We are lead to following condition
\begin{definition} \label{cond. C}
We shall say that the market satisfies condition (C) if there exists a
deterministic portfolio $\theta_{t}^{0}$ which is admissible and satisfies
$B_{t}' \theta_{t}^{0}=\gamma_t,$ i.e.
\begin{equation}
<\theta_{t}^{0}\,,\left(  \mathcal{L}_{t}p_{0}\right)  \sigma_{t}%
^{i}>\ =\gamma_{t}^{i}, \label{77}%
\end{equation}
for each $i\in\mathbb{I}$ and $t$.
\end{definition}
Condition $C$ is then equivalent to $\gamma_{t} \in \mathcal{R}(B_{t}'),$
the range of $B_{t}'.$ %
In the case when $\mathbb{I}$ is finite, there is never
uniqueness in the choice of $\theta_{t}^{0}.$ \\

In the case when $\mathbb{I}$ is finite, we know that condition (C) is satisfied
and it can easily  be verified, with $n$ elements say,
 by picking $n$ maturities $0<S_{1}<...<S_{n}$ and by
seeking $\theta_{t}^{0}$ as a linear combination of rollovers:\ $\theta
_{t}^{0}=\sum x_{t}^{i}\delta_{S_{i}}$. Condition (\ref{77}) then reduces to a
system of $n$ linear equations with $n$ unknowns which determines the
$x_{t}^{i}$.

In the case when $\mathbb{I=\mathbb{N}}$, condition (C) may not be satisfied. We will
be content with %
reminding that the left-hand side of equation (\ref{77}) is
meaningful, since $\left(  \mathcal{L}_{t}p_{0}\right)  \sigma_{t}^{i}$
belongs to the space $E^{s}.$ %

If condition (C) is satisfied, equation (\ref{76}) becomes:%
\begin{equation*}
\begin{split}
<\theta_{t},\,\,p_{t}\sigma_{t}^{i}>\ & =-\hat{\lambda}E_{Q}[  \xi
_{\timeh}I^{\prime}(\hat{\lambda}\xi_{\timeh})\ |\ \mathcal{F}_{t}]  <\theta
_{t}^{0},\frac{\,\mathcal{L}_{t}p_{0}}{p_{t}}p_{t}\sigma_{t}^{i}> \\
&=<-\hat{\lambda}E_{Q}[  \xi
_{\timeh}I^{\prime}(\hat{\lambda}\xi_{\timeh})\ |\ \mathcal{F}_{t}] \
      \frac{\,\mathcal{L}_{t}p_{0}}{p_{t}}\  \theta_{t}^{0} \ , \ p_{t}\sigma_{t}^{i}>
\end{split}
\end{equation*}
and an obvious  solution $\theta_{t}=\bar{\theta}_{t}$ (the risky part of the optimal portfolio) is
given by:%
\[
\bar{\theta}_{t}=-\hat{\lambda}E_{Q}[  \xi
_{\timeh}I^{\prime}(\hat{\lambda}\xi_{\timeh})\ |\ \mathcal{F}_{t}] \
      \frac{\,\mathcal{L}_{t}p_{0}}{p_{t}}\  \theta_{t}^{0}.
\]
Applying Lemma \ref{lm hedge eq}, with $x$ defined by (\ref{sol hedge eq1}),
we obtain a hedging portfolio
 $\hat{\theta}=x \delta_{0}+\bar{\theta}$ of $\hat{X},$ where $\bar{\theta}$ is as above, and:
\[
x_{t}=\frac{1}{p_{t}\left(  0\right)  }\left(  E_{Q}\left[  I\left(
\hat{\lambda}\xi_{\timeh}\right)  \ |\ \mathcal{F}_{t}\right]  -<\bar{\theta}%
_{t},p_{t}>\right).
\]
To sum up, in the case when the %
$m_{s}$ and the $\sigma_{s}^{i},i\in\mathbb{I}$, are deterministic,
with $\sigma_{t}^{i}\left(  0\right)=0,$
 with condition $(C)$ and equation (\ref{rel drift-sigma})  satisfied,
  an optimal admissible and self-financing
portfolio is given by
\begin{equation} \label{sum up}
 \hat{\theta}_{t}=x_{t}\delta_{0}+\bar{\theta}_{t}, \; \; \text{where} \; \; 
\bar{\theta}_{t}\ =y_{t}\ \frac{(\mathcal{L}_{t}p_{0})}{p_{t}} \ \theta_{t}^{0}
\end{equation}
and where the coefficients $x_{t}$ and $y_{t}$ are real-valued progressively measurable processes
given by
\begin{align}
y_{t}  &  = -E_{Q}[\hat{\lambda}\xi_{\timeh}I^{\prime}(\hat{\lambda}\xi
_{\timeh})\,|\,\mathcal{F}_{t}]\label{prtf explicit 1}\\
x_{t}  &  =(p_{t}(0))^{-1}\left(  E_{Q}[I(\hat{\lambda}\xi_{\timeh}%
)\,|\,\mathcal{F}_{t}]-y_{t}<\theta_{t}^{0}\,,\,\mathcal{L}_{t}p_{0}>\right).
\label{prtf explicit 2}%
\end{align}

This leads immediately to a mutual fund theorem: whatever the utility function
and the initial wealth, the optimal portfolio at time $t$ is a linear combination of the current
account $\delta_{0}$ and the portfolio $f \mapsto <\theta_{t}^{0}%
,\frac{\,\ltrans{t}p_{0}}{p_{t}}f>,$ i.e. the portfolio
$\frac{\,\ltrans{t}p_{0}}{p_{t}}\theta_{t}^{0}.$
This portfolio is in general not self-financed, so it can
not be given the status of a \textit{market portfolio}. However we can easily reformulate
the result with a self-financed portfolio. In fact, chose an admissible utility
function, with $a=0,$ according to Definition \ref{85}. For this utility function,
let $\Theta$ be the optimal portfolio given by (\ref{sum up}), with unit initial wealth.
Obviously
$\frac{\,\ltrans{t}p_{0}}{p_{t}}\theta_{t}^{0}$ is a linear combination of
$\delta_{0}$ and $\Theta_{t}.$ This gives us:
\begin{theorem}[Mutual fund theorem] \label{Th mutual fund}
The optimal portfolio $\Theta$ %
has the following properties:  \\

\noindent i) $\Theta$ is an admissible self-financing portfolio,
with unit initial value, i.e. $\sesq{\Theta_{0}}{\zcpxd{0}}=1,$
and the value at each time $t \in \mathbb{T}$ is strictly positive, i.e.  $\sesq{\Theta_{t}}{\zcpxd{t}}>0.$ \\

\noindent ii) For each utility function $U,$ admissible according to Definition \ref{85}
and each initial wealth $v \in \;]a,\infty [\,,$
there exist two real valued processes $c$ and $d$ such that if
$\hat{\theta}_{t} =c_{t}\delta_{0}+d_{t}\Theta_{t},$
then $\hat{\theta}$ is an optimal self financing portfolio for $U,$
i.e. a solution of problem  ($P_{0}$).
\end{theorem}
\subsubsection{Stochastic $m$ and $\sigma$}
We shall here concentrate on the case of an approximately complete market, which
is equivalent to that the volatility operator is  non-degenerated. In fact,
according to $iii)$ of Theorem \ref{Th D0 non complete and approx complete},
the market is approximately complete if and only if $\vol{}{t}(\omega)$ has
a trivial kernel a.e. $(t,\omega) \in \mathbb{T} \times \Omega.$ We remind that
the market of price process $\gamma$ is unique in this case.

In the case of a finite number of Bm. we obtain easily from Lemma \ref{X in Lp}
and  Theorem \ref{th completeness R^m} the following result
(see Theorem 3.6 of \cite{I.E.-E.T bond th}):
\begin{theorem}\label{th opt port R^m}
Let $\mathbb{I}$ be a finite set, let $U$ be admissible in the sens of Definition %
\ref{85} %
and let  $i)$ of Condition \ref{uniform cond sigma} %
be satisfied.
The problem ($P_{0}$) %
then has a solution $\hat{\theta}.$ One solution $\hat{\theta}=a \delta_{0}+\bar{\theta} \in \sfprtfs$
is  given by (\ref{prtf eq 4}). %
\end{theorem}
In the case of an infinite number of Bm. we shall impose Malliavin differentiability
properties on the market price of risk $\gammapx{}{}.$ To this end we introduce
the space $\derprod{1}{s},$ for $s>0$ by
\begin{equation} \label{def D^{1}_{s}}
\derprod{1}{s} =\{X \in \derprod{}{0} \; | \; DX \in \derprod{}{0}(L^{2}(\mathbb{T},\ell^{s,2})) \}.
\end{equation}
We can now state a result in the case of an infinite number of Bm., quite analog
to the case of a finite number of Bm. (see Theorem 4.5 of \cite{E.T Bond  Completeness}):
\begin{theorem}\label{th opt port l^2}
Let $\mathbb{I}=\mathbb{N},$
let $U$ be admissible in the sens of Definition \ref{85}, %
let $ii)$ of Condition \ref{uniform cond sigma} %
be satisfied
and let $\ln(\xi_{\timeh}) \in \derprod{1}{s},$ where $s>0$ is given by $ii)$ of Condition \ref{uniform cond sigma}.
The problem ($P_{0}$) %
then has a solution $\hat{\theta}.$ One solution $\hat{\theta}=a \delta_{0}+\bar{\theta} \in \sfprtfs$
 is  given by (\ref{prtf eq 4}). %
\end{theorem}
\begin{proof}
We only consider the case of $U'>0,$ since the case of $U'(x)=0$ for some $x$ is so similar.
Let the hypotheses of the theorem be satisfied. The portfolio $\hat{\theta}$ is a solution of
equation ($P_{0}$), if $\hat{\theta} \in \sfprtfs$ and if it  hedges
$\hat{X}$ given by Proposition \ref{exist unique X}.
(See  Corollary 3.4 of \cite{I.E.-E.T bond th}).
It is enough to verify that Theorem \ref{th completeness l^2} applies to 
$\hat{X}=I (\hat{\lambda} \xi_{\timeh})$ for a certain given $\hat{\lambda} >0.$

$I $ is $C^{1},$  so
$\mder{t}\hat{X}=\lambda \xi_{\timeh}\varphi' (\lambda \xi_{\timeh})\mder{t}\ln(\xi_{\timeh}).$
Since $\ln(\xi_{\timeh}) \in \derprod{1}{s},$ this gives
$\|\mder{}\hat{X}\|_{L^{2}( \mathbb{T},\ell^{s,2})}
=|\lambda \xi_{\timeh}\varphi' (\lambda \xi_{\timeh})| \, \|\mder{}\ln(\xi_{\timeh})\|_{L^{2}( \mathbb{T},\ell^{s,2})}.$
The inequality in 5) of Definition \ref{85} gives  %
$\|\mder{}\hat{X}\|_{L^{2}( \mathbb{T},\ell^{s,2})}
\leq C ((\lambda \xi_{\timeh})^{p}
+(\lambda \xi_{\timeh})^{-p}) \|\mder{}\ln(\xi_{\timeh})\|_{L^{2}( \mathbb{T},\ell^{s,2})},$
for some $p \geq 1.$
Condition (\ref{gamma strong}) of Definition \ref{market cond}
 shows that $(\lambda \xi_{\timeh})^{p}+(\lambda \xi_{\timeh})^{-p} \in L^{q}(\Omega, P),$
for all $q \geq 1.$ By hypothesis
$\|\mder{}\ln(\xi_{\timeh})\|_{L^{2}( \mathbb{T},\ell^{s,2})} \in \derprod{}{0},$
so H\"older's inequality now gives that
$\|\mder{}\hat{X}\|_{L^{2}( \mathbb{T},\ell^{s,2})} \in \derprod{}{0},$
i.e. $\mder{}\hat{X} \in \derprod{}{0}(L^{2}( \mathbb{T},\ell^{s,2})).$
By Lemma \ref{X in Lp}, $\hat{X} \in \derprod{}{0}.$
It  follows that $\hat{X} \in \derprod{1}{s}.$
We can now apply Theorem \ref{th completeness l^2}, which proves the existence of $\hat{\theta}.$
%\AAAqed
\end{proof}
\subsubsection{Examples.} \label{exampl}
We now give some examples of optimal bond portfolios for logarithmic and
quadratic utility functions $U.$ Other examples can be found in \cite{I.E.-E.T bond th}.

First we assume the drift function $m_{t}$ and the volatility operator $\sigma_{t}$
to be deterministic.
We shall therefore suppose that the market satisfy condition
$(C),$ of Definition \ref{cond. C}, so the market prices of risk $\gamma$
is deterministic and satisfy  condition (\ref{77}).
We shall  derive the optimal portfolio directly, going through
the steps leading to the general solution (\ref{sum up}).

Secondly we study the general case of stochastic drift function $m_{t}$ and volatility operator $\sigma_{t}$
for the logarithmic utility function.

The final optimal discounted wealth is $\hat{X}=I(\hat{\lambda}\xi_{\timeh})$. The
corresponding optimal discounted wealth process $Y$is given by
$Y_{t}=E_{Q}[I(\hat{\lambda}\xi_{\timeh})\,|\,\mathcal{F}_{t}].$ The initial wealth
$Y_{0}=v$ determines $\hat{\lambda}$ by the equation
\begin{equation}
v=Y_{0}=E_{Q}[I(\hat{\lambda}\xi_{\timeh})]. \label{81}%
\end{equation}
We recall that $(p_{t})^{-1}\mathcal{L}_{t}p_{0}\in E^{s}$ a.s and that
$p_{t}(0)>0$ a.s.
\paragraph{Logarithmic utility (deterministic $m$ and $\sigma$)}
Let
\begin{equation}
U(x)=\ln(x). \label{util-fnct-log}
\end{equation}
We have $I(x)=1/x,$ and $\hat{X}=(  \hat{\lambda}\xi_{\timeh})^{-1},$
so that equation (\ref{81}) gives:
\[
v=E_{Q}[1/(\hat{\lambda}\xi_{\timeh})]=E_{P}[\xi_{\timeh}/(\hat{\lambda}\xi_{\timeh}%
)]=1/\hat{\lambda}.
\]
Then using the expression (\ref{75}) for $\xi_{t}$ and $\tilde{W}_{t}^{i}$ we have:%
\begin{equation}
\frac{1}{\xi_{t}}=\exp\left(  -\frac{1}{2}\int_{0}^{t}\sum_{i\in\mathbb{I}%
}\left(  \gamma_{s}^{i}\right)  ^{2}ds+\int_{0}^{t}\sum_{i\in\mathbb{I}}%
\gamma_{s}^{i}d\tilde{W}_{s}^{i}\right).  \label{82}%
\end{equation}
The right-hand side is a $Q$-martingale, then so is $1/\xi_{t}$. It follows
that the optimal discounted wealth at $t$ is%
\[
Y_{t}=E_{Q}[I(\hat{\lambda}\xi_{\timeh})\,|\,\mathcal{F}_{t}]=\frac{1}{\hat
{\lambda}\xi_{t}}=\frac{v}{\xi_{t}}.
\]
Since $d(1/\xi_{t})=\sum_{i\in\mathbb{I}}(\gamma_{t}^{i}/\xi_{t})d\tilde{W}_{t}^{i}$
and $\hat{X}  =Y_{\timeh},$
it then follows that: %
\begin{equation}
\hat{X}  =v\left(  1+\sum_{i\in\mathbb{I}}\int_{0}^{\timeh}\gamma_{t}^{i}\frac{1}{\xi_{t}
}d\tilde{W}_{t}^{i}\right).
\end{equation}
The hedging equation (\ref{hedge eq 1}) and the above formula give:
\begin{equation} \label{log 1}
\forall i\in\mathbb{I},\ \  <\theta_{t}\left(  \omega\right)  ,\ p_{t}\left(
\omega\right)  \,\sigma_{t}^{i}\left(  \omega\right)  >\ =\  \frac{v}{\xi
_{t}\left(  \omega\right)  }\gamma_{t}^{i}%
\end{equation}
By condition (C)  we find a portfolio $\theta^{0}$ satisfying
   $\gamma_{t}^{i}=\ <\theta_{t}^{0}\,,\left(  \mathcal{L}_{t}p_{0}\right)\sigma_{t}^{i}>,$
so
\begin{equation} \label{log 1.1}
\gamma_{t}^{i}=\ <\left(  \mathcal{L}_{t}p_{0}\right) \theta_{t}^{0}\,,\sigma_{t}^{i}>.
\end{equation}
Substituting this expression of $\gamma$ into (\ref{log 1}) we obtain:
\begin{equation} \label{log 2}
\forall i\in\mathbb{I},\ \ \
< p_{t}\left(\omega\right) \theta_{t}\left(\omega\right)
  -\frac{v}{\xi_{t}(\omega)}(\mathcal{L}_{t}p_{0}) \ \theta_{t}^{0}
 ,\ \,\sigma_{t}^{i}\left(  \omega\right)  >\ = 0.
\end{equation}
One solution of this equation is obviously given by $\theta=\bar{\theta},$ where
\begin{equation} \label{log 3}
\bar{\theta}_{t}(\omega)
=y_{t}(\omega) \ \frac{(\mathcal{L}_{t}p_{0}) }{ p_{t}(\omega)} \ \theta_{t}^{0}, \; \; 
y_{t}(\omega)= \frac{v}{\xi_{t}(\omega)}.
\end{equation}
The discounted value of $\bar{\theta}$ at time $t$ in state $\omega$ is then
\begin{equation} \label{log 4}
(V_{t}(\bar{\theta}))(\omega)=<\bar{\theta})_{t}\,, p_{t}>
=\frac{v}{\xi_{t}(\omega)} \ <\theta_{t}^{0}\,,  \mathcal{L}_{t}p_{0}>.
\end{equation}
The optimal portfolio $\hat{\theta}$ is now obtained by using Lemma \ref{lm hedge eq}:
$\hat{\theta}_{t}  =x_{t}\delta_{0}+\bar{\theta}_{t},$ where
\begin{equation} \label{log 5}
x_{t}  =\frac{1}{p_{t}(0)}\frac{v}{\xi_{t}}(1-<\theta_{t}^{0}\,,\hat{\theta
}_{t}\left(  \omega\right)  \mathcal{L}_{t}p_{0}>).
\end{equation}
As it should, the discounted value of $\hat{\theta}$ is then
$V_{t}(\hat{\theta})=Y_{t}=v/\xi_{t}.$

We note the following useful property:
the ratio of the investment in bonds with time to maturity $S>0$ to the total
investment is deterministic. In fact this ratio is simply price
at $t=0,$ of a zero-coupon bond with time to maturity $S+t:$
\begin{equation} \label{log 6}
\frac{\bar{\theta}_{t}(S,\omega) \ p_{t}(S,\omega)}{(V_{t}(\bar{\theta}))(\omega)}
=p_{0}(S+t).
\end{equation}
\paragraph{Quadratic utility (deterministic $m$ and $\sigma$)}
Let the utility function be:
\[
U\left(  x\right)  =\mu x-\frac{1}{2}x^{2}%
\]
As in $i)$ of Example \ref{U example}, we find that
\[
I(y) =\mu -y.
\]
The final discounted optimal wealth is $\hat{X}=I(\hat{\lambda}\xi_{\timeh}),$ so
\[
\hat{X} =\mu  -\hat{\lambda}\xi_{\timeh}.
\]
We determine $\hat{\lambda}$ by the condition:%
\begin{equation} \label{quadr 1}
v    =E_{Q}\left[  \hat{X}\right]  =E_{Q}\left[ \mu 
-\hat{\lambda}\xi_{\timeh} \right]
  =\mu  -\hat{\lambda}E_{Q} \left[  \xi_{\timeh}\right] .
\end{equation}
Set
\[
Z_{t}=\exp{ \left( -\frac{1}{2}\int_{0}^{t}\sum_{i \in \mathbb{I}} (\gamma^{i}_{s})^{2}ds
    -\int_{0}^{t}\sum_{i \in \mathbb{I}} \gamma^{i}_{s} d\wienerq{i}{s})\right)}.
\]
Then $Z$ is a martingale with respect to $Q$ and formula (\ref{W Q}) gives
\begin{equation} \label{quadr 2}
\xi_{t}=Z_{t} \exp{ \left(\int_{0}^{t}\sum_{i \in \mathbb{I}} (\gamma^{i}_{s})^{2}ds\right)}.
\end{equation}
We have, by substitution into (\ref{quadr 1}):
\[
v =\mu  -\hat{\lambda}E_{Q} \left[  \xi_{\timeh}\right] 
=\mu  -\hat{\lambda}
  \exp{ \left(\int_{0}^{\timeh}\sum_{i \in \mathbb{I}} (\gamma^{i}_{s})^{2}ds\right)}.
\]
This gives
\begin{equation} \label{quadr 4}
\hat{\lambda}
=\left( \mu -v \right)
  \exp{ \left( - \int_{0}^{\timeh}\sum_{i \in \mathbb{I}} (\gamma^{i}_{s})^{2}ds\right)}.
\end{equation}
It now follows from (\ref{quadr 2}) that
\begin{equation} \label{quadr 5}
\hat{X} =\mu  -\hat{\lambda}\xi_{\timeh}
=\mu + \left( v -\mu \right) \ Z_{\timeh}
\end{equation}
and the optimal discounted wealth at $t$ is%
\[
Y_{t}=E_{Q}[I(\hat{\lambda}\xi_{\timeh})\,|\,\mathcal{F}_{t}]
=\mu + \left( v -\mu \right) \ Z_{t}.
\]
Since $d Z_{t}=- Z_{t} \sum_{i \in \mathbb{I}} \gamma^{i}_{t} d\wienerq{i}{t},$
 we have that
\begin{equation*} %
\hat{X}
=\mu - \left( v -\mu \right) \
    \int_{0}^{\timeh} \sum_{i \in \mathbb{I}}Z_{t}  \gamma^{i}_{t} d\wienerq{i}{t}
=\mu +
  \int_{0}^{\timeh} \sum_{i \in \mathbb{I}} (\mu -Y_{t})  \gamma^{i}_{t} d\wienerq{i}{t},
\end{equation*}
so the hedging equation reads (see (\ref{hedge eq 1})):
\begin{equation} \label{quadr 7}
\forall i\in\mathbb{I},\ \
<\theta_{t}\left(  \omega\right)  ,\ p_{t}\left(\omega\right)  \,\sigma_{t}^{i}\left(  \omega\right)  >\ =\
   -      \left(\mu -Y_{t} (\omega) \right) \ \gamma^{i}_{t}.
\end{equation}
As usually, condition (C) gives a portfolio $\theta^{0}$ satisfying
   $\gamma_{t}^{i}=\ <\theta_{t}^{0}\,,\left(  \mathcal{L}_{t}p_{0}\right)\sigma_{t}^{i}>,$
which together with (\ref{quadr 7}) gives:
\begin{equation*} %
\forall i\in\mathbb{I},\ \ \
< p_{t}\left(\omega\right) \theta_{t}\left(\omega\right)
  +  \left( Y_{t}(\omega) -\mu \right)  \  (\mathcal{L}_{t}p_{0}) \ \theta_{t}^{0}
 ,\ \,\sigma_{t}^{i}\left(  \omega\right)  >\ = 0.
\end{equation*}
One solution of this equation is $\theta=\bar{\theta},$ where
\begin{equation*} %
\bar{\theta}_{t}(\omega)=y_{t}(\omega) \ 
    \frac{(\mathcal{L}_{t}p_{0}) }{ p_{t}(\omega)} \ \theta_{t}^{0}, \; \; 
y_{t}(\omega)= \mu -Y_{t}(\omega) .
\end{equation*}
$\bar{\theta}$  gives the risky part of the optimal portfolio.

Applying Lemma \ref{lm hedge eq}  we obtain the  optimal portfolio
$\hat{\theta}_{t}  =x_{t}\delta_{0}+\bar{\theta}_{t},$ where
\begin{equation} \label{quadr 9}
x_{t}= (\zcpxd{t}(0))^{-1} (Y(t)-(\mu -Y(t))\sesq{\theta_{t}^{0}}{\ltrans{t}\zcpx{0}}).
\end{equation}

\paragraph{Logarithmic utility (stochastic $m$ and $\sigma$)}
We assume that the conditions of Definition \ref{market cond} are satisfied.
We chose $\gamma_t(\omega)$ to be orthogonal to the kernel of
$\vol{}{t}(\omega),$ a.e. $(t,\omega).$ This $\gamma$ satisfies the conditions of
Definition \ref{market cond}. Formulas (\ref{util-fnct-log})--(\ref{log 1}) then
still hold true. As in the discussion preceding the condition $(C),$
of Definition \ref{cond. C} it follows that $\gamma_t(\omega)$ is a.s. in the closure
of the range of $B_{t}'(\omega).$ Therefore, in this example,
the natural generalization of the condition $(C)$
to the stochastic case is simply to impose the same condition (\ref{77})
of Definition \ref{cond. C} to be satisfied with a stochastic portfolio $\theta^{0} \in \prtfs.$
Formulas (\ref{log 1.1})--(\ref{log 6}) are then also true statements and it
follows using Theorem \ref{Th price compatible strong cond.} that $\hat{\theta} \in \sfprtfs.$
In particular the ratio of the investment in bonds with time to maturity $S>0$ to the total
investment is deterministic.

\subsection{The H-J-B approach}
When $m_{t}$ and $\sigma_{t}^{i}$ are given functions
$m_{t}(p_{t})$ and $\sigma_{t}^{i}(p_{t})$ of the price $p_{t},$ for every
$t,$ then the optimal portfolio problem ($P_{0}$) can be considered within a
Hamilton-Jacobi-Bellman approach. In this subsection we illustrate this
approach, without being rigorous and we suppose that the utility function $U$
satisfies the conditions of Definition \ref{85}. For notational simplicity
we exclude the price argument in $m_{t}$ and $\sigma_{t}^{i}.$

The optimal value function, here denoted by $F,$ then only depends of time
$t,$ of the value of the discounted wealth $w$ and the discounted price
function $f\in E^{s}$ of Zero-Coupons at time $t:$
\[
F(t,w,f)=\sup\{E[U(V_{\timeh}(\theta))\;|\;V_{t}(\theta)=w,\;p_{t}=f] \;|\;\theta
\in\mathsf{P}_{sf}\}.
\]
The derivative  $DG(f;g)$ of a function $E^{s} \ni f \mapsto G(f)$ in the direction
$g \in E^{s},$ is as usually defined by
\[
DG(f;g)=\lim_{\epsilon \rightarrow 0} \frac{G(f+\epsilon g)-G(f)}{\epsilon}.
\]
Suppose that $G$ is $C^{2}.$
Writing $DG(f)$ for the map $g \mapsto DG(f;g)$ and $D^{2}G(f)$ for the map
$g_{1} \times g_{2} \mapsto DG(f;g_{1}, g_{2}),$ we have that %
$DG(f)$ is a linear continuous form on $E^{s}$ and $D^{2}G(f)$
is a bi-linear continuous form.

Let us first consider the case of a volatility operator $\vol{}{}$ with
trivial kernel, i.e. for every strictly positive price (function) $f \in E^{s},$
the kernel of the linear map $\vol{}{t}: \ell^2(\mathbb{I}) \rightarrow E^{s}$
is trivial a.s. According to Definition \ref{market cond} there is then a
unique market of price process $\gamma.$ Define the Hamiltonian $H(t,w,f,x)$ by: %
\begin{equation}
\begin{split}
&H(t,w,f,x)=  \sum_{i \in \mathbb{I}} x^{i}(t,w,f)\gamma^{i}_{t}  \frac{\partial F}{\partial w} \, (t,w,f)
     +DF(t,w,f; \partial f +\sum_{i \in \mathbb{I}} \gamma^{i}_{t} \vol{i}{t}f) \\
      & + \sum_{i \in \mathbb{I}}  \bigl( \frac{1}{2}(x^{i}(t,w,f))^{2}\frac{\partial^{2} F}{\partial w^{2}}  \, (t,w,f) %
       +x^{i}(t,w,f) \frac{\partial}{\partial w} DF(t,w,f;  \vol{i}{t}f) \\
      &+\frac{1}{2} D^{2}F(t,w,f;  \vol{i}{t}f, \vol{i}{t}f) \bigr).
\end{split}
 \label{HJB hamiltonian eq inf dim}
\end{equation}
In that formula, $x=\left(  x^{i}\right)  _{i\in\mathbb{I}} \in \ell^{2}$ is the control,
which is related to the optimal terminal wealth by formula
(\ref{mart decomp 1}). A control $x$ is called \emph{admissible} if
\begin{equation}
x^{i}(t,V_{t}(\theta),p_{t})=<\theta_{t}\,,\,p_{t}\sigma_{t}^{i}>
\label{HJB control x}
\end{equation}
for all $\theta \in\mathsf{P}_{sf}$. In other words, $x^{i}$ can be
interpreted as the value invested in the $i$-th source of noise. Using
the Ito formula, one derives the (formal) HJB equation:
\begin{equation}
\frac{\partial F}{\partial t}\,(t,w,f)+\sup_{x} H(t,w,f,x)=0,
\label{HJB eq inf dim}%
\end{equation}
with the boundary condition
\begin{equation}
F(\timeh,w)=U(w).
 \label{HJB bound cond inf dim}
\end{equation}
The optimal control $\hat{x},$ solution of the optimization problem
\[
\sup_{x} H(t,w,f,x),
\]
is given by
\begin{equation}  \label{opt contr}
\hat{x}^{i}(t,w,f)=-\left(  \frac{\partial^{2}F}{\partial w^{2}}\right)
^{-1}\left(  \gamma_{t}^{i}\frac{\partial F}{\partial w}
  +(D\frac{\partial F}{\partial w})(t,w,f; \sigma_{t}^{i}f)\right)  ,\;i\in\mathbb{I}.
\end{equation}
Now, substitution of $H(t,w,f,\hat{x}(t,w,f))$ into equation (\ref{HJB eq inf dim}) gives:%
\begin{equation}  \label{HJB eq 2}
\begin{split}
&  \frac{\partial^{2}F}{\partial w^{2}}(t,w,f)
     \Big(\frac{\partial F}{\partial t}(t,w,f)
    +DF(t,w,f;\partial f+m_{t}f) \\ &+\frac{1}{2}\sum_{i\in\mathbb{I}}
    D^{2}F(t,w,f;\sigma_{t}^{i}f,\sigma_{t}^{i}f) \Big)
  =\frac{1}{2}\sum_{i\in\mathbb{I}}\left(\gamma_{t}^{i}\frac{\partial F}{\partial
w}+(D\frac{\partial F}{\partial w})(t,w,f;\sigma_{t}^{i}f)\right)^{2}.
\end{split}
\end{equation}
Once the solution $F$ of (\ref{HJB eq 2}),
with boundary condition (\ref{HJB bound cond inf dim}), is found,  the optimal control $\hat{x}$
is given by (\ref{opt contr}).
Any optimal portfolio $\hat{\theta}$ is then a solution of the equation:%
\[
\hat{x}^{i}(t,V_{t}(\hat{\theta}),p_{t})=<\hat{\theta}_{t}\,,\,p_{t}\sigma_{t}^{i}>, \; \;
\forall \ \ i\in\mathbb{I}, \ \ t \in \mathbb{T}. %
\]

Next we consider the case of a volatility operator, which does not necessarily
have a trivial kernel. Once more we define the Hamiltonian $H(t,w,f,x,\gamma)$ by
formula (\ref{HJB hamiltonian eq inf dim}), which now also depends on the control
$\gamma,$ a  $\ell^2(\mathbb{I})$ valued function of $(t,w,f).$ A control $(x,\gamma)$
is admissible if condition (\ref{HJB control x}) is satisfied and if the conditions
of Definition \ref{market cond}  are satisfied, so writing out the price argument
$f \in E^s$ in $m_{t}$ and $\sigma_{t}^{i}:$
\begin{equation}
m_{t}(f)=\sigma_{t}(f) \gamma_{t}(w,f).
\label{HJB control gamma}
\end{equation}
The optimal control $\hat{\gamma}$ is determined by conditions (\ref{HJB control x})
 and (\ref{HJB control gamma}). This can be seen as follows. Let $\gamma^{\perp}(f)$
be the unique solution of (\ref{HJB control gamma}) such that $\gamma^{\perp}(f)$
is in the orthogonal complement $(\mathcal{K}(\sigma_{t}(f)))^{\perp}$
of the kernel $\mathcal{K}(\sigma_{t}(f)),$ let
$\hat{\alpha} =\hat{\gamma} -  \gamma^{\perp}$ and let
$P_t(f)$ be the orthogonal projection on  $\mathcal{K}(\sigma_{t}(f)).$ 
Condition (\ref{HJB control x})
implies that $\hat{x} \in (\mathcal{K}(\sigma_{t}(f)))^{\perp}.$
According to (\ref{opt contr}), this can only be satisfied if
\begin{equation}  \label{opt contr gamma}
\hat{\gamma} =  \gamma^{\perp} +\hat{\alpha}
\;\; \text{and} \;\;
\hat{\alpha}_t(w,f) \frac{\partial F}{\partial w}= P_t(f) \nu_t(w,f),
\end{equation}
where $\nu_t^i (w,f)=(D\frac{\partial F}{\partial w})(t,w,f; \sigma_{t}^{i}f).$
So in the general the case of a volatility operator, which does not necessarily
have a trivial kernel, the H-J-B approach leads to the equation (\ref{HJB eq 2}),
with $\gamma$ replaced by $\hat{\gamma}$ defined by formula (\ref{opt contr gamma}).

In the case when $m_{t}$ and $\sigma_{t}^{i}$ are independent of $p_{t},$ then
the $\hat{x}^{i}$ are independent of $f,$ $\gamma = \gamma^{\perp}$
 and the above equations simplify:%
\[
\frac{\partial F}{\partial t}\frac{\partial^{2}F}{\partial w^{2}}=\frac{1}%
{2}\left(  \sum_{i\in\mathbb{I}}\Vert\gamma_{t}^{i}\Vert^{2}\right)
(\frac{\partial F}{\partial w})^{2},
\]
with the boundary condition
\[
F(\timeh,w)=U(w),\;w\in\mathbb{R}.
\]
Each self financing portfolio $\hat{\theta}\in\mathsf{P}_{sf},$ such that
\[
<\hat{\theta}_{t}\,,\,p_{t}\sigma_{t}^{i}>=-\gamma_{t}^{i}\left(
\frac{\partial F}{\partial w}\right)  \left(  \frac{\partial^{2}F}{\partial
w^{2}}\right)  ^{-1},\; \; \forall \ i\in\mathbb{I}, \ t \in \mathbb{T},
\]
where $w=V_{t}(\hat{\theta}),$ is then a solution of problem ($P_{0}$).
The solutions in the examples in \S \ref{exampl}, as well as the general solution
(\ref{sum up}) for deterministic $m$ and $\sigma,$ are easily obtained by solving
these equations.

\appendix
\section{Appendix}
\label{App}
In this appendix, we reproduce results (proved in the appendix of \cite{I.E.-E.T bond th}),
used in this article, concerning existence
of solutions of some SDE's and $L^{p}$ estimates of these solutions.
The notations $\mathbb{T}=[0, \timeh],$ $W^{i},$ $\mathbb{I}$ and
$(\Omega,P,\mathcal{F},\mathcal{A})$
are defined in \S \ref{non-param frame}.
Through the appendix  $\drift{}$ and $\vol{i}{},$  $i \in \mathbb{I},$
are $\mathcal{A}$-progressively measurable $\multsp$-valued processes satisfying
\begin{equation}
 \int_{0}^{\timeh}(\|\drift{t} \|_{\multsp}+\sum_{i \in \mathbb{I}}\| \vol{i}{t} \|^{2}_{\multsp})dt < \infty, a.s.
 \label{App cond 1}
\end{equation}
The $\multsp$-valued semi-martingale $L$ is given by
\begin{equation}
 L(t)=\int_{0}^{t}(\drift{s}ds + \sum_{i \in \mathbb{I}}\vol{i}{s} d\wienerp{i}{s}), \quad \text{if} \; 0 \leq t \leq \timeh
 \label{A. notation 1}
\end{equation}
and by $L(t)=L(\timeh),$ if $t > \timeh.$ 
We introduce, for $t \geq 0,$ the random variable
\begin{equation}
 \mu(t)=t+\int_{0}^{t}(\| \drift{s} \|_{\multsp}+\sum_{i \in \mathbb{I}}\| \vol{i}{s} \|^{2}_{\multsp})ds, \quad \text{if} \; 0 \leq t \leq \timeh
 \label{A. notation 2}
\end{equation}
and  $\mu(t)=t-\timeh+\mu(\timeh)$ if $t > \timeh.$
$\mu$ is a.s. strictly increasing, absolutely continuous and on-to $[0,\infty[.$ The inverse $\tau$
of $\mu$ also have these properties and $\tau (t) \leq t.$
For a continuous $\multsp$-valued processes $Y$ on $[0,\timeh]$ we introduce
\begin{equation}
\rho_{t}(Y)=(E[\sup_{s \in [0,t]} \|Y(\tau(s))\|_{\multsp}^{2}])^{1/2},
 \label{A. semi norm}
\end{equation}
for $ t \in [0, \infty[,$ where we have defined $Y(t)$ for $t > \timeh$ by $Y(t)=Y(\timeh).$
We note that $\rho_{t}(Y) \leq (E[\sup_{s \in [0,t]} \|Y(s)\|_{\multsp}^{2}])^{1/2},$
since $\tau (t) \leq t.$
\begin{lemma}\label{existence lemma}
If condition (\ref{App cond 1}) is satisfied and if $Y$ is an  $\mathcal{A}$-progressively measurable  $\multsp$-valued continuous process on $[0,\timeh],$
satisfying  $\rho_{t}(Y) < \infty,$ for all $t \geq 0,$ then the equation
\begin{equation}
 X(t)= Y(t)+\int_{0}^{t}\ltransm{t-s}X(s)(\drift{s}ds + \sum_{i \in \mathbb{I}}\vol{i}{s} d\wienerp{i}{s}),
 \label{linear SDE}
\end{equation}
$t \in [0,\timeh],$ has a unique solution $X,$ in the set of $\mathcal{A}$-progressively measurable $\multsp$-valued continuous process
satisfying:
\begin{equation}
  \int_{0}^{\timeh}(\| X(t) \|_{\multsp}+\| X(t) \drift{t} \|_{\multsp}
+\sum_{i \in \mathbb{I}}\| X(t) \vol{i}{t} \|^{2}_{\multsp})dt < \infty \; \text{a.s.}
 \label{X mild}
\end{equation}
Moreover this solution satisfies: \\
$i)$ If $\int_{0}^{\timeh}(\| \drift{t} \|_{\multspd}+\sum_{i \in \mathbb{I}}\| \vol{i}{t} \|^{2}_{\multspd})dt < \infty$
and $ Y$ is a continuous $\multspd$-valued process with $\rho_{t}(\partial Y) < \infty,$ for all $t \geq 0,$ then
$ X$ is a continuous $\multspd$-valued process. \\
$ii)$ If $(i)$ is satisfied and if $Y$ is a semi-martingale, then $X$ is a semi-martingale. \\
$iii)$ If $Y$ is $H^{s}$-valued, then $X$ is $H^{s}$-valued. \\
\end{lemma}
The next lemma establish conditions under which the solution of equation (\ref{linear SDE}) is in
$L^{p},$ $p \in [0, \infty[ \, .$ The notation $\tilde{\mathcal{E}}$ was introduced in (\ref{exp marting}).
\begin{lemma}\label{Lp norms}
Let condition (\ref{App cond 1}) be satisfied and let
$(i)$
$$ E[ \exp(p \int_{0}^{\timeh}(\| \drift{t} \|_{\multsp}+\sum_{i \in \mathbb{I}}\| \vol{i}{t} \|^{2}_{\multsp})dt)] < \infty,$$
for each $p \in [1,\infty[.$ Suppose that $Y$ in Lemma \ref{existence lemma}
satisfies $(ii)$
$$E[\sup_{t \in \mathbb{T}}\|Y(t) \|_{\multsp}^{p}] < \infty,$$ for each $p \in [1,\infty[.$
Then the unique solution $X$ of equation (\ref{linear SDE}) in Lemma \ref{existence lemma} satisfies
\begin{equation}
 E[\sup_{t \in \mathbb{T}}\|X(t) \|_{\multsp}^{p}] < \infty, \; \forall p \in [1,\infty[ \,.
 \label{Lp norms 1}
\end{equation}
Moreover if $(iii)$
$$E[(\int_{0}^{\timeh}(\| \drift{t} \|_{\multspd}+\sum_{i \in \mathbb{I}}\| \vol{i}{t} \|^{2}_{\multspd})dt)^{p}]
    < \infty$$
and $(iv)$
$$E[\sup_{t \in \mathbb{T}}\| Y(t) \|_{\multspd}^{p}] < \infty,$$
 for each $p \in [1,\infty[,$ then also
\begin{equation}
 E[\sup_{t \in \mathbb{T}}\| X(t) \|_{\multspd}^{p}] < \infty,  \; \forall p \in [1,\infty[ \,.
 \label{Lp norms 3}
\end{equation}
In particular, estimates (\ref{Lp norms 1}) and (\ref{Lp norms 3}) applies to
$X=\tilde{\mathcal{E}}(L).$
\end{lemma}

\noindent \textbf{Note added in the proofs:}
Since the preparation of this paper, the optimal bond portfolio problem has further been
studied in various directions:
\begin{itemize}
\item[1] The reference De Donno, M. and Pratelli, M.:
 \textit{A theory of stochastic integration for bond markets},
Ann. Appl. Probab. \textbf{15}, 2773--2791 (2005)
considers the optimal bond portfolio problem in a more general semi martingale bond market.
Existence of optimal wealth strategies is established and existence of optimal portfolios is studied.
\item[2] The reference Ringer, N. and Tehranchi, M.:
 \textit{Optimal portfolio choice in the bond market},
Finance Stoch. \textbf{10}, 553--573 (2006) considers the optimal bond portfolio problem
in a Markovien setting of local volatility operators with full range and
which are globally Lipschitzien. More precisely it is assumed,
with our notations and limiting us to the time homogeneous case,
that the function $C: E \rightarrow \mathcal{HS}(\ell^{2},E),$ where
$C(f) = f\vol{}{}(f),$ is globally Lipschitzien and that for all strictly positive $f \in E$
the closure of the range $\mathcal{R}(C(f))$ is the subset of elements $g \in E$ such that
$g(0)=0$. If moreover (the unique) market price of risk is
globally Lipschitzien then they establish the existence of a solution to
the optimal portfolio problem.
We note that the proof of this boils down to the verification of properties of
the Malliavin derivative of $\ln(\xi_{\timeh})$ as was already the case in
Theorem 4.5 of \cite{E.T Bond  Completeness} (see Theorem \ref{th opt port l^2}).
We also note that their Gaussian example, of course satisfies our condition (C) of Definition \ref{cond. C},
so it is covered by our treatment.
\end{itemize}


\begin{thebibliography}{9}

\bibitem{Adams 03} Adams, R.A. and Fournier, J.J.F.: \emph{Sobolev Spaces},
Academic Press 2003.

\bibitem{Bj-Ka-Ru97} Bj\"ork, T., Kabanov, Y. and Runggaldier, W.: \emph{Bond
market structure in the presence of marked point processes}, Mathematical
Finance, \textbf{7}, 211--239 (1997).

\bibitem{Bj-Ma-Ka-Ru97} Bj\"ork, T., Masi, G., Kabanov, Y. and Runggaldier,
W.: \emph{Toward a general theory of bond markets}, Finance and Stochastics,
\textbf{1}, 141--174 (1997).

\bibitem{Bj-Sv01}
Bj\"ork, T. and Svensson, L.: {On the Existence of Finite Dimensional Realizations
for Nonlinear Forward Rate Models}, Mathematical Finance, {\bf 11}, 205--243 (2001).

\bibitem{Calderon}
Calderon, A.P.: {\em Lebesgue spaces of differentiable functions and distributions},
Proc. Symp. Pure Math. IV, AMS 1961, 33--49.

\bibitem{Carmona-Tehr}
Carmona, R. and Tehranchi, M.: {\em A Characterization of Hedging Portfolios for
Interest Rate Contingent Claims}, Preprint March 24, 2003.

\bibitem{DaPrato-Zabczyk} Da Prato, G. and Zabczyk, J.: \emph{Stochastic
Equations in Infinite Dimensions}, Encyclopedia of Mathematics and its
Applications, Cambridge University Press, 1992.

\bibitem{DeDonno Pratelli}
De Donno, M. and Pratelli, M.: {\em On the use of measure-valued strategies in bond markets},
Finance and Stochastics, {\bf 8}, 87--109 (2004).

\bibitem{I.E.-R.T} Ekeland, I. and T\'emam, R.: {\em Convex Analysis and Variational
Problems}, Classics in Applied Mathematics 28, SIAM 1999.

\bibitem{I.E.-E.T bond th} Ekeland, I. and Taflin, E.: \emph{A Theory of Bond
Portfolios},  Ann. Appl. Probab. \textbf{15}, 1260--1305 (2005). Also
%
\texttt{http://arxiv.org/abs/math.OC/0301278}

\bibitem{Filipovic} Filipovi\'c, D.: {\em Consistency Problems for HJM Interest Rate Models},
Phd thesis, Dep. Math. ETH, Z\"urch 2000
Preprint 2001.

\bibitem{HJM92} Heath, D.C., Jarrow, R.A. and Morton, A.: \emph{Bond pricing
and the term structure of interest rates: a new methodology for contingent
claim valuation}, Econometrica, \textbf{60}, 77--105 (1992).

\bibitem{Horm}  H\"ormander, L.: {\em The analysis of  linear  partial  differential
operators}, Vol. I, Springer-Verlag 1985.

\bibitem{Kall-Xiong}  Kallianpur, G., and J. Xiong,
Stochastic Differential Equations in Infinite Dimensional Spaces, Lecture
Notes-Monograph Series, Institute of Mathematical Statistics, 1995.

\bibitem{K-S 99}
Karatzas, I. and Shreve, S.E.: {\em Methods of Mathematical Finance}, Applications of Mathematics, Volume 9,
Springer-Verlag 1999.

\bibitem{Kato66} Kato, T. \emph{Perturbation Theory for Linear Operators}, Die
Grundleheren der mathematischen Wissenschaften, Volume 132, Springer-Verlag,
New York 1966.

\bibitem{Kr-Scha} Kramkov, D. and Schachermayer W.: \emph{The Asymptotic
Elasticity of Utility Functions and Optimal Investment in Incomplete Markets},
Annals Appl. Probability, \textbf{9}, 904--950 (1999).

\bibitem{Lax 02} Lax, P.D.: \emph{Functional Analysis},
Wiley-Interscience 2002.

\bibitem{Lintner65} Lintner, J.: \emph{The Valuation of Risk Assets and the
Selection of Risky Investments in Stock Portfolios and Capital Budgets},
The Review of Economics and Statistics, {\bf 47}, 13--37 (1965).

\bibitem{Markowitz52}
Markowitz, H.: {\em Portfolio Selection},
Jour. Finance, {\bf 7}, 77--91 (1952).

\bibitem{Mert69}
Merton, R.: {\em Lifetime Portfolio Selection Under Uncertainty: The
Continuous-Time case},
Rev. Economics and Stat. {\bf 51}, 247--257 (1969).

\bibitem{Mert71}
Merton, R.: {\em Optimum Consumption and Portfolio Rules in a Continuous Time
Model},
Jour. Economic Theory, {\bf 3}, 373--413 (1971).

\bibitem{Mikul-Rozo98}
Mikulevicius, R. and Rozovskii, B.L.:
{\em Normalized stochastic integrals in topological vector spaces},
Seminaire de Probabilites XXXII, LNM, Springer-Verlag, 1998

\bibitem{Mikul-Rozo99}
Mikulevicius, R. and Rozovskii, B.L.:
{\em Martingale problems for SPDE's},
Stochastic Partial Differential Equations: Six Perspectives, Ed:
R. Carmona and BL Rozovskii, AMS, Mathematical Surveys and Monographs, 1999

\bibitem{Musi93}
Musiela, M., {\em Stochastic PDEs and term structure models}, Journ\'ees Internationales de Finance, IGR-AFFI, La Baule, 1993.

\bibitem{Nual71}
Nualart D.: {\em The Malliavin Calculus and Related Topics}, Probability and its Applications,
Springer-Verlag, 1991.

\bibitem{Pham 2003}
Pham, H.: {\em A predictable decomposition in infinite asset model with jumps.
Application to hedging and optimal investment}, Stochastics and Stochastic Reports, {\bf 5}, 343--368 (2003).

\bibitem{Pliska86}
Pliska, S.R.: {\em A stochastic calculus model of continuous trading: optimal portfolios},
Math. Operations Research {\bf 11}, 371--382 (1986)

\bibitem{Revuz-Yor} Revuz, D. and Yor, M.: \emph{Continuous Martingales and
Brownian Motion}, Grundlehren der mathematischen Wissenschaften, Band 293, Spriner-Verlag

\bibitem{Rudin 1}
Rudin, W.: {\em Real and Complex Analysis}, 3rd edition,  McGraw-Hill, 1986.

\bibitem{Rudin 2}
Rudin, W.: {\em Functional Analysis}, 2nd edition, McGraw-Hill, 1991.

\bibitem{Rutkowski99}
Rutkowski, R.: {\em Self-financing Trading Strategies for Sliding, Rolling-horizon,
and Consol Bonds}, Math. Finance {\bf 5}, 361--385 (1999)

\bibitem{Sharp64}
Sharp, W.F.: {\em Capital Asset Prices: A Theory of Market Equilibrium under Conditions of Risk},
The Journal of Finance, {\bf 19}, 425--442 (1964).

\bibitem{E.T Bond  Completeness} Taflin, E.: \emph{Bond Market Completeness
and Attainable Contingent Claims},  Fin. Stoch.   {\bf9}, 429--452 (2005). Preprint \\
\texttt{http://arxiv.org/abs/math.OC/0402364}

\bibitem{Yosida}
Yosida, K.: {\em Functional Analysis}, Grundlehren der mathematischen Wissenschaften,
Band 123, Springer-Verlag.

\end{thebibliography}
\end{document}